\newif\ifmydraft
\newif\iflinenumbers
\linenumberstrue
\documentclass[]{tran-l}   % from template
\usepackage{paralist}
\usepackage{wrapfig}
\numberwithin{equation}{section} % from template
\iflinenumbers
 \usepackage[pagewise,modulo]{lineno}
\else
  \newcommand{\linenumbers}{}
\fi
\def\getRCS$#1 #2 #3 #4 #5${%
\def\RCSfile{#2}\def\RCSver{#3}\def\RCSdate{#4}}
\getRCS$Id: iomega2.tex,v 2.19 2007/04/30 18:08:01 mitchell Exp $
\ifmydraft                        % \iffalse for distribution!!!!
 \usepackage[active]{srcltx}   
\fi
\usepackage{amssymb}

\usepackage{verbatim}           % only for comments
\newenvironment{inparaenumi}{\begin{inparaenum}[\upshape(i)]}
  {\end{inparaenum}}
\newtheorem{theorem}{Theorem}[section]
\newtheorem{proposition}[theorem]{Proposition}
\newtheorem{lemma}[theorem]{Lemma}
\newtheorem{corollary}[theorem]{Corollary}

\newtheorem*{CLAIM}{Claim}

\theoremstyle{definition}
\newtheorem{definition}[theorem]{Definition}

\theoremstyle{remark}

\newtheorem*{remark*}{Remark}

\theoremstyle{plain}

\newcommand\CASE[1]{\smallskip\par\noindent(\emph{Case #1})\;}

\newcommand\note[2][]{\ifmydraft[[#1 --- #2]]\fi}

\newcommand\deq{\mathrel{:=}}
\newcommand{\model}{model}
\newcommand{\simple}{simple}

\newcommand{\ga}{\alpha}

\renewcommand\gg{\gamma}
\newcommand{\gd}{\delta}

\newcommand{\gk}{\kappa}
\newcommand{\gl}{\lambda}

\newcommand{\gth}{\theta}
\newcommand{\gw}{\omega}

\newcommand{\gS}{\Sigma}
\newcommand\ca{\mathcal A}

\newcommand{\ns}{\text{NS}}
\newcommand{\smc}{strongly generic condition}
\newcommand{\lex}{\le^*}
\newcommand{\eqx}{=^*}
\newcommand\inx{\in^{*}}          % For ``effectively in because of
\newcommand\ninx{\notin^{*}}

\newcommand\1{\mathbf{1}}
\newcommand\ps{\mathcal{P}}
\newcommand\fequiv{\equiv}
\newcommand\forces{\Vdash}
\newcommand\nforces{\nVdash}

\newcommand\decides{\parallel}
\newcommand\seq[1]{\langle\,#1\,\rangle}

\newcommand\set[1]{\{\,#1\,\}}
\newcommand\sing[1]{\{#1\}}
\newcommand\card[1]{\left|#1\right|}
\newcommand\restrict{{\upharpoonright}}

\newcommand\image{\text{``}}

\DeclareMathOperator{\dinter}{\bigtriangleup}
\DeclareMathOperator{\cof}{cf}
\DeclareMathOperator{\Cof}{Cof}
\DeclareMathOperator{\domain}{domain}

\DeclareMathOperator{\otp}{otp}

\newcommand\nothing{\varnothing}
\renewcommand{\emptyset}{\varnothing}

\newcommand\union{\bigcup}

\newcommand\II{I}
\newcommand\OO{O}
\newcommand{\OOa}[2][\gamma]{\OO_{#1,(#2',#2]}}
\newcommand{\OOO}[1][\eta]{\OOa{#1}}
\newcommand\CC{C}
\newcommand\DD{D}               % The generic sequence of closed, unbounded sets.
\newcommand{\RR}{R}               % ``complete'' requirements.

\newcommand\cut{{|}}
\newcommand\icut[2]{#1|^*#2}
\newcommand\Pkp{P^*}            % Forcing of section 3.

\DeclareMathOperator{\cpx}{cp}
\newcommand{\cp}[2]{\cpx^{#1}(#2)}

\usepackage[%
 hyperref%  For xdvi or dvips with -z.
]{hyperref}
\begin{document}
\title{$I[\gw_2]$ can be the nonstationary ideal on $\Cof(\gw_1)$}
\author{William J.  Mitchell}
\address{Department of Mathematics\\
  PO Box 118105\\
  University of Florida\\
  Gainesville, Florida 32611-8105\\USA}
\email{mitchell@math.ufl.edu}
\urladdr{http://www.math.ufl.edu/~mitchell}
\thanks{I would like to thank Matt
 Foreman, Bernard Koenig and the referee of this paper for valuable remarks and
 corrections, .  In 
 addition
 I would like to thank Matt Foreman and Martin Zeman for inviting me
 to a weeks visit to the University of California at Irvine for a week
 during which I gave an extended
 exposition 
of this work.   Suggestions made during this visit lead
 directly to a dramatically improved revision of this paper.
 \endgraf The writing of this paper was partially supported
%%% TO-DO: Fill in grant number.
 by grant number DMS-0400954 from the National Science Foundation.  
}

%\tableofcontents
\maketitle
\linenumbers
\begin{abstract}
  We answer a question of Shelah by showing that it is consistent that
  every member of 
  $I[\gw_2]\cap\Cof(\gw_1)$ is nonstationary if and only if it is
  consistent  that there is a $\gk^+$-Mahlo cardinal $\gk$.
\end{abstract}
%%% === TO-DO DELETE THIS ===
\tableofcontents

\section{Introduction}
In \cite[definition 2.1]{shelah:rsss}
Shelah  defined the following ideal $I[\gk^+]$:
\begin{definition}
  \label{def:Ikplus}
  Define, for any sequence 
  $A=\seq{a_\ga:\ga<\gk^+}$, the set $B(A)$  to be the set of ordinals
  $\nu<\gk^+$ such that there is a set $c\subset\nu$ with
  $\otp(c)=\cof(\nu)$, $\union c=\nu$, and
  $\set{c\cap\xi:\xi<\nu}\subset\set{a_\ga:\ga<\nu}$.   Then $I[\gk^+]$
  is the set of subsets of $\gk^+$ which are contained, up to a
  nonstationary set, in some set $B(A)$. 
\end{definition}

He proved in \cite[theorem~4.4]{shelah:rsss} that
$\gk^+\cap\set{\nu:\cof(\nu)<\gk}\in I[\gk^+]$ for all regular
cardinals $\gk$, 
and observed that it is consistent that the restriction of
$I[\kappa^+]$ to $\set{\nu<\kappa^+:\cof(\nu)=\kappa}$ is generated by
a single stationary, costationary set.   He 
asked whether it is consistent that every subset of
$\set{\nu<\kappa^+:\cof(\nu)=\kappa}$ in 
$I[\kappa^+]$ is nonstationary.   We answer this in the affirmative for
$\kappa=\omega_1$. 
\begin{theorem}
  \label{thm:mainthm}
  If it is consistent that there is a cardinal $\gk$ which is
  $\gk^+$-Mahlo then it is consistent that 
  $I[\gw_2]$ does not contain any stationary subset of
  $\set{\nu<\omega_2:\cof(\nu)=\omega_1}$. 
\end{theorem}

The fact that a $\gk^+$-Mahlo cardinal $\gk$ is necessary is due to
Shelah;  a proof is given in \cite[theorem 13]{mitchell.wvsi}.

Our proof of theorem~\ref{thm:mainthm} uses forcing  to 
add a sequence $\seq{D_\alpha:\alpha<\kappa^+}$ of closed
unbounded subsets of $\kappa$, in the process collapsing the cardinals
between $\omega_1$ and $\kappa$ onto $\omega_1$ so that $\kappa$ becomes
$\omega_2$.
In the resulting model there is, for every set of the form $B(A)$,  some
ordinal $\alpha<\kappa^+$ such that $B(A)\cap
D_{\alpha}$ does not contain any ordinal of cofinality $\omega_1$.  
Thus every set in 
$I[\omega_2]\restrict\set{\nu<\omega_2:\cof(\nu)=\omega_1}$ is
nonstationary. 

Section~\ref{sec:oneclub} introduces the basic construction by showing
how to add a single new closed unbounded set.   This serves as a
warm-up for section~\ref{sec:moreclubs},  introducing ideas of the
forcing in a simpler context, and 
also is used in section~\ref{sec:moreclubs} as the basic component of
the forcing used to prove theorem~\ref{thm:mainthm}.

\medskip

Most of our notation is standard.  
We write $\lim(X)$ for the set of limit ordinals $\ga$ such that $\ga\cap X$
is cofinal in $\ga$, and $\overline X$ for the topological closure,
$X\cup\lim(X)$, of $X$, and 
we write $\Cof(\lambda)$ for $\set{\nu:\cof(\nu)=\lambda}$.
%% We use the nonstandard notation $\sup^X(A)$
%% to mean $\min(X\setminus\sup(A\cap X))$.

The reader of this paper may find it helpful to also consult the
expository paper 
\cite{mitchell.acus}, which discusses some  of the material covered in
this paper along with related topics.

A basic ingredient of the forcing in this paper is the idea of forcing
with models as side 
conditions.    This idea,  in the form used in this paper, was
discovered independently 
by the author but the general
technique method was originally introduced 
and has been extensively
investigated by Todorcevic.
His original applications concerned properties of 
$\omega_1$ and used forcing notions which collapsed $\omega_2$, but in
later
applications such as \cite{todorcevic.dsct} he used a form, related to
that used in this paper, which did not collapse $\omega_2$.   Koszmider
\cite{koszmider:scuf}
has  
developed a 
modification of Todorcevic's technique which uses a previously given
morass
to simplify the actual forcing.   Koszmider's method is arguably
simpler, but it is not suitable for the present construction:
it
would require a morass on $\omega_2$ of the generic extension, which
is the inaccessible cardinal $\kappa$ of the ground model. 

A forcing essentially identical to that described in
section~\ref{sec:oneclub} was discovered independently by Sy
Friedman~\cite{friedman.ffc}.    
The presentation in \cite{friedman.ffc} does not collapse any
cardinals, instead adding a closed unbounded subset
of the $\omega_2$ of the ground model; however this difference is due
to the difference in the intended applications of the forcing rather
than any intrinsic difference in the forcing itself.

\section{Adding a single closed unbounded set}
\label{sec:oneclub}
\subsection{The forcing notion}

In this section we define a new forcing $P_B$ which adds a closed
unbounded subset $D$ of the set $B^*\deq B\cup\Cof(\gw)$, where $B$ is
a stationary subset of $\set{\lambda<\kappa:\cof(\lambda)>\omega}$ for
a regular cardinal $\kappa$.  The forcing preserves
$\gw_1$, while (if
$\kappa>\omega_2$) collapsing the intervening cardinals to make $\gk=\gw_2$.
This forcing  
serves both as a warm-up for and as the basic building block of the
forcing in section~\ref{sec:moreclubs} which adds $\gk^+$ many closed
unbounded 
sets to construct a model in which  $I[\gw_2]\cap\Cof(\gw_1)$  is the
nonstationary ideal. 
As another application, 
it will be observed later in this section that this forcing gives a
new construction of a model with  no special $\aleph_2$-Aronszajn trees
(or, starting from a weakly compact cardinal, no $\aleph_2$-Aronszajn
trees), and it is shown in \cite{mitchell.acus} that stripping this
forcing down to its basic  technique yields a construction of such a
model which is  much simpler
than any of those which were previously known.
\medskip{}

The forcing $P_B$ is based on the standard finite forcing
$P_{\omega_1}$, introduced by Baumgartner
in \cite[page~926]{baumgartner.apfa}, for adding a closed unbounded
subset of $\omega_1$.  In order to motivate the definition of $P_B$ we
give a brief description of this forcing $P_{\omega_1}$, show how a
straightforward attempt to apply it to $\omega_2$ fails, and describe
the new technique which we use to make it succeed.

The presentation of $P_{\omega_1}$ which we will give is a variant of
a version,  discovered by U. Abraham  
\cite{abraham-shelah.fcus},  of Baumgartner's forcing. 
The set $D$ constructed by this forcing is not generic for
Baumgartner's forcing as originally described in 
\cite{baumgartner.apfa}, since $D$ has the 
property that $\liminf_{\ga<\nu}\otp(D\cap(\nu\setminus\ga))$ is as 
large as possible for any limit ordinal  $\nu\in D$; 
however Zapletal  \cite{zapletal.ccf} has shown 
that the two forcings are equivalent. 

The conditions in the forcing $P_{\omega_1}$ are finite sets of
symbols which we call \emph{requirements}.   There are two types of
these requirements: $\II_{\lambda}$ for ordinals $\lambda<\omega_1$,
and $\OO_{(\eta',\eta]}$ for pairs of ordinals 
$\eta'<\eta<\omega_1$.   
Two requirements $\II_{\lambda}$ and $\OO_{(\eta',\eta]}$ are
\emph{incompatible} if $\eta'<\lambda\le\eta$; any other two
requirements are \emph{compatible}.   A condition in $P_{\omega_1}$ is a
finite set of 
requirements, any two of which are compatible, and the ordering of
$P_{\omega_1}$ is by superset: $p'\le p$ if $p'\supseteq p$.

If  $G$ is a generic subset of $P_{\omega_1}$ then we define
$D\deq\set{\lambda<\omega_1:\II_{\lambda}\in\bigcup G}$.   A little 
thought shows that 
\begin{equation}\label{eq:3}
\forall\lambda<\omega_1\;\left(\lambda\notin
  D\iff
  \exists\eta',\eta\;\bigl(\OO_{(\eta',\eta]}\in\bigcup G\And
  \eta'<\lambda\le\eta\bigr) 
%%   \exists\OO_{(\eta',\eta]}\in\bigcup
%%   G\bigr)\;(\eta'<\lambda\le\eta)
\right),
\end{equation}
and it follows that  $D$ is a closed and unbounded subset of $\omega_1$.

The cardinal $\omega_1$ is preserved by the forcing $P_{\omega_1}$
because the forcing is 
proper; indeed it has the stronger property that if $M$ is any
countable elementary substructure of $H_{\omega_1}$ and
$\lambda=\sup(M)$ then the condition 
$\sing{\II_{\lambda}}$ not only forces that $G\cap M$ is $M$-generic,
but actually forces that $G\cap M$ is a $V$-generic subset of
$P_{\omega_1}\cap M$.
Note that for this property it is sufficient to take $M\prec
H_{\omega_{1}}$ rather than $H_{\omega_2}$: since the relevant dense
sets are taken from $V$, rather than 
from $M$, it is not necessary that $P\in M$. 

In order to define a similar forcing $P_{\omega_2}$
adding a 
new closed unbounded subset of $\omega_2$, one could naively try to use the 
same  definition, but with requirements $\II_{\lambda}$
for any $\lambda<\omega_2$ and $\OO_{(\eta',\eta]}$ for
any $\eta'<\eta<\omega_2$; however this forcing is not proper and does
collapse $\omega_1$.
To simplify notation we will show why this is true below the condition
$\sing{\II_{\gw_1\cdot\gw}}$, which  forces
$\gw_1\cdot\gw\in D$.   For each $n<\gw$ let
$\xi_n=\sup\set{\xi<\gw_1:\gw_1\cdot n+\xi\in D}$, so that
$0\le\xi_n\le\gw_1$.
If $p\le\sing{\II_{\gw_1\cdot\gw}}$ and $\xi<\omega_1$ then for any
sufficiently large $n<\omega$ the set 
$p\cup\sing{\II_{\gw_1\cdot n+\xi},\OO_{(\gw_1\cdot
    n+\xi,\gw_1\cdot(n+1)]}}$
is a condition extending $p$ which forces that $\xi_n=\xi$.
It follows that
$\omega_1=\set{\xi_n:n<\omega\And\xi_n<\omega_1}$; thus $\omega_1$ is
collapsed in $V[D]$.
\medskip{}

In order to avoid this problem we will use a third type of
requirement in the definition of $P_{\omega_2}$.
This new requirement, which we write as  $C_{M}$ for any
countable $M\prec H_{\omega_2}$, 
is intended to play the same role as the requirement $I_{\lambda}$
plays in the proof that $P_{\omega_1}$ is proper:
the 
condition $\sing{\CC_{M}}$ will force that $G\cap M$ is a
$V$-generic subset of $P_{\omega_2}\cap M$.
This will be accomplished by finding, for each condition
$p\le\sing{\CC_M}$,
a condition $p\cut M\in P_{\omega_2}\cap M$ with the property that every
condition 
$q\le p\cut M$ in 
$P_{\omega_2}\cap M$ is compatible
with $p$: thus the condition $p\cut M\in M$ will
capture all of the 
influence which $p$ has on the forcing $P_{\omega_2}\cap M$. 

To see how this works, consider a set
$p=\sing{\CC_M,\OO_{(\eta',\eta]}}$.    If $(\eta',\eta]\cap M=\nothing$
then $\OO_{(\eta',\eta]}$ will have no effect on the forcing inside $M$,
and we will take $p\cut M=\nothing$.    If $\eta'$ and $\eta$ are in
$M$ then $(\eta',\eta]$ is a member of $M$, and we will take $p\cut
M=\sing{\OO_{(\eta',\eta]}}$.  In either case $p$ will be a condition,
but if neither of these holds---if 
$(\eta',\eta]$ intersects $M$ but is not a member of $M$---then there
is no requirement inside $M$ which will have the same effect on $G\cap
M$ as the requirement
$\OO_{(\eta',\eta]}$ does, and in this case we will say that
$\OO_{(\eta',\eta]}$ and 
$\CC_{M}$ are not compatible, and hence $p$ is not a condition.

To see how this will block the collapse of $\omega_1$ described
earlier for the
naive version of the forcing, let $M$ be any countable elementary
substructure of $H_{\omega_2}$ with $\omega_1\cdot\omega\in M$.
The pair $\sing{\CC_M,\II_{\omega_1\cdot\omega}}$ will be a a
condition, and as with $P_{\omega_1}$ it will force
$\omega_1\cdot\omega\in D$.
Now suppose that $p\le\sing{\CC_M,\II_{\omega_1\cdot\omega}}$ is a condition
which forces, for some $n<\omega$, that $\xi_n<\omega_1$, that is,
that $D\cap\omega_1\cdot(n+1)$ is bounded in $\omega_1\cdot(n+1)$.   
By using \eqref{eq:3} (which we will show to hold for $P_{\omega_2}$) we can see
that this implies that there is a requirement $\OO_{(\eta',\eta]}\in
p$ with $\eta'<\omega_1\cdot(n+1)\le\eta$. Now $(\eta',\eta]\cap
M\not=\nothing$, since $\omega_1\cdot(n+1)$ is in the intersection, so the
compatibility of $\CC_M$ with $\OO_{(\eta',\eta]}$ implies that
$\OO_{(\eta',\eta]}\in M$ and in particular that
$\eta'<\sup(M\cap\omega_1\cdot(n+1))$.   Hence
$p\forces\dot\xi_n\le\sup(M\cap\omega_1)$, and since $n$ was arbitrary
it follows that
$\sing{\CC_M,\II_{\omega_1\cdot\omega}}\forces
\set{\dot\xi_n:n<\omega\And\dot\xi_n<\omega_1}=M\cap\omega_1$. 

\bigskip{}

We are now ready to give the definition of the forcing $P_{B}$.   We
assume that $B$ is a stationary subset of an inaccessible cardinal
$\kappa$ and that  every member $\lambda$ of $B$ is a cardinal with
uncountable cofinality such that  $H_{\lambda}\prec H_{\kappa}$.   
This definition can be easily adapted (assuming that
$2^{\omega}=\omega_1$) to the case $\kappa=\omega_2$, discussed
previously as $P_{\omega_2}$, by replacing the models
$H_{\lambda}$ in the definition with structures
$L_{\lambda}[A]$, where $A\subset\omega_2$ enumerates
$[\omega_2]^{\omega}$.   Friedman \cite{friedman.ffc} has pointed 
out that the assumption $2^{\omega}=\omega_1$ can be weakened, provided
that there exists a stationary set 
$S\subset[\omega_2]^{\omega}$ such that $\card{\set{x\cap\nu:x\in
    S}}=\omega_1$ for all $\nu<\omega_2$.

We also assume that $H_{\kappa}$ has definable Skolem functions, so
that $M\cap N\prec H_{\kappa}$ whenever $M\prec H_{\kappa}$ and
$N\prec H_{\kappa}$.    This assumption can be avoided by replacing
$H_{\kappa}$ with a structure which does have Skolem functions. 

%%%% NOTE: In the submitted paper I required $H_\lambda\prec H_\kappa$
%%%% for $\lambda\in B^*$ of cofinality $\omega$.   I can't remember
%%%% why I did this, but I've been all through section2 and confirmed
%%%% that it wa not used there, at least.  Thus this definition should
%%%% be okay.   It remains to see whether I can make the same change
%%%% in section 3.

We write $B^*=B\cup\set{\lambda<\kappa:\cof(\lambda)=\omega}$.
% We write $B^*$ for $B\cup\set{\lambda<\kappa:\cof(\lambda)=\omega\And
% \card{\lambda}=\lambda\And
% H_{\lambda}\prec H_{\kappa}}$.
The forcing will add a new closed
unbounded subset of $B^*$.

The forcing $P_B$ uses three types of \emph{requirements}:
\begin{enumerate}
\item $\II_{\lambda}$, for any $\lambda\in B^*$,
\item $\OO_{(\eta,\eta']}$, for any interval with $\eta<\eta'<\kappa$, and 
\item $\CC_{M}$, for any countable set $M\prec H_{\kappa}$.
\end{enumerate}

These symbols $\II_{\lambda}$, $\OO_{(\eta',\eta]}$ and $\CC_{M}$
are used for convenience; since the subscripts are distinct we
can take each requirement to be equal to its subscript, that is,
$\II_{\lambda}=\lambda$, $\OO_{(\eta',\eta]}=(\eta',\eta]$ and
$\CC_M=M$. 

We first specify which pairs of requirements are compatible.   
The first clause is the same as for $P_{\omega_1}$, and 
an explanation of the second clause has already been given.
Clauses~3 and~4 similarly assert that $\CC_{M}$ is compatible with
$\II_{\lambda}$ or $\CC_N$ if and only if there there is a condition
$\II_{\lambda}\cut M$ or $\CC_N\cut M$ which is a member of $M$ and
reflects the effect which the requirement $\II_{\lambda}$ or $\CC_{N}$
in $P_B$ has on  the forcing $P_B\cap M$.
This will be made precise in lemma~\ref{thm:compat}.

\begin{definition}\label{def:Pcompat}
  \begin{enumerate}
  \item \label{item:OI}
    Two requirements $\OO_{(\eta',\eta]}$ and $\II_{\lambda}$ are
    \emph{incompatible} if $\eta'<\lambda\le\eta$; otherwise they are
    \emph{compatible}. 
  \item\label{item:OM}
    Two requirements $\OO_{(\eta',\eta]}$ and $\CC_{M}$ are
    \emph{compatible} if either $\OO_{(\eta',\eta]}\in M$ or every
    requirement $\II_{\lambda}\in M$ is compatible with
    $\OO_{(\eta',\eta]}$. 
  \item\label{item:MI}
    \begin{enumerate}
    \item \label{item:MIf}
      An \emph{$M$-fence for a requirement $\II_{\lambda}$} is a requirement
      $\II_{\lambda'}\in M$ such that any requirement
      $\OO_{(\eta',\eta]}$ in $M$ incompatible with $\II_{\lambda}$ is
      also incompatible with $\II_{\lambda'}$.
    \item\label{item:MIc}
      Two requirements $\CC_M$ and $\II_{\lambda}$ are \emph{compatible} if
      either $\lambda\ge\sup(M)$ or there exists a $M$-fence for $\II_{\lambda}$.
    \end{enumerate}
  \item\label{item:MN}
    \begin{enumerate}
      \item\label{item:MNf}
        An \emph{$M$-fence for a requirement $\CC_{N}$} is a finite
        set $x\in M$ of requirements $I_\lambda$, with $\lambda\in B$,
        with the following property: Let 
        $\OO_{(\eta',\eta]}\in M$ be any requirement which is
        incompatible with $N$, 
        and which has $\eta'\geq\sup(M\cap N)$ if $M\cap N\in M$.
        Then 
        there is some $\II_{\lambda}\in x$ which is incompatible with
        $\OO_{(\eta',\eta]}$. 
      \item\label{item:MNc}
        Two requirements $\CC_M$ and $\CC_N$ are \emph{compatible} if
        the following clauses hold both as stated and with $M$ and
        $N$ switched: 
        \begin{enumerate}
        \item\label{item:MNci}
          Either $M\cap N\in M$ or $M\cap N=M\cap H_{\sup(M\cap N)}$.
        \item\label{item:MNcii}
          There is a $M$-fence for $\CC_N$.
        \end{enumerate}
    \end{enumerate}
  \end{enumerate}    
\end{definition}

\begin{definition}
  \label{def:PB}
  A condition $p$ in the forcing $P_B$ is a finite set of requirements 
  such that each pair of requirements in $p$ is compatible.   The set
  $P_B$ is ordered by reverse inclusion: $p'\le p$ if $p'\supseteq p$.
\end{definition}

\begin{proposition}
  If $\CC_M$ is a requirement and $p\in M\cap P_B$ then $p\cup
  \sing{\CC_M}$ is a condition.
\end{proposition}
\begin{proof}
  In verifying that $\CC_M$ is compatible with any requirement in $p$,
  notice that any requirement $\II_{\lambda}\in M$ is its
  own $M$-fence.

  For each requirement $C_N\in p$, the model $N$ is a member of $M$ and hence
  $\nothing$ is a $M$-fence for   $C_N$. 
\end{proof}

Notice that Clauses~(\ref{item:MNf}) and (\ref{item:MNci}) of 
Definition~\ref{def:Pcompat} imply
that $\sup(M\cap N)\notin M$ unless $M\cap N\in M$.   Suppose to the
contrary that $M\cap N\notin M$ but
$\lambda\deq\sup(M\cap N)\in M$.
Then  for any
$\eta'\in M\cap\lambda$ the  requirement
$\OO_{(\eta',\lambda]}$ is in $M$ and is incompatible with $N$.   The only way that a fence $x\in M$ could be
incompatible with all such requirements $\OO_{(\eta',\lambda]}$ would
be if $\II_{\lambda}\in x$, but $\lambda\notin B$ since
$\cof(\lambda)=\omega$.

In the case $M\cap N\in M$, the requirement $\CC_{M\cap N}$ will be
used in section~\ref{sec:smc} 
to augment the $M$-fence for $\CC_N$:
Any requirement
$\OO_{(\eta',\eta]}\in M$ with $\eta'<\sup(M\cap N)$ which is
incompatible with $\CC_N$ will be incompatible with $\CC_{M\cap N}\in M$.

Definition~\ref{def:Pcompat} described $M$-fences 
in terms of their function.   We now give alternate structural
characterizations and note that the fences are unique:   
\begin{proposition}\label{thm:MIcompat}
  The requirements $\II_{\lambda}$ and $\CC_M$ are compatible if either
  $\lambda\ge\sup(M)$ or else $\min(M\setminus\lambda)\in B^*$; in the
  later case $\II_{\min(M\setminus\lambda)}$ is the unique $M$-fence
  for $\II_{\lambda}$. 
\end{proposition}
\begin{proof}
  Set $\lambda'=\min(M\setminus\lambda)$.   
  If $\lambda'\in B^*$ then $\II_{\lambda'}$ is a
  requirement, and it is easy to see that it is a $M$-fence for
  $\II_{\lambda}$.
  
  To see that it is the only possible $M$-fence for
  $\II_{\lambda}$, note that if $\eta\in M\cap\lambda'$ then 
  the requirement
  $\OO_{(\eta,\lambda']}$ is a member of $M$ 
  and is incompatible with $\II_{\lambda}$.
  However any requirement $\II_{\lambda''}\in M$ with
  $\lambda''\not=\lambda'$ will be compatible with
  $\OO_{(\eta,\lambda']}$, provided that  $\eta>\lambda''$ in the case
  that $\lambda''<\lambda'$.
\end{proof}
The structural characterization of an $M$-fence for $\CC_N$ is slightly
 more complicated:

\begin{proposition}\label{thm:MNcompat}
  Suppose that $\CC_{M}$ and $\CC_{N}$ are requirements, and let $y$
  be the set of ordinals $\lambda\in M$ such that 
  $\lambda>\sup(M\cap N)$ and $\lambda=\min(M\setminus\eta)$ for some
  $\eta\in N$.

  Then there is $M$-fence for $\CC_{N}$ if and only if $y$ is finite,
  $y\subset B$, 
  and if   $M\cap N\notin M$ and $M\nsubseteq N$ then 
  $\min(M\setminus\sup(M\cap N))\in y$.
  In this case
  $x\deq\set{\II_{\lambda}:\lambda\in y}$ is a $M$-fence
  for $\CC_N$, and $x$ is minimal in the sense that it is a subset of any
  other 
  $M$-fence for $\CC_{N}$. 
\end{proposition}
  
 Proposition~\ref{thm:MNcompat} asserts
that two compatible requirements 
$\CC_M$ and $\CC_{N}$ divide each other into finitely
many blocks:  a common block below $\sup(M\cap N)$, followed by 
finitely many disjoint blocks alternating between $M$ 
and $N$.    Each block lies inside a gap
in the other model, the upper end of which is delineated by a member
of the $M$-fence for $\CC_N$ 
or the $N$-fence for $\CC_M$. 

\begin{wrapfigure}{r}{3cm}
    \setlength\unitlength{.8cm}
  \newcommand\dotsize{.2}
  \linethickness{1.5pt}
  \hskip 3pt
    \begin{picture}(3.5,3.5)(0,0)
      \put(.8,0){$M$}
      \put(1,.5){\line(0,1){1}}
      \put(1,2){\line(0,1){.3}}
      \put(1,2){\circle*{\dotsize}}
      \put(1,3){\line(0,1){.3}}
      \put(1,3){\circle*{\dotsize}}
      \put(2.0,0){$N$}
      \put(2.3,.5){\line(0,1){1}}
      \put(2.3,2.5){\line(0,1){.3}}
      \put(2.3,2.5){\circle*{\dotsize}}
      \put(2.3,3.5){\line(0,1){.3}}
      \put(2.3,3.5){\circle*{\dotsize}}
    \end{picture}
  \caption{fences for $\CC_M$ and $\CC_N$ when $M\cap N\notin M$ and
    $M\cap N\notin N$.}
  \label{fig:PB}
\end{wrapfigure}
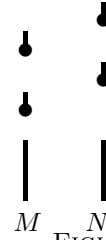
\noindent This is illustrated by
figure~\ref{fig:PB}, where the solid dots show the required fences for
compatibility of $\CC_M$ and $\CC_N$ in the case
where $M\cap N=M\cap H_{\sup(M\cap N)} = N\cap H_{\sup(M\cap N)}$.

The cases $M\cap N\in M$ and $M\cap N\in N$ are similar, except that
if, say $M\cap N\in M$, then $\sup(M\cap N)\in M$, that is, the bar in
$M$ is
longer, and the smallest fence is in $N$.

\begin{proof}[Proof of Lemma~\ref{thm:MNcompat}]
  To see that any $M$-fence $x'$ for $\CC_N$ must contain $x$, and
  that 
  therefore the existence of such a fence implies that $y$ is finite
  and $y\subset B$,
  suppose that $\lambda\in y$ and $\II_{\lambda}\notin x'$ and consider a
  requirement $\OO_{(\eta,\lambda]}$ where $\eta\in M\cap\lambda$,
  $\eta\geq\sup(M\cap N)$ if $M\cap N\in M$,  and
  $\eta>\max(\set{\tau<\lambda:\II_{\tau}\in x'})$.
  Then $\OO_{(\eta,\lambda]}$ is a member of $M$ which is compatible with
  $x'$; however it is 
  incompatible with $\CC_N$ since $(\eta,\lambda]\cap N\not=\nothing$
  because $\sup(N\cap\lambda)\ge\sup(M\cap\lambda)>\eta$ and
  $\OO_{(\eta,\lambda]}\notin N$ because $\lambda\notin N$. 
  
  To see that the stated conditions imply that $x$ is a $M$-fence for
  $\CC_N$, suppose that  $\OO_{(\eta',\eta]}\in M$ is incompatible
  with $C_N$.   If there is any ordinal $\gamma\in(\eta',\eta]\cap N$
  with $\gamma\geq\sup(M\cap N)$ then $\OO_{(\eta',\eta]}$ is
  incompatible with $\II_{\min(M\setminus \gamma)}\in x$, so we can
  assume that  $(\eta',\eta]\cap N\subseteq\sup(M\cap N)$.   Then
  $\eta'<\sup(M\cap N)$, so  according to
  definition~\ref{def:Pcompat}(\ref{item:MNf}) we need only consider
  the case $M\cap N\notin M$.  
  In this case $M\cap \sup(M\cap N)\subseteq N$, so $\eta<\sup(M\cap
  N)$ would imply $\sing{\eta',\eta}\subset N$, contradicting the
  assumption that 
  $\OO_{(\eta',\eta]}$ is incompatible with  $\CC_N$.   Thus we have
  $\eta'<\sup(M\cap 
  N)<\eta$.
  Then the statement of the lemma requires
  $\lambda\deq\min(M\setminus\sup(M\cap N))\in y$ so
  $\eta'<\lambda\leq\eta$, and
  $\OO_{(\eta',\eta]}$ is incompatible with $\II_{\lambda}\in x$.
\end{proof}

Whenever we refer to a $M$-fence for $\CC_N$ we will 
mean the minimal fence $x$ described in
proposition~\ref{thm:MNcompat}.
We will also refer to any of the individual
requirements in this minimal 
fence as an $M$-fence for $\CC_N$.

The fact that any 
superset of the minimal $M$-fence for $\CC_N$ is, according to
definition~\ref{def:Pcompat}(\ref{item:MNf}), also an $M$-fence  for
$\CC_N$ is something of an anomaly; however alternate definitions
which avoid this seem, at least in the forcing of
section~\ref{sec:moreclubs}, to be significantly more complicated.

\begin{corollary}\label{thm:nolim}
  If the requirements  $\CC_M$ and $\CC_{M'}$ are compatible then
  $\lim(M\cap M')=\lim(M)\cap\lim(M')$.
  \qed
\end{corollary}

Since the forcing $P_B$ is not separative, it will be convenient to
define notation for the equivalent separative forcing: if $\dot G$ is
a name for the 
generic set then we will say that $p'\lex p$ if
$p'\forces p\in\dot G$ and $p'\eqx p$ if $p'\lex p$ and $p\lex p'$.
The goal in this subsection is to prove the following lemma:

\begin{lemma}
  \label{thm:PBcomplete}
  Suppose that $p$ is a condition, and let $X$ be the finite set
  of ordinals $\lambda$ such that either 
  \begin{inparaenumi}
  \item\label{item:1}
    $\II_\lambda$ is one of the fences required for compatibility of two
    requirements in $p$,
  \item
    $\lambda=\sup(M\cap\lambda')$ for some
    $\CC_M\in p$ and some $\II_{\lambda'}$ which is either in $p$ or
    included in $X$ by clause~(\textup{\ref{item:1}}), or 
%%    $\II_{\lambda'}$
%%    $\sing{\II_{\lambda'},\CC_M}\subseteq p$, or
  \item $\lambda=\sup(M)$ for some $\CC_M\in p$.    
  \end{inparaenumi}
  Then $p'=p\cup\sing{\II_{\lambda}:\lambda\in X}$ is a condition and 
  $p'\eqx p$.   Furthermore, 
  $
    \forall\lambda<\kappa\;(
      (p\forces\II_\lambda\in\bigcup\dot G)\iff \II_{\lambda}\in p')
    $.
\end{lemma}

We first consider the fences:
\begin{lemma}\label{thm:compPfen}
  Suppose that $p\in P_B$, and $p'$ is the set obtained by adding to
  $p$ each of the fences required for compatibility of requirements in
  $p$.  Then $p'\in P_B$ and $p'\eqx p$, and every fence required for
  compatibility of members of $p'$ is a member of $p'$.
\end{lemma}
\begin{proof}
  We must show that each of the fences $\II_{\lambda}\in p'\setminus
  p$ is compatible with any requirement in $p$, and that any fence
  required for this compatibility is already a member of $p'$. 
  Suppose that $\CC_M\in p$ and $\II_{\lambda}$ is a $M$-fence for
  one of the requirements $\II_{\tau}$ or $\CC_N$ in $p$.
  \smallskip{}

  First we show that $\II_\lambda$
  is compatible with any requirement $\OO_{(\eta',\eta]}\in p$.
  Suppose to the contrary that ${\lambda}\in(\eta',\eta]$.
  Since $\II_\lambda\in M$, the compatibility of $\OO_{(\eta',\eta]}$
  with $\CC_M$ implies that 
  $\OO_{(\eta',\eta]}\in M$, so that 
  $\eta'<\sup(M\cap\lambda)<\lambda\le\eta$.
  If $\II_{\lambda}$ is a $M$-fence for $\II_{\tau}\in p$ then
  $\eta'\in M$ implies that $\eta'<\tau<\eta$, contradicting the
  compatibility of $\II_{\tau}$ and $\OO_{(\eta',\eta]}$.   If
  $\II_{\lambda}$ is a $M$-fence for $\CC_{N}\in p$ then $(\eta',\eta]\cap
  N\not=\nothing$, so $\OO_{(\eta',\eta]}$ is a member of $N$ as
  well as of $M$; however 
  this is impossible since $\eta\ge\lambda>\sup(M\cap N)$.
  \smallskip{}

  It remains to show that $\II_{\lambda}$ is compatible with any
  requirement $\CC_{M'}\in p$.
  Let $\lambda'=\min(M'\setminus\lambda)$, so that $\II_{\lambda'}$ is the
  $M'$-fence for $\II_\lambda$ required for compatibility of
  $\II_{\lambda}$ and $\CC_{M'}$.     We need to show that
  $\II_{\lambda'}\in p'$.

  If $\lambda\in {M'}$ then $\lambda'=\lambda$ and so
  $\II_{\lambda'}=\II_{\lambda}\in 
  p'$.   If $\lambda\notin M'$ and $\lambda\ge\sup(M\cap {M'})$ then
  $\II_{\lambda'}$ is a member of 
  $p'$ because it is a ${M'}$-fence for $\CC_M$.   Thus we can assume that
  $\lambda<\sup(M\cap {M'})$ and $\lambda\notin {M'}$, so that
  $\lambda<\lambda'<\sup(M\cap N)$.  
  Since $\lambda\in M\setminus M'$ it
  follows that $M\cap {M'}\in M$.

  If $\II_{\lambda}$ is an $M$-fence for $\II_{\tau}\in p$, then
  $(\tau,\lambda]\cap {M'}\subseteq(\tau,\lambda]\cap M=\nothing$, so
  $\lambda'=\min({M'}\setminus\lambda)=\min({M'}\setminus\tau)$ and hence
  $\II_{\lambda'}$ is in 
  $p'$ as the ${M'}$-fence for $\II_{\tau}$.

  Now suppose that $\II_{\lambda}$ is an
  $M$-fence for $\CC_{N}\in p$, that is,  
  $\lambda>\sup(M\cap {N})$ 
  and 
  $\lambda=\min(M\setminus\eta)$ for some $\eta<\lambda$ in  $N$.
  In the case that $\lambda\ge\sup({M'}\cap N)$ we claim that
  $\sup(N\cap\lambda')>\sup(M'\cap\lambda')$, so that 
  $\II_{\lambda'}$ is in $p'$ as an $M'$-fence for $\CC_N$.
  We have
  $\sup(N\cap{\lambda'})\ge
  \eta\geq\sup(M\cap\lambda)$, since $\lambda'>\lambda$.
  However $\sup(M\cap\lambda)>\sup({M'}\cap\lambda)$ since
  $\lambda<\sup(M\cap M')$,
  $\cof(\lambda)>\omega=\cof(\sup(M'\cap\lambda))$,  and $M\cap
  M'=M'\cap H_{\sup(M\cap M')}\in M$. 
  Finally, $M'\cap\lambda=M'\cap\lambda'$ since
  $\lambda'=\min(M'\setminus\lambda)$.  Hence
  $\sup(N\cap\lambda')>\sup(M'\cap\lambda')$, as claimed.  
  
  We will now show that the remaining case, $\lambda<\sup({M'}\cap
  N)$, is not possible.   If it did hold then we would have
  $\sup({M'}\cap N)>{\lambda'}$.
  Now ${\lambda'}\in M\cap {M'}$, since $M\cap {M'}$ is an initial
  segment of ${M'}$. 
  It follows that $N\cap {M'}$ is not an initial segment of ${M'}$, as this
  would imply that ${\lambda'}\in N$, contradicting the fact that
  ${\lambda'}\geq\lambda>\sup(M\cap N)$.  Hence $N\cap {M'}=N\cap H_{\sup(N\cap
    {M'})}\in {M'}$, and it follows that 
  $N\cap\lambda\in M$.
  Then $\sup(N\cap\lambda)<\sup(M\cap\lambda)$, but this contradicts
  the assumption that $\lambda$ is $M$-fence for $\CC_N$. 
  \smallskip{}

  This completes the proof that every fence required for the compatibility
  of requirements of $p'$ is already a member of $p'$,
  and hence that $p'\in P_B$.
  To see that $p'\eqx p$, let  $q\le p$ be arbitrary and let $q'\le q$
  be obtained from $q$  as in the lemma by adding to $q$ all of the
  fences required for the compatibility of requirements
  in $q$.   Then $q'\supseteq p'$, and hence $q'\le q$ forces that
  $p'\in\dot G$.
\end{proof}
\begin{lemma}\label{thm:compPsup}
  Suppose that $p\in P_B$, and $p'$ is obtained from $p$ by adding
  those requirements $\II_{\lambda}$ such that there is some $\CC_M\in
  p$ such that either $\lambda=\sup(M)$ or
  $\lambda=\sup(M\cap\lambda')$ for some $\II_{\lambda'}\in p$.  Then
  $p'\in P_B$ and $p'\eqx p$.

  Furthermore if $\CC_N\in p$ and $\II_{\lambda}\in p'\setminus p$
  with $\lambda<\sup(N)$ then the $N$-fence for $\II_{\lambda}$ either
  is equal to $\II_{\lambda}$ or else is an $N$-fence for some
  $\II_{\lambda'}\in p$.
\end{lemma}
\begin{proof}
  Again we need to show that every requirement $\II_{\lambda}\in p'$ is
  compatible with every requirement $\OO_{(\eta',\eta]}$ or $\CC_{{M'}}$
  in $p$.

  In order to  show that any requirement $\II_{\lambda}$ as specified in
  the lemma is compatible with any requirement $\OO_{(\eta',\eta]}\in
  p$, we will assume that $\II_{\lambda}$ is incompatible with  
  $\OO_{(\eta',\eta]}$ and show that $\OO_{(\eta',\eta]}$ is
  incompatible with $\CC_{M}$ or $\II_{\lambda'}$, contradicting the
  assumption that it is 
  in $p$.  Now $(\eta',\eta]\cap M\not=\emptyset$ since
  $\lambda\in\lim(M)$, so $\OO_{(\eta',\eta]}$ is incompatible with
  $\CC_M$ unless $\OO_{(\eta',\eta]}\in M$.   Since $\eta\geq\lambda$,
  this is impossible if
  $\lambda=\sup(M)$.   If $\lambda=\sup(M\cap\lambda')$ then
  $\OO_{(\eta',\eta]}\in M$ 
  implies that $\eta\geq\lambda'$, so $\OO_{(\eta',\eta]}$ is
  incompatible with $\II_{\lambda'}$.

  Now we show that $\II_{\lambda}$ is compatible with every
  requirement $\CC_{M'}\in p$.   This is immediate if
  $\lambda\ge\sup({M'})$.   If $\sup(M\cap {M'})\le\lambda<\sup({M'})$ then by
  proposition~\ref{thm:MNcompat} the fence
  $\II_{\min({M'}\setminus\lambda)}$ is a required ${M'}$-fence for $\CC_M$.
  Finally suppose that $\lambda<\sup({M'}\cap M)$.   If ${M'}\cap M\in {M'}$
  then $\lambda\in {M'}$ and hence $\II_{\lambda}$ is its own
  ${M'}$-fence.   Otherwise ${M'}\cap\sup({M'}\cap M)\subset M$ so  
  $\lambda\leq\lambda'\le\min(M\setminus\lambda)\le\min({M'}\setminus\lambda)$,
  so the ${M'}$-fence $\II_{\min({M'}\setminus\lambda)}$ for
  $\II_\lambda$ is the 
  same as the ${M'}$-fence $\II_{\min({M'}\setminus\lambda')}$ for
  $\II_{\lambda'}$. 

  This completes the proof that $p'\in P_B$ and that any nontrivial
  fences for members of $p'\setminus p$ are already fences for members
  of $p$.   To see that $p\lex p'$, notice that for any condition
  $q\le p$ we have $q'\supseteq p'$ and hence $q'\le p'$.   Thus
  $p\forces p'\in\dot G$.
\end{proof}

Let us call a condition $p\in P_B$ \emph{complete} if every fence
required for compatibility of requirements in $p$ is a member of $p$,
and if $\II_{\lambda}\in p$ whenever there is $\CC_M\in p$ such that
$\lambda=\sup(M)$ or $\lambda=\sup(M\cap\lambda')$ for some
$\II_{\lambda'}\in p$.

\begin{corollary}
  For any condition $p$ there is a complete condition $p'\eqx p$.
\end{corollary}
\begin{proof}
  Begin by using lemma~\ref{thm:compPfen} to add to $p$ all fences
  required for compatibility of $p$, and then use
  lemma~\ref{thm:compPsup} to add requirements of the form
  $\II_{\sup(M)}$ or $\II_{\sup(M\cap\lambda)}$. 
\end{proof}

\begin{lemma}\label{thm:OOout}
  Suppose that $p$ is a complete condition and $\lambda<\kappa$ is an
  ordinal such that $\II_{\lambda}\notin p$.   Then there is a
  requirement $\OO_{(\eta',\eta]}$ incompatible with $\II_{\lambda}$
  such that $p\cup\sing{\OO_{(\eta',\eta]}}\in P_B$.
\end{lemma}
\begin{proof}
  We may assume that there is $\CC_M\in p$ with $\sup(M)>\lambda$, for
  otherwise we could take $\OO_{(\eta',\eta]}=\OO_{(\eta',\lambda]}$
  where $\eta'$ is any sufficiently large ordinal less than $\lambda$.
  Since $\II_{\sup(M)}\in p$ for each $\CC_M\in p$, 
  it follows that there is some ordinal $\tau>\lambda$ with
  $\II_{\tau}\in p$.  Let $\tau$ be the least such.  
  
  If $\CC_M\in
  p$ then either $\sup(M\cap\tau)<\lambda$ or $\tau\in\lim(M)$, for
  otherwise 
  we would have $\lambda\leq\sup(M\cap\tau)<\tau$ and 
  $\II_{\sup(M\cap\tau)}\in p$, contradicting the choice of either
  $\lambda$ 
  or $\tau$.     Let $Y=\set{M:\CC_M\in p\And\tau\in\lim(M)}$. 
  Then
  $Y\not=\emptyset$, since otherwise we could take
  $\OO_{(\eta',\eta]}=\OO_{(\eta',\lambda]}$ 
  for any sufficiently large $\eta'<\lambda$.
  
  I claim that $\set{M\cap\tau:M\in Y}$ is linearly ordered by
  $\subseteq$.  To see this, note that
  $\tau\in\lim(M)\cap\lim(M')=\lim(M\cap M')$, so that $\sup(M\cap
  M')\geq\tau$.   The claim then follows from the condition
  Definition~\ref{def:Pcompat}(\ref{item:MNci}) for compatibility of
  $\CC_M$ and $\CC_{M'}$. 

  Now pick $M\in Y$ so that $M\cap\tau$ is as small as possible, and
  set 
  $\eta=\min(M\setminus\lambda)$.  If $\eta'$ is any member of
  $M\cap\lambda$ then $\OO_{(\eta',\eta]}\in M'$ for all $M'\in Y$.
  I claim that there is $\eta'\in M\cap\lambda$ such that
  \begin{equation*}
 \eta'>\max\bigl(\set{\xi<\lambda:\II_{\xi}\in
     p}\cup\set{\sup(M'\cap\tau):\CC_{M'}\in p\And M'\notin Y}\bigr).
  \end{equation*}
  It will follow that $\OO_{(\eta',\eta]}$ is compatible with every
  requirement in $p$, so that $p\cup\sing{\OO_{(\eta',\eta]}}\in P_B$. 

  To prove the claim we need to show that
  $\sup(M\cap\lambda)>\xi$ for all
  $\II_{\xi}\in p$ with $\xi<\lambda$, and 
  $\sup(M\cap\lambda)>\sup(M'\cap\lambda)$ for all $\CC_{M'}\in p$
  with $M'\notin Y$.

  If $\xi<\lambda$ and $\II_{\xi}\in p$ then
  $\II_{\min(M\setminus\xi)}\in p$.    Since $\II_{\eta}\notin p$ and
  $\eta=\min(M\setminus\lambda)$ it
  follows that $\xi\leq\min(M\setminus\xi)<\lambda$.   Hence
  $\sup(M\cap\lambda)>\xi$. 

  Now suppose that $\CC_{M'}\in p$ but $M'\notin Y$.  
  If $\sup(M'\cap M)\geq\tau$ then $\eta\in M\setminus M'$ implies
  that $M\cap M'\in M$, so $\sup(M'\cap\lambda)=\sup((M\cap
  M')\cap\lambda)\in M$ and hence
  $\sup(M\cap\lambda)>\sup(M'\cap\lambda)$. 
  Thus we can assume that $\sup(M'\cap M)<\tau$, so that $\sup(M'\cap
  M)<\lambda$.  If  $\xi$ is any member of $M'\cap\lambda$ with
  $\xi\geq\sup(M\cap M')$ then 
  $\II_{\min(M\setminus\xi)}$ is an $M$-fence for $M'$ and hence is in
  $p$.   Since $\II_{\eta}\notin p$ it follows that
  $\min(M\setminus\xi)<\lambda$.   Thus we can assume that
  $M'\cap\lambda\subset\sup(M\cap M')$.   If $M'\cap M\in M$ this
  implies 
  $\sup(M'\cap\lambda)=\sup(M\cap M')\in M$, while if $M'\cap M\notin M$
  then $\II_{\min(M\setminus\sup(M\cap M'))}\in p$ as a $M$-fence for
  $M'$, and as before this implies $\min(M\setminus\sup(M\cap
  N))<\lambda$.   Thus in any case we have
  $\sup(M\cap\lambda)>\sup(M'\cap\lambda)$. 
\end{proof}
\begin{proof}[Proof of lemma~\ref{thm:PBcomplete}]
  We already know that $p$ can be extended to a complete condition
  $p'$ so that $p'\eqx p$.    By lemma~\ref{thm:OOout}, if 
  $\II_{\lambda}\notin p'$ then there is $q\le p'$  so that
  $q\forces\II_{\lambda}\notin\bigcup\dot G$.
\end{proof}
\begin{definition}
  If $G$ is a generic subset of $P_B$ then we write $D$ for the set of
  $\lambda<\kappa$ such that $\II_{\lambda}\in\bigcup G$.
\end{definition}
\begin{corollary}
  The set $D$ is a closed and unbounded subset of $B^*$.   
\end{corollary}
\begin{proof}
  That $D$ is a subset of $B^*$ follows from the fact that
  $\II_{\lambda}$ is a requirement only if $\lambda\in B^*$.
  To see that $D$ is unbounded, suppose that $p\in P_B$ and
  $\eta<\kappa$.   Let $\lambda$ be any member of
  $\kappa\setminus\eta$ such that $\cof(\lambda)=\omega$ and $p\in
  H_\lambda$;  then   $\II_{\lambda}$ is a requirement which is
  compatible with $p$ and which forces that $\lambda\in D\setminus\eta$.

  Finally, let $\lambda<\kappa$ and $p\in P_B$ be arbitrary such that
  $p\forces\lambda\in\lim(\dot D)$.    Then $p$ is incompatible with
  any requirement $\OO_{(\eta',\eta]}$ with $\eta'<\lambda\le\eta$,
  and it follows by  lemma~\ref{thm:OOout} that
  $p\forces\lambda\in\dot D$.   Hence $D$ is closed.
\end{proof}

\subsection{Strongly generic conditions}\label{sec:smc}
It was pointed out in the discussion preceding the definition
of $P_B$ that the forcing $P_{\omega_1}$ satisfies a property
stronger than 
that of being proper, and it was stated as part of the motivation for
$P_{\omega_2}$ and hence for $P_B$ that these forcings would satisfy
the same property.   We now make this notion precise:
\begin{definition}\label{def:strgen}
If $P$ is a forcing notion and $X$ is a set then we say that $p$ is
\emph{strongly $X,P$-generic} if $p\forces_{P}\lq\lq\dot G\cap X$ is a
$V$-generic 
subset of $P\cap X$'' where $\dot G$ is a name for the generic set.%
\footnote{It should be noted that strong genericity as defined here 
  is not related to the notion which  Foreman,
  Magidor and Shelah \cite{formagshe.martins-max-i} call strong
  genericity.} 
 \end{definition}

Being strongly $X,P$-generic is stronger than Shelah's notion of 
a $P,X$-generic 
condition $p$, which only needs to force that $\dot G\cap X$ is a
$X$-generic 
subset of $P\cap X$.
Also, the existence of a strongly $X,P$-generic condition does not
require that $P\in X$, as does the existence of a $X,P$-generic
condition in Shelah's sense.

Definition~\ref{def:strgen} can be restated: $p_0$ is strongly
$X,P$-generic if, below the condition 
$p_0$, the forcing $P$ can be written as a two stage iteration.    If
we write $P/p_0$ for the forcing $P$ below the condition $p_0$,
Definition~\ref{def:strgen} 
implies that (with some abuse of notation) there is a $(P\cap X)$-term
$\dot R$ such that $P/p_0\fequiv ((P/p_0)\cap X)*\dot R$.   The
following equivalent definition of strong  genericity clarifies the
meaning of the notation $(P/p_0)\cap X$: 

%
%% \footnote{I have learned from Todorcevic that
%%   Baumgartner has  
%%   introduced, in unpublished lectures, a related but weaker notion
%%   which he calls
%%   \emph{strongly proper}.   He used this notion to show that countable
%%   support iteration of Cohen reals yields a model with no
%%   $\gw_2$-Aronszajn trees.}

\begin{proposition}
A condition $p_{0}\in P$ is strongly $X,P$-generic if and only if
\begin{inparaenumi}
\item 
  if 
  $p, q_0$ and $q_1$ are any conditions such that 
  $p\leq p_0$, $p\leq q_0,q_1$, and $\sing{q_0,q_1}\subset X$ then
  $q_0$ and $q_1$ are compatible in $P\cap X$, and 
\item 
  for every $p\le p_0$ in $P$ there is a
  condition $p\cut X\in P\cap X$ such that
  any condition $q\le p\cut X$
  in $X$  is compatible with $p$.  
\end{inparaenumi}
\end{proposition}
   
\begin{proof}
  First assume that $p_0$ is strongly $X,P$-generic.   If $p,q_0$ and
  $q_1$ are as in clause~(i) then $p\forces{q_0,q_1}\subset \dot G\cap
  X$, and since $p$ also forces that $\dot G\cap X$ is a generic
  subset of $P\cap X$ it follows that $q_0$ and $q_1$ are compatible
  in $P\cap X$.   For clause~(ii), suppose that $p\leq p_0$ and let
  $D$ be the set of $q\in P\cap X$ such that either $q$ is 
  incompatible with $p$ or every $q'\leq q$ in $P\cap X$ is compatible
  with $p$.  Then $p$ forces that $\dot G\cap D\not=\emptyset$, so
  there is some $q\in D$ which is compatible with $p$.
  This condition $q$  is  a
  suitable choice for $p\cut X$. 

  Now suppose that $p_0$ satisfies clauses~(i) and~(ii).
  First suppose that  $q_0,q_1$ are members of $P\cap X$ and $p\leq
  p_0$ forces that $\sing{q_0,q_1}\subseteq\dot G$.  We can assume, by
  extending $p$ if necessary, that $p,q_0$ and $q_1$ satisfy the
  hypothesis of clause~(i), which  implies that $q_0$ and $q_1$ are
  compatible in $P\cap X$.   Hence $p_0$ forces that $\dot G\cap X$ is
  a pairwise compatible subset of $P\cap X$.
  Now suppose that $D$ is a dense subset of $P\cap X$ and
  $p\leq p_0$.  Then there is some $q\leq
  p\cut X$ such that $q\in D$.  It follows that $p$ and $q$ are compatible in
  $P$, and  if $p'$ is any common extension of $p$
  and $q$ then $p'\forces q\in D\cap (\dot G\cap X)$.
  Thus $p_0$ forces that $\dot G\cap X$ is a generic subset of $P\cap X$.
\end{proof}

All of the forcing notions $P$ used in this paper will satisfy that
$q_0\cup q_1=q_0\wedge q_1$ for all compatible conditions $q_0,q_1\in
P$, and    hence clause~(i) will be satisfied by any set $X$ which is
closed under finite unions.   Thus we will only need to consider clause~(ii).
A function $p\mapsto p\cut X$ satisfying clause~(ii) will be called a
\emph{witness to the strong  $P,X$-genericity of $p_0$}.

We will usually omit $P$, writing ``strongly $X$-generic'' instead of
``strongly $X,P$-generic'', when it is clear which forcing notion is
meant.

We will say that a model $X$ has \smc{s} if for every $p\in P\cap X$
there is a strongly $X$-generic condition $p'\le p$.   In many
applications, including all the 
examples in this paper, there is a single strongly $X$-generic
condition $p_0$  which is compatible with every condition
$q\in P\cap X$.

The next two definitions are standard:
\begin{definition}\label{def:stationary}
  A set $Y\subset\ps(I)$ is \emph{stationary} if for every
  structure $\mathcal A$ with universe $I$ and in a countable
  language, there is a set $M\in Y$ with $\mathcal A|M\prec \mathcal A$.
\end{definition}
  
We will not normally specify the index set $I$.
Notice that the property of being a strongly $P,M$-generic condition
(unlike the property of being a $P,M$-generic condition) depends only on
$M\cap P$; hence the set $I$ can be taken to be the set $P$ of
conditions.   However we will also take advantage of the well known
fact that if $I'\supset I$ and $Y$ is a stationary subset of $\ps(I)$,
then $\set{x\subset I':x\cap I\in Y}$ is a stationary subset of
$\ps(I')$.    This observation makes it possible to apply the
stationarity of a given class 
$Y\subset \ps(P)$
to obtain an elementary substructure $M\prec \mathcal A$
where $\mathcal A$ is a model with universe properly containing $P$. 

\begin{definition}\label{def:presat}
  A forcing notion $P$ is said to be \emph{$\delta$-presaturated} if
  for any set $A\subset V$ in $V[G]$ with $\card{A}^{V[G]}\leq\delta$,
  there is a set $A'\supset A$ in $V$
  such that $\card{A}^V<\delta$.
\end{definition}
We use $\delta$-presaturation as a local version of the $\delta$-chain
condition: it is equivalent to the statement that for every collection
$\ca$ of fewer than $\delta$ antichains in $P$  there is a dense set
of conditions $p$ such that the set 
of conditions in $\bigcup \ca$ which are compatible with $p$ has size
less than 
$\delta$.   This ensures that forcing with $P$ does not collapse
$\delta$. 

The following is a well known observation. 

\begin{lemma}\label{thm:smccc}
  Suppose that $P$ is a forcing notion such that for  stationarily
  many models $M$ of
  size less than $\gd$ there is, for each $q\in M$, a $M$-generic
  condition $p\le q$.    Then $P$ is
  $\gd$-presaturated.
\end{lemma}

\begin{proof}
  Assume that $\dot A$ is a $P$-name for a subset $A$ of $V$ in
  $V[G]$ such that $\mu\deq\card{A}^{V[G]}<\delta$, and let $\dot k$ be a
  $P$-name such that $p\forces\dot k\colon\mu\xrightarrow{\text{onto}}\dot
  A$.  For any sufficiently large cardinal  $\theta$, 
  pick a model $M\prec H_{\gth}$ of size less than $\gd$   
  such that $\sing{\dot k,\dot
  A,p,P}\cup\mu\subset M$ and such that there is a $M$-generic
  condition  $p_0<p$.   Then $p_0$
  forces that for every $\xi<\mu$ there is $q\in M\cap\dot G$ and
  $x\in M$ such that $q\forces\dot k(\xi)=x$, and hence $p_0$ forces
  that $\dot A\subset M$.  
\end{proof}
\begin{corollary}\label{thm:smccc1}
  If $P$ is a forcing notion such that the trivial condition $\1^{P}$
  is $M$-generic for stationarily many sets $M$ of size less than
  $\delta$, then $P$ has the  $\delta$-chain condition.
\end{corollary}
\begin{proof}
  Let $\mathcal{A}$ be a maximal antichain in $P$, and apply the proof of
  the lemma with the singleton $G\cap\mathcal{A}$  as 
  $A$ and with $p_0=\1^{P}$.
\end{proof}

We say that a forcing notion $P$ \emph{has meets} if any
compatible pair $p,q$ of conditions has a greatest lower bound $p\land q$.
The following definition states another property shared
by all \smc{s} in this paper: 
\begin{definition}\label{def:tidy} 
  If $P$ is a notion of forcing with meets and $X$ is a set then 
  we say that a strongly $X,P$-generic condition $p$ is 
  \emph{tidy} if  
  there is a 
  function $q\mapsto q\cut X$ witnessing the 
  strong $X$-genericity of $p$ such that 
  $(q\land q')\cut X=q\cut X\land q'\cut X$ whenever 
  $q, q'\le p$ are compatible.
\end{definition}

\begin{proposition}\label{thm:tidy-easy}
  Suppose that a strongly $X$-generic condition is tidy with
  witnessing function $q\mapsto q\cut X$.   Then
  \begin{inparaenumi}
  \item $q'\cut X\leq q\cut X$ for all $q'\leq q\leq p$, and 
  \item $q\lex q\cut X$ for all $q\leq p$.
  \end{inparaenumi}
\end{proposition}
\begin{proof}
  For clause~(i), we have $q'\cut X=(q'\wedge q)\cut X=(q'\cut
  X)\wedge(q\cut X)\leq q\cut X$.   For clause~(ii), if $q\not\lex
  q\cut X$ then there is $q'\leq q$ such that $q'\forces q\cut
  X\notin\dot G$; however $q'\wedge (q'\cut X)\leq q'$  and by
  clause~(i) $q'\cut X\leq q\cut X$.
\end{proof}

The next lemma states the critical fact which makes the existence of
strongly generic condition necessary to the constructions in this
paper. 

\begin{lemma}\label{thm:approx}
  Suppose that $p$ is a tidy strongly $X,P$-generic condition, 
  and that stationarily many models $M$ of size $\gd$ have \smc{s} for
  $P$.
  Let $G$ be a generic subset of $P$ with $p\in G$, and suppose $k\in
  V[G]$ is a function 
  with domain $\mu\in V$ such that $k\restrict x\in V[G\cap X]$ for
  each $x\in([\mu]^{\gd})^{V}$.   Then $k\in V[G\cap X]$.
\end{lemma}
Recall that the strong $X,P$-genericity of $p$ forces that $\dot G\cap
X$ is a $V$-generic subset of $P\cap X$.   Thus any two conditions
$q,q'\leq p\cut X$ are compatible in $P\cap X$ if and only if they
are compatible if $X$, and if $q\leq p\cut X$ is in $X$ and $\phi$ is
any formula then $q\forces_{P\cap X}\phi(\dot G\cap X)$ if and only if
$p\wedge q\forces_{P}\phi(\dot G\cap X)$.
\begin{proof}
  Let $\dot k$ be a name for $k$ and let $p_0\le p$ be a
  condition which forces that $\dot k$ satisfies the hypothesis of the
  lemma. 
  Let $\theta$ be a cardinal larger than $\kappa$ and 
  pick a model $M\prec H_{\gth}$ of size $\gd$  such that
  $\set{P_B,p_0,X,\dot k}\subset M$, the function $q\mapsto q\cut X$
  is in $M$,  and there is a strongly
  $M$-generic condition $p_1\le p_0$.
  Let $p_2\le p_1$ be a condition
  such that $p_2\forces \dot k\restrict(\mu\cap M)=\dot s$ for some $P\cap
  X$-term~$\dot s$.
  Note that if $r\le p_2$ is any 
  condition such that $r\forces\dot s(\nu)=x$ for some $\nu\in \mu\cap
  M$ then 
  $r\cut X\forces\dot s(\nu)=x$, as otherwise there would be $r'\le 
  r\cut X$ in $X$ such that $r'\forces\dot s(\nu)\not=x$, which is
  impossible since $r'$ and $r$ are compatible. 

  We will show that 
  \begin{equation}\label{eq:6}
    M\models \forall q\le p\wedge(p_2\cut
    M)\;\forall \nu\in\mu\;\bigl(q\decides \dot k(\nu)
    \implies  (p\wedge(p_2\cut M)\land
    (q\cut X))\decides \dot k(\nu)\bigr).
  \end{equation}
  Here $p\wedge(p_2\cut M)$ and $q\cut X$ are compatible since $q$ and
  $q\cut X$ are compatible, and $q\land q\cut X\le q\le p\wedge(p_2\cut M)$.
  Furthermore, since the three conditions $q$, $p_2\cut M$ and $q\cut X$
  are compatible the  condition $p\wedge (p_2\cut M)\wedge (q\cut
  X)$ 
  in formula~\eqref{eq:6} must decide $\dot k(\nu)$ in the same way as
  $q$ does.
  It may also be noted that in the forcings used in this paper, and in
  most likely applications of lemma~\ref{thm:approx}, the inclusion of $p$ in
  formula~\eqref{eq:6} is unnecessary, as $p_2\leq p\in M$ implies $p_2\cut
  M\leq p$.

  Suppose to the contrary that formula~\eqref{eq:6} is not valid, so
  that  there are 
  $q\le p\wedge(p_2\cut M)$ in $M$, 
  $\nu\in M\cap\mu$ and $x\in M$ such that 
  $q\forces\dot k(\nu)=x$ but for some 
  $r\le p\wedge(p_2\cut M)\land (q\cut X)$ in $M$ we have $r\forces\dot
  k(\nu)\not=x$. 
  Then $q\land p_2\forces\dot s(\nu)=\dot k(\nu)=x$, so 
  $(q\land p_2)\cut X\forces\dot s(\nu)=x$.
  Now $r\le p_2\cut M$ implies that $r\land p_2$ is a condition, and
  $r\land p_2\forces\dot s(\nu)=\dot k(\nu)\not=x$, so $(r\land
  p_2)\cut X\forces\dot s(\nu)\not=x$.  Thus $(r\land p_2)\cut X$ is
  incompatible with $(q\land p_2)\cut X$.

  Now $r\cut X\lex q\cut X$ in $P\cap X$, since otherwise there is some
  $r'\le r\cut X$
  in $X$ which is incompatible with $q\cut X$, but then $r'$ is
  compatible with $r\le q\cut X$ and hence with $q\cut X$.
  Thus 
  $(r\land p_2)\cut X=r\cut X\land p_2\cut X\lex q\cut X\land p_2\cut
  X=(q\land p_2)\cut X$, again in $P\cap X$.    Hence 
  $(r\land p_2)\cut X$ and  $(q\land p_2)\cut X$ are compatible, and 
  this contradiction completes the proof of
  formula~\eqref{eq:6}. 
  \medskip{}

  By elementarity $V$ also satisfies the right side of
  formula~\eqref{eq:6}.
  Since 
  $q\forces q\cut X\in\dot G$ for any $q\le p$ it follows
  that $p\wedge(p_2\cut M)$ forces that $\dot k\in V[\dot G\cap X]$:
  \begin{equation*}
    \forall\nu<\mu\forall x\;\bigl(k(\nu)=x\iff\exists  q'\in
      (G\cap X)\quad( p\wedge (p_2\cut M)\wedge q' )\forces\dot k(\nu)=x\bigr).
  \end{equation*}
  To see this, suppose $V[G]\models k(\nu)=x$.  Then there is $q\leq
  p\wedge (p_2\cut M)$ in $G$ such that $q\forces \dot k(\nu)=x$, but
  then $q'=q\cut X\in G\cap X$ is as required.  
\end{proof}
Another application of the  idea of this proof is given in
\cite{mitchell04:nhal}, where it is used to give 
an easier proof of the 
main lemma of~\cite{mitchell.atit} and of a related lemma of
Hamkins~\cite{hamkins:ggft}.   

\subsection{Strongly generic Conditions in $P_B$}

\begin{lemma}\label{thm:kappa-smc}
  If $\gl\in B$ then the condition
  $\sing{\II_{\gl}}$ is tidily strongly $H_\lambda$-generic.
\end{lemma}
\begin{proof}
  Define $p\cut H_\lambda$ for $p\le\sing{\II_{\gl}}$ to be 
  $p\cut H_\lambda\deq (p\cap H_\lambda)\cup\set{\CC_{M\cap
  H_\lambda}:\CC_M\in p}$.  

  It is straightforward to verify that $p\cut H_\lambda$ is a condition, and
  it is clearly tidy since each member of $(p\cut H_\lambda)\setminus p$ is
  determined by a single member of $p$ other than $\II_{\lambda}$.
  \medskip{}
  
  To see that the function $p\mapsto p\cut H_\lambda$ witnesses the strong
  $H_\lambda$-genericity 
  of $\sing{\II_\lambda}$, suppose that $q\le p\cut H_\lambda$ is in
  $H_{\lambda}$.   We need to show that the requirements in $p\cup q$
  are compatible.    We will show that any requirement in $p$ is
  compatible with each requirement in $q$.   

  Any requirement $\II_{\tau}\in p$ with $\tau\ge\lambda$ is compatible
  with any requirement in $H_{\lambda}$ and in particular with any requirement
  in $q$; while any requirement $\II_{\tau}\in p$ with $\tau<\lambda$
  is in 
  $H_{\lambda}$ and hence is a member of $q$.  Similarly, 
  the assumption that $p\le\sing{\II_\lambda}$ 
  ensures that any requirement of the form $\OO_{(\eta',\eta]}\in p$
  either satisfies $\lambda\le\eta'$, in which case it is compatible with
  any condition in $H_{\lambda}$, or else it satisfies $\eta<\lambda$,
  in which 
  case $\OO_{(\eta',\eta]}\in H_{\lambda}$ and hence 
  $\OO_{(\eta',\eta]}\in p\cut H_\lambda\subset q$.  
  
  In the case of a requirement $\CC_N\in p$ we have $\CC_{N\cap H_\lambda}\in
  p\cut H_\lambda\subseteq q$.    
  Any requirement $\OO_{(\eta',\eta]}\in H_\lambda$ which is compatible
  with $\CC_{N\cap H_\lambda}$ is also compatible with $\CC_N$.
  A requirement 
  $\II_{\tau}\in H_{\lambda}$ which is 
  compatible with $\CC_{N\cap H_{\lambda}}\in p\cut H_\lambda$ must be
  compatible with $\CC_N$, using the same fences, unless
  $\sup(N)>\lambda>\tau\geq\sup(N\cap\lambda)$, and in 
  that case the required $N$-fence is
  $\II_{\min(N\setminus\tau)}=\II_{\min(N\setminus\lambda)}$, which is
  required by the compatibility of $\II_{\lambda}$ and $\CC_{N}$.

  Finally, if $\CC_{N'}\in q$ then $N'\cap N=N'\cap(N\cap H_\lambda)$,
  and so $\CC_{N}$ and $\CC_{N'}$ satisfy clause~\ref{item:MNci} of
  definition~\ref{def:Pcompat} in the same way that
  $\CC_{N'}$ and $\CC_{N\cap H_{\lambda}}$ do.   The $N'$-fence for
  $\CC_{N\cap H_\lambda}$ is also a $N'$-fence for $\CC_{N}$, and a
  $N$-fence for $\CC_{N'}$ is given by the $N\cap
  H_\lambda$-fence for $\CC_{N'}$ together with
  $\II_{\min(N\cap\lambda)}$. 
\end{proof}
\begin{corollary}\label{thm:kappa-cc}
  If $\gk$ is inaccessible and $B$ is stationary in $\gk$, then $P_B$
  is $\gk$-presaturated and hence preserves all cardinals greater than
  or equal to $\gk$.
\end{corollary}
\begin{proof}
  This is immediate from lemmas~\ref{thm:kappa-smc} and~\ref{thm:smccc}.
\end{proof}
As was pointed out earlier, lemma~\ref{thm:compat} below, like
lemma~\ref{thm:kappa-smc} above, is a  variation of the proof of
properness for $P_{\gw_1}$.  Lemma~\ref{thm:compat} replaces the
countable set 
$M\prec H_{\gw_1}$ with a countable set $M\prec H_{\gk}$.
 
\begin{lemma}
  \label{thm:compat}
  If $\CC_M$ is a requirement then $\sing{\CC_M}$ is tidily strongly
  $M$-generic.   
\end{lemma}
\begin{proof}
  The proof is similar to that of lemma~\ref{thm:kappa-smc}, but is
  more complicated because $M$ is not transitive.
  For a condition $p\le\sing{\CC_M}$, let $\cp{M}p$ be the set of all
  $M$-fences required for the compatibility 
  of $\CC_M$ with other members of $p$.    We define the map $p\mapsto
  p\cut M$ by 
  \begin{equation*}% \label{eq:compat}
      p\cut M=(p\cap M)\cup\cp{M}p\cup\set{\CC_{N\cap M}:\CC_N\in
      p\And N\cap M \in M}. 
  \end{equation*}
  To see that $p\cut M$ is a condition, note that $(p\cap M)\cup
  \cp{M}p$ is a 
  condition because it is a subset of the condition $p'\supset p$,
  given by 
  lemma~\ref{thm:compPfen}, which contains all $\CC_N$-fences for all
  $\CC_N\in p'$.
  Since $N\cap M$ is an initial segment of $N$ for all $\CC_N\in p$
  with $N\cap M\in M$, it is easy to see that the result of adding the
  requirements $\CC_{N\cap M}$ is still a requirement.

  The function $p\mapsto \cp{M}p$ is tidy, since each member of
  $\cp{M}p\setminus p$ is determined  by a single
  member of $p$ other than $\CC_M$.   Each member $\CC_{N\cap M}$ of
  $(p\cut 
  M)\setminus\cp{M}p$ is also determined by the single member $\CC_N$
  of $p$, and hence the full map $p\mapsto p\cut M$ is tidy.
  \medskip{}

  In order to show that the function 
  $p\mapsto p\cut M$ witnesses the strong $M$ genericity of
  $\sing{\CC_M}$,   
  we need to show that if $q\le p\cut M$ is in $M$ then 
  any requirement in $p$ is compatible with every
  requirement in $q$. 

  First consider a requirement $\II_{\tau}\in p$.  If $\tau\ge\sup(M)$
  then $\II_{\tau}$ is compatible with any requirement in $M$, and if
  $\tau\in M$ then $\II_{\tau}\subset p\cap M\subset q$.   Hence we
  can assume that $\tau\in\sup(M)\setminus M$.   Set
  $\tau'=\min(M\setminus\tau)$.
  Then any requirement  $\OO_{(\eta',\eta]}$ or $\CC_N$ in $M$ which
  is incompatible with $\II_{\tau}$ would also be incompatible with
  $\II_{\tau'}\in\cp{M}p\subset q$, so $\II_{\tau}$ is compatible with
  every requirement in $q$. 

  \smallskip{}

  Any requirement $\OO_{(\eta',\eta]}\in p$ is  compatible with
  $\CC_M$, and thus either $(\eta',\eta]\cap M=\nothing$, in which case
  $\OO_{(\eta',\eta]}$ is compatible with any requirement in $M$, or
  else $\OO_{(\eta',\eta]}\in M$, in which case 
  $\OO_{(\eta',\eta]}\in p\cap M\subseteq p\cut M\subseteq q$. 
  \smallskip{}

  The case of a requirement $\CC_N\in p$ is somewhat more
  complicated than the previous two.   We first show that every
  requirement  $\II_{\tau}\in 
  q$ is compatible with $\CC_{N}$.   If $\tau\ge\sup(N)$ then 
  $\II_{\tau}$ is compatible with $\CC_N$, and if
  $\sup(N)>\tau\ge\sup(M\cap N)$ then 
  the required $N$-fence for $\II_{\tau}$ is a member of the 
  $N$-fence for $\CC_M$.  Thus we can suppose that  $\tau<\sup(M\cap
  N)$. 
  If $M\cap N$ is an initial segment of $N$ then it follows that
  $\tau\in N$, so we can also suppose that $M\cap N\in M$.   Then
  $\CC_{M\cap N}\in q$, and the required $N$-fence for $\II_{\tau}$ is 
  the same as the $(M\cap N)$-fence for $\II_{\tau}$ required for the
  compatibility of $\sing{\CC_{M\cap N},\II_{\tau}}\subseteq q$.

  Now we show that any requirement $\OO_{(\eta',\eta]}\in q$ is
  compatible with $\CC_N$.      If $(\eta',\eta]\cap N=\nothing$ then
  this is immediate, so we can assume that there is some $\xi\in
  (\eta',\eta]\cap N$.
  We cannot have $\xi>\sup(M\cap N)$, since in that case
  $\OO_{(\eta',\eta]}$ would be incompatible with a member of the
  $M$-fence for $\CC_N$, which is contained in $\cp{M}p\subseteq q$.
  Thus  
  $\eta'<\delta\deq\sup(M\cap N)$.   
  If $M\cap N\notin M$ then 
  $\eta<\delta$ as well, as otherwise $\OO_{(\eta',\eta]}$ would be
  incompatible with $\II_{\min(M\setminus\delta)}$, which is  a member
  of the $M$-fence for $\CC_N$.  
  But $M\cap N\notin M$ implies $M\cap
  \delta\subseteq N$, so $\OO_{(\eta',\eta]}\in N$ and thus
  $\OO_{(\eta',\eta]}$ is compatible with
  $\CC_N$.
  If, on the other hand, $M\cap N\in M$ then $\CC_{M\cap N}\in q$, and
  the compatibility of $\OO_{(\eta',\eta]}$ with $\CC_{N}$ follows from
  its compatibility with $\CC_{M\cap N}$.
  
  Finally we show that $\CC_N$ is compatible with any requirement
  $\CC_{N'}\in q$.     We verify clause~\ref{item:MNci} first.  In the
  case that $M\cap N$ is an initial segment of $M$, the set  $N\cap
  N'$ is also an initial 
  segment of $N'$.  On the other hand $N\cap N'$ is a countable subset of
  $H_{\delta}$ in $M$, and since  the cardinal
  $\delta'\deq\min(M\setminus\delta)$ is in $B$ it 
  follows that $N\cap N'\in M\cap H_{\delta'}=M\cap H_{\delta}=M\cap
  N$.   Thus $N\cap N'\in N$.
 
  In the  other case, when $M\cap N\in M$, we have $N\cap N'=(N\cap M)\cap
  N'$.   Since $q$ is a condition this is 
  an initial segment of one of $N'$ and  $N\cap M$, and 
  either an initial segment or
  a member of the other.    Now $N\cap N'$ will stand
  in the same relation to $N$ as it does to $N\cap M$.
  Thus $\CC_N$ and $\CC_{N'}$ satisfy clause~\ref{item:MNci}.

  It remains to verify that the necessary fences exist. 
  If $M\cap N$ is an initial segment of $M$, then any $N$-fence for
  $\CC_M$ is also a $N$-fence for $\CC_{N'}$.   Otherwise the union of
  a $N$-fence for $\CC_M$ with a $M\cap N$-fence for $\CC_{N'}$ gives
  a $N$-fence for $\CC_{N'}$.
  
  If $M\cap N$ is an initial segment of $M$ then a $N'$-fence for
  $\CC_N$ can be obtained by taking the set of all $N'$-fences for
  members of the $M$-fence for $\CC_N$, and otherwise the $N'$-fence
  for $\CC_N$ can be obtained by taking the union of this set with a
  $N'$-fence for $\CC_{N\cap M}$.

  This concludes the proof that any requirement in $M$ which is
  compatible with $p\cut M$ is compatible with a requirement $\CC_N\in
  p$, and hence finishes the proof of lemma~\ref{thm:compat}.
\end{proof}
\begin{corollary}\label{thm:proper}
  The forcing $P_B$ is proper.
  \qed
\end{corollary}

\begin{lemma}
  \label{thm:collapse}
  If $B$ is stationary and  $G$ is a $V$-generic subset of $P_B$ then
  $\gw^{V[G]}_1=\gw_1$, 
  $\gw_2^{V[G]}=\gk$, and all larger cardinals are preserved. 
\end{lemma}
\begin{proof}
  Corollary~\ref{thm:proper} implies that $\gw_1$ is preserved, and 
  corollary~\ref{thm:kappa-cc} implies that $\gk$ is preserved.
  All larger cardinals are preserved since $\card{P_B}=\gk$.

  Thus we only need
  to show that each cardinal 
  $\gl$ in the interval $\gw_1<\gl<\gk$ is collapsed.
  To see this, let 
  $Y\deq \set{M\cap\lambda:\lambda\in M\And \CC_M\in \bigcup G}$.
  If $\CC_M$ and
  $\CC_{M'}$ are compatible and $\lambda\in M\cap M'$ then 
  clause~\ref{item:MNci} of
  definition~\ref{def:Pcompat} implies that either $M\cap
  \lambda\subset M'$ or $M\cap \lambda\subset M$, so $Y$ is linearly
  ordered by subset.   Since each member of $Y$ is countable, it
  follows that $\card{Y}\le\omega_1$ and hence $\card{\bigcup
  Y}=\omega_1$ in $V[G]$.   But $\bigcup Y=\lambda$, 
  since for any condition $p\in P_B$ and
  any ordinal $\xi<\lambda$ we can find a countable set $M\prec
  H_{\kappa}$ with $\sing{p,\xi,\lambda}\subset M$, so that 
  $p\cup\sing{\CC_M}$ 
  is a condition extending $p$ which forces that $\xi\in
  M\cap\lambda\in Y$.
\end{proof}

\begin{comment}
\note[4/10/03]{The following lemma is not certainly not true of
  \emph{all} inaccessible $\gl<\gk$.   For example, the 
  first member of $B\cap D$ will be the $\gw_1$st member of $D$.  (And
  if $\gl$ is any non-Mahlo cardinal, with say $E\subset\gk$ closed, unbounded with
  no inaccessible cardinals, then $D\cap E$ has ordertype $\gw_1$ and
  is cofinal in $\gl$.)

  Can
  the hypothesis be weakened?  Is it true, for example, of the first
  Mahlo cardinal.}  
\note[4/10/03]{QUESTION: is $V[D]=V[G]$?   Note that 
  $\gd\deq\min(B\cap D)$ satisfies that $\set{M\cap \gd:\exists p\in
  G\;\CC_M\in p}$ is well quasi-ordered by $\in$.   In particular,
  this gives an explicit collapse of $\gd$ onto $\gw_1$.   Is this
  collapse in $V[D\cap\gk]$ (or even $V[D]$)?  (Of course $\gd$ is
  collapsed in $V[D]$ by the covering lemma, if for example V~=~L.)}

\note[4/13/03]{Contrary to the last notes, This lemma is true of the
  first member of $D\cap B$, I think, if things are defined
  right.   I'm not sure if it is now (presumably taking
  $G\cut\gl=\set{p\cut\gl:p\in G}$ for $D\cap\gl$).   However if I
  weaken the definition of a conditions, requiring the $M$ occurring in
  a requirement $\CC_M$ to be only $M\prec_{\gS_0}H_\gk$ instead of
  $M\prec H_\gk$, then the proof below works fine for any $\gl\in
  B\cap D$.    One could also ask that $M$ be closed under $\pi$,
  satisfy ZFC, and so on with no problems.

  Note that if $\gl\in D\cap B\setminus\lim(B)$ then the collapse map
  of $\gl$ \emph{is} in $V[D\cap\gl]$, so this is no counterexample.}
\end{comment}

\begin{lemma}
  \label{thm:noboundedcollapse}
  If  $\gl\in D\cap B$   then
  every function $\tau\colon\gw_1\to V$ in $V[G]$ such that
  $\forall\xi<\gw_1\;(\tau\restrict\xi\in V[G\cap H_\gl])$ is in $V[G\cap%
  H_\gl]$. 
\end{lemma}
\begin{proof}
  This is immediate from lemmas~\ref{thm:approx}, \ref{thm:kappa-smc}
  and~\ref{thm:compat}. 
\end{proof}

The following observation explains why this forcing is relevant to the
ideal $I[\omega_2]$:

\begin{proposition}
  Suppose that $B\subset\kappa$ is a set of inaccessible cardinals in
  $V$ and that $G$ is a generic subset of $P_B$.   Then in $V[G]$ the
  restriction of the ideal 
  $I[\omega_2]$ to ordinals of cofinality $\omega_1$ is 
  generated by the nonstationary ideal on $\omega_2$ together with the
  single  set
  $S=\set{\lambda\in \kappa:B\cap\lambda \text{ is nonstationary
  in }\lambda}$. 

  Furthermore, any stationary subset of $B\setminus S$ in $V$ 
  remains stationary in $V[G]$.
\end{proposition}
\begin{proof} 
  To see that $S\in I[\omega_2]$,
  pick for each $\lambda\in S$ a
  closed unbounded set
  $E_\lambda\subset\lambda$ in $V$ which is disjoint from $B$.   If
  $\lambda\in \lim(D)\cap S$ and $\cof(\lambda)>\omega$ then the set
  $c_{\lambda}\deq E_\lambda\cap D$ is cofinal in $\lambda$, but has order
  type $\omega_1$ since any member of $D$ of uncountable cofinality is
  in $B$.
  
  Let $A_{\nu}=\set{a_\nu\cap D:\nu<\kappa}$ where
  $\set{a_\nu:\nu<\kappa}$ enumerates in $V$ the bounded subsets of
  $\kappa$, and let $F$ be the closed unbounded set of
  $\lambda<\kappa$ such that every bounded subset of $\lambda$ in $V$
  is a member of $\set{a_\nu:\nu<\lambda}$.    
  Then  $S\cap (\lim(D)\cap F)\subset B(A)$.
  \smallskip{}

  To see that no stationary subset of $B\setminus S$ is in 
  $I[\omega_2]$, let
  $A=\seq{a_{\nu}:\nu<\kappa}$ be an arbitrary sequence in $V[G]$, and
  let 
  $\dot A$ be a name for $A$.  
  Fix a continuous increasing elementary chain
  $\seq{X_{\nu}:\nu<\kappa}$ of elementary substructures of $H_{\kappa^+}$
  with $\dot A\in X_0$, and let $F$ be the closed unbounded set of
  cardinals $\lambda<\kappa$ such that $X_{\lambda}\cap
  H_\kappa=H_\lambda$.

  We will show
  that $F\cap D\cap (B\setminus S)$ is disjoint from $B(A)$.
  Suppose to the
  contrary that $\lambda\in F\cap D\cap(B\setminus S)$ and let
  $c\subset\lambda$ witness that $\lambda\in B(A)$.
  Since the strongly $X_{\lambda}$-generic condition $\sing{\II_{\lambda}}$
  is in $G$, the set $a_\nu=\dot a_\nu^{G}\in V[G\cap H_\lambda]$ for
  all $\nu<\lambda$.  Hence 
  $c\cap\nu\in V[G\cap H_{\lambda}]$ for each $\nu<\lambda$, and it
  follows by lemma~\ref{thm:noboundedcollapse} that $c\in V[G\cap
  H_{\lambda}]$.   However this is impossible:
  $G\cap
  H_{\lambda}$ is a generic subset of $P_B\cap H_{\lambda}=P_{B\cap
  \lambda}$ and  $B\cap\lambda$ is a  stationary subset of
  $\lambda$, so lemma~\ref{thm:collapse} implies that  $\lambda$ is not
  collapsed in $V[G\cap H_{\lambda}]$.
  \smallskip{}

  To see that any stationary subset of $B\setminus S$ remains
  nonstationary in $V[G]$, let $T\subseteq B\setminus S$ be stationary
  and let 
  $\dot E$ be a name for a closed unbounded subset $E$ of $\kappa$.
  Now 
  pick a continuously increasing sequence of elementary substructures
  $X_{\nu}$ of $H_{\kappa^+}$ with $\dot E\in X_0$ and let
  $F$ be the 
  closed unbounded set of cardinals $\lambda\in\kappa$ such that
  $X_\lambda\cap H_\kappa=H_\lambda$.
  Then $T\cap F$ is unbounded in $\kappa$, and since $T\subset B$ any
  condition $p\in P_B$ is compatible with $\sing{\II_\lambda}$ for
  some $\lambda\in T\cap F$.
  Since  $\sing{\II_\lambda}$ is strongly $H_\lambda$-generic, it
  forces that  $\dot E$ is unbounded in $\lambda$, and hence that
  $\lambda\in\dot E$.   Thus $p\cup\sing{I_\lambda}\forces\lambda\in
  T\cap\dot 
  E$.
\end{proof}
One other application of this forcing is of interest: like the forcing
described in \cite{mitchell.atit} it gives a model with no
special $\omega_2$-Aronszajn trees if $\kappa$ is Mahlo in $V$, and no 
$\omega_2$-Aronszajn trees if $\kappa$ is weakly compact in $V$.
The proof is the same as in \cite{mitchell.atit}, with
lemmas~\ref{thm:approx}, \ref{thm:kappa-smc} and~\ref{thm:compat}
taking the place of the main lemma in that paper. 
It would perhaps be hard to argue that 
this construction is simpler than that of \cite{mitchell.atit},
especially in view of the fact that (as is pointed 
out in \cite{mitchell04:nhal}) the proof of main
lemma of \cite{mitchell.atit} can be substantially simplified by using
the idea of the proof of 
lemma~\ref{thm:approx}.    However it is shown in 
\cite{mitchell.acus} that if the current forcing is
simplified by eliminating requirements of the forms $\II_{\lambda}$
and $\OO_{(\eta',\eta]}$, and using clause~\ref{item:MNci} of
definition~\ref{def:Pcompat} as the 
only compatibility condition,  then the generic extension is still a
model with no 
special $\omega_2$-Aronszajn trees, or no $\omega_2$-Aronszajn trees.  This
is certainly 
the simplest construction known of such a model, and is likely the
simplest possible.  

\section{Adding $\gk^+$ closed, unbounded subsets of $\gk$}
\label{sec:moreclubs}

We will now extend the forcing from section~\ref{sec:oneclub} in order
to construct  a
sequence $\seq{D_{\ga}:\ga<\gk^+}$ of closed, unbounded subsets  of
$\gk$.  This  sequence will be continuously diagonally decreasing,
which 
means that $D_{\ga+1}\subset D_\ga$ for all $\ga$, and that if $\ga$ is a
limit ordinal then 
$D_{\ga}$ is  equal to the diagonal intersection
$\dinter_{\ga'<\ga}D_{\ga'}=
\set{\nu:(\forall\ga'\in\pi_{\ga}\image\nu)\;\nu\in D_{\ga'}}$. 
The definition of this diagonal 
intersection will depend on a choice of maps
$\pi_{\ga}\colon\gk\cong\ga$.   In addition, the sets $D_\ga$ will be
subsets  of $B^*_\ga$, where $B_\ga=\set{\nu<\gk:\nu\text{ is 
  }f_{\ga}(\nu)\text{-Mahlo}}$, and the definition of the set $B_\ga$
depends on the 
choice of the function $f_\ga$ representing $\ga$ in the nonstationary ideal.
The first
subsection  describes how to use $\square_\gk$ to define the functions
$\pi_\ga$ and~$f_{\ga}$.

We assume throughout this section that 
$\square_{\kappa}$ holds. 
We also 
assume throughout the section that $\kappa$ is inaccessible and that
$2^{\kappa}=\kappa^{+}$, but only 
in the final subsection~\ref{sec:finish} will we make use of the
assumption that 
$\kappa$ is $\kappa^+$-Mahlo.
\subsection{Using $\square_\gk$}
\label{sec:square}
Let $\seq{C_{\ga}:\ga<\gk^+}$ be a  $\square_\gk$ sequence.
This means that if $\alpha<\kappa$ then  $C_{\alpha}$ is a closed
unbounded subset of $\alpha$ 
with ordertype at most $\kappa$, 
and if $\beta$ is a limit point of $C_{\alpha}$ then
$C_{\beta}=C_{\alpha}\cap\beta$.  
It will be convenient to assume that $C_{\ga+1}=\sing{\ga}$ for all
$\ga$, that $C_{\ga}=\ga$ for limit $\ga\le\gk$, and that 
$\min(C_{\alpha})=\gk$ for all limit $\ga>\gk$.
We will write $c_{\ga,\xi}$ for the $\xi$th member of $C_{\ga}$.

The desired functions $\pi_{\ga}$ and $f_\ga$ will be defined by
writing $\ga$ as a union $\union_{\xi<\gk}A_{\ga,\xi}$ of sets
$A_{\ga,\xi}$ of size less than $\gk$:

\begin{definition}\label{def:A}
  We define $A_{\alpha,\xi}$ for $\alpha<\kappa^+$ and $\xi<\kappa$ by
  recursion on $\alpha$:
  \begin{enumerate}
  \item If $\alpha=\eta+1$ then $A_{\alpha,0}=\nothing$  and 
    $A_{\alpha,\xi}=A_{\eta,\xi}\cup\sing{\eta}$ for $0<\xi<\kappa$. 
  \item If $\alpha$ is a limit point of $\lim(C_\alpha)$ then
    $A_{\alpha,\xi}=\bigcup\set{A_{\eta,\xi}:\eta\in\lim(C_\alpha)}$.
  \item If $\alpha$ is a limit ordinal but $\lim(C_{\alpha})$ is
    bounded in $\alpha$ then  
    set $\bar\alpha=\sup(\sing{0}\cup\lim(C_\alpha))$,  and let
    $\set{\alpha_n:n<\omega}$ enumerate  
    $C_\alpha\setminus\alpha_0$ in increasing
    order.
    Thus $\bar\alpha=\alpha_0$ if $\otp(C_\alpha)>\omega$, and $\bar\alpha=0$
    otherwise.  Then 
    \begin{equation*}
      A_{\alpha,\xi}= A_{\bar\alpha,\xi}
        \cup
        \bigcup_{n<k}\left(A_{\alpha_n,\xi}\cup\sing{\alpha{_n}}\right)
    \end{equation*} 
    where  $k\le\omega$ is least such that either 
    \begin{inparaenumi}
    \item 
      $\otp(C_{\alpha_0})+k\ge\xi$,
    \item $k>0$ and $\alpha_{k-1}\notin
      A_{\alpha_k,\xi}$, or
    \item $k=\omega$.
    \end{inparaenumi}
  \end{enumerate}
\end{definition}
\begin{proposition}\label{thm:Agaxi-easy}
  \begin{enumerate}
  \item If $\xi'<\xi<\kappa$ then $A_{\alpha,\xi'}\subseteq A_{\alpha,\xi}$.
  \item If $\xi<\otp(C_{\alpha})$
    then $A_{\alpha,\xi}\subseteq c_{\alpha,\xi}$.
  \item  $\lim(C_\alpha)\cap c_{\alpha,\xi}\subset A_{\alpha,\xi}$.
  \item $\bigcup_{\xi<\kappa}A_{\alpha,\xi}=\alpha$.
  \end{enumerate}
\end{proposition}
\begin{proof}
  Each of the four clauses in this proposition is proved by induction
  on $\alpha$.  In the successor case~\ref{def:A}(1) all clauses of
  this lemma follow from the induction hypothesis applied to
  $A_{\alpha-1,\xi}$, so we only need to consider
  cases~\ref{def:A}(2,3). 

  For  
  clause~1, the induction argument follows easily from an
  inspection of the terms of the definition.   
  \smallskip{}

  In the case that $\alpha$ falls into
  case~\ref{def:A}(2), clause~2 follows immediately from the induction
  hypotheses together with the fact that $C_{\eta}=C_{\alpha}\cap\eta$
  for all $\eta\in\lim(C_\alpha)$.
  In the case that $\alpha$ falls into
  case~\ref{def:A}(3), it follows similarly by applying  the induction
  hypothesis to $A_{\bar\alpha}$ when $c_{\alpha,\xi}\leq\bar\alpha$,
  and it follows from clause~(i) in the definition of $k$ for larger
  $\xi$. 
  
    In the case that $\alpha$ falls into
  case~\ref{def:A}(2), clause~3 follows from the induction hypothesis
  in the same way as did clause~2.
  Also similarly, the induction hypothesis applied to
  $A_{\bar\alpha,\xi}$ verifies clause~3 when
  $c_{\alpha,\xi}\leq\bar\alpha$, and the definition of $k$ ensures
  that $\alpha_0\in A_{\alpha,\xi}$ when $c_{\alpha,\xi}>\alpha_0$.  
  \smallskip

  To prove clause~4 
  in the case that $\alpha$ falls into
  case~\ref{def:A}(2), we have
  $\bigcup_{\xi<\kappa}A_{\alpha,\xi}=
  \bigcup_{\alpha'\in\lim(C_{\alpha}}\bigcup_{\xi<\kappa}A_{\alpha',\xi}$
  and by the induction hypthesis
  $\alpha'=\bigcup_{\xi<\kappa}A_{\alpha',\xi}$ for all
  $\alpha'\in\lim(C_{\alpha})$. 
  In case~\ref{def:A}(3) we have
  $\bar\alpha=\bigcup_{\xi<\kappa}A_{\bar\alpha,\xi}\subseteq
  \bigcup_{\xi<\kappa}A_{\alpha,\xi}$ and
  $\alpha_n=\bigcup_{\xi<\kappa}A_{\alpha_n,\xi}$ for each $n<\omega$ 
  by the induction hypothesis.  To complete the proof it will be sufficient 
  to show that for each $n<\omega$ there is an ordinal $\xi<\kappa$ such
  that $k>n$, where $k$ is the integer used in case~\ref{def:A}(3)
  to define $A_{\alpha,\xi}$. 
  For $n=0$ this is true for $\xi=\otp(C_{\bar \alpha})+1$.   Assume
  as an induction hypothesis that there is $\xi_0$ such that $k> n$
  for $\xi\geq \xi_0$.   By the induction hypothesis on $\alpha$ there
  is $\xi_1$ such that $\alpha_{n}\in A_{\alpha_{n+1},\xi_1}$, and
  then $k> n+1$ whenever
  $\xi\geq\max(\xi_0,\xi_1,\otp(C_{\bar\alpha})+n+2)$.  
\end{proof}
The next lemma states the most important property of the sets
$A_{\alpha,\xi}$: 
\begin{lemma}\label{thm:Agaxi-cohere}
  If $\gamma\in
  A_{\alpha,\xi}\cup\lim(A_{\alpha,\xi})\cup\lim(C_{\alpha})$ then
  $A_{\gamma,\xi}=A_{\alpha,\xi}\cap\gamma$.
\end{lemma}
\begin{proof}
  Again we prove this lemma by induction on $\alpha$,and   the successor
  case~\ref{def:A}(1) is straightforward.

  When $\alpha$ falls into 
  case~\ref{def:A}(2), we first observe that if $\gamma<\gamma'$ are
  in $\lim(C_{\alpha})$ then
  $A_{\gamma,\xi}=A_{\gamma',\xi}\cap\gamma$ by the induction
  hypothesis, and it follows that
  $A_{\gamma,\xi}=A_{\alpha,\xi}\cap\gamma$ for all
  $\gamma\in\lim(C_{\alpha})$.   If $\gamma\in
  A_{\alpha,\xi}\cup\lim(A_{\alpha,\xi})$ then pick
  $\gamma'\in\lim(C_{\alpha})\setminus\gamma$.   Then by the induction
  hypothesis
  $A_{\gamma,\xi}=A_{\gamma',\xi}\cap\gamma=A_{\alpha,\xi}\cap\gamma$.  

  Now suppose $\alpha$ falls into case~\ref{def:A}(3).  Then the lemma
  holds for $\gamma\leq\bar\alpha$ by the same argument.
  For $\gamma>\bar\alpha$, note that if $k$ is as used in the
  definition of $A_{\alpha,\xi}$ then for any $n<n'<k$ we have 
  $A_{\alpha_n,\xi}=A_{\alpha_{n'},\xi}\cap\alpha_{n}=A_{\alpha,\xi}\cap\alpha_n$.
  Now for any $\gamma\in A_{\alpha,\xi}\cup\lim(A_{\alpha,\xi})$ we
  must have $\alpha_{n}\geq\gamma$ for some $n<k$, and the lemma then 
  follows in the same way as in case~\ref{def:A}(2).
\end{proof}
\begin{corollary}\label{thm:Asmall}
  If $\omega\leq\xi<\kappa$ then $\card{A_{\alpha,\xi}}\le\card\xi$.
\end{corollary}
\begin{proof}
  The proof is by induction on $\alpha$.   The only problematic case
  is \ref{def:A}(3), in which case $A_{\alpha,\xi}$ is defined as a
  union of $\card{\lim(C_\alpha)}$ many sets.
  However lemma~\ref{thm:Agaxi-cohere} and
  proposition~\ref{thm:Agaxi-easy}(2) imply that in this case
  $A_{\alpha,\xi}=\bigcup\set{A_{\eta,\nu}:\eta\in\lim(C_\alpha)\cap
    c_{\alpha,\xi}}$, a union of  $\card{\xi}$ many sets.   Since the
  induction hypothesis implies that each of these sets $A_{\eta,\xi}$
  has size at most 
  $\card{\xi}$, it follows that $\card{A_{\alpha,\xi}}\leq\card{\xi}$. 
\end{proof}

\begin{corollary}\label{thm:Alimits}
  If $\gamma\in\lim(A_{\alpha,\xi})\cap\alpha$ then
  $\xi\geq\otp(C_\gamma)$ and 
  $\gamma\in A_{\alpha,\xi+1}$.   Furthermore $\gamma\in
  A_{\alpha,\xi}$ unless $\xi=\otp(C_{\gamma})$. 
\end{corollary}
\begin{proof}
  The proof is by induction on $\alpha$.   The conclusion follows immediately
  from the induction hypothesis and lemma~\ref{thm:Agaxi-cohere}
  unless $\gamma=\sup(A_{\alpha,\xi})$.  It also follows
  easily from the induction hypothesis if $\alpha$ falls into one
  of the first two cases of Definition~\ref{def:A}, so we can assume
  that 
  $\alpha$ falls into case~\ref{def:A}(3).   If $k>0$ then
  $\sup(A_{\alpha,\xi})=\alpha_{k-1}\in A_{\alpha,\xi}$, and if $\gamma<\bar\alpha$
  then the conclusion follows from the induction hypothesis.   This
  only leaves the case $\gamma=\bar\alpha=\alpha_0$.  
  Now $\alpha_0\in\lim(A_{\alpha_0,\xi})$
  implies that $\xi\geq\otp(C_{\alpha_0})$ by
  Lemma~\ref{thm:Agaxi-easy}(2),
  and if $\xi\geq\otp(C_{\alpha_0})+1$ then $k>0$ and so  $\alpha_0\in
  A_{\alpha,\xi}$.
\end{proof}
% \marginpar{$C^*_{\alpha,\xi}=A_{\alpha,\xi}\cap C_{\alpha}$, if it is
%   ever needed.}

\begin{corollary}
  Suppose $\gamma\in\alpha\setminus C_{\alpha}$, and let
  $\bar\gamma=\min(C_{\alpha}\setminus\gamma)$.   Then $\gamma\in
  A_{\alpha,\xi}$ if and only if $\gamma\in A_{\bar\gamma,\xi}$ and
  $\bar\gamma\in A_{\alpha,\xi}$. \qed
\end{corollary}

The following corollary, giving some other useful properties of the sets
$A_{\ga,\gl}$, is easily proved using the definition and previous results:

\begin{corollary}
  \label{thm:Auseful}
  \begin{enumerate}
  \item \label{item:useful1} Suppose that  $\nu\deq\sup(A_{\ga,\gl}\cap
    A_{\ga',\gl'})<\min(\alpha,\alpha')$ and $\gl>\gl'$.
    Then $\nu\in A_{\ga,\gl}\cap
    A_{\ga',\gl'+1}$.
  \item If $\gl$ is a limit ordinal then
    $A_{\ga,\gl}=\union_{\gl'<\gl}A_{\ga,\gl'}$. 
  \item If $\alpha'<\alpha$ and $A_{\ga',\gl}\subset A_{\ga,\gl}$ then
    $A_{\ga',\gl}=A_{\alpha,\lambda}\cap\sup(A_{\alpha',\gl})$.  
%   \item If $\alpha\geq\kappa$ or $\alpha$ is a limit ordinal then
%     $A_{\alpha,\lambda}\cap\kappa=\min(\alpha,\lambda)$. \marginpar{Is
%       item~4 ever used?} 
  \end{enumerate}
  \begin{proof}
    For clause~(1), we have $\nu\in\lim(A_{\alpha',\lambda'})\cap
    \lim(A_{\alpha,\lambda})$.   Then corollary~\ref{thm:Alimits}
    implies that $\nu\in A_{\alpha',\lambda'+1}$.   Furthermore it
    implies $\otp(C_{\nu})\leq\lambda'<\lambda$, so $\nu\in
    A_{\alpha,\lambda}$ by the second sentence of
    corollary~\ref{thm:Alimits}. 

    The other clauses of corollary~\ref{thm:Auseful} are
    straightforward. 
  \end{proof}
\end{corollary}
\begin{definition}\label{def:Aga}
  \begin{enumerate}
  \item We define $f_\ga(\gl)=\otp(A_{\ga,\gl})$.
  \item We write $B_\ga$ for the set of cardinals $\gl$ which are
    $f_\ga(\gl)$-Mahlo.    Thus $\gk$ is $\ga+1$-Mahlo if and only if
    $B_\ga$ is stationary.
  \item We write $\pi_{\alpha}(\eta,\lambda)$ for the $\eta$th member of
    $A_{\alpha,\lambda}$, if $\otp(A_{\alpha,\lambda})>\eta$, and
    otherwise $\pi_{\alpha,}(\eta,\lambda)$ is undefined.
  \end{enumerate}
\end{definition}
\begin{proposition}
  \begin{inparaenumi}
  \item  $[f_\ga]_{\ns}=\ga$ for all $\alpha<\kappa^+$.    
  \item  If $\alpha'\in A_{\alpha,\lambda}$ then
    $\pi_\alpha(\eta,\lambda')=\pi_{\alpha'}(\eta,\lambda')$ for all
    $\lambda'\ge\lambda$ and $\eta<\otp(A_{\alpha',\lambda'})$.
  \item If $\alpha'\in\lim(C_{\alpha})$ then
    $\pi_{\alpha}(\eta,\lambda)=\pi_{\alpha'}(\eta,\lambda)$ for all
    $\lambda$ and all $\eta<\otp(A_{\alpha',\lambda})$.
  \end{inparaenumi}\qed{}
\end{proposition}

We will normally write $\pi_{\ga}\image X$ instead of the 
correct, but cumbersome, expression
$\pi_{\ga}\image\bigl(X^2\cap\domain(\pi_{\ga})\bigr)$. 

\begin{proposition}\label{thm:gaps}
  Suppose that $X\prec (H_{\gk^+},\vec C)$ and $\ga\in
  X\setminus\lim(X)$.  Then,  $\ga'\deq\sup(X\cap\ga)$ is a limit point of
  $C_{\ga}$, and $X\cap\lim(C_{\ga})$ is cofinal in $\ga'$.

  Hence $C_{\ga'}=C_{\ga}\cap \ga'$,
  $A_{\ga',\xi}=A_{\ga,\xi}\cap\ga'$ for every $\xi<\gk$, and
  $\pi_{\ga'}=
  \pi_{\ga}\restrict\set{(\eta,\lambda):\pi_{\alpha}(\eta,\lambda)<\alpha'}$.  
\end{proposition}
\begin{proof}
  By elementarity we have $\alpha'\in\lim(C_{\alpha})$, and a second
  application of elementarity shows that  $\lim(C_{\alpha})\cap X$ is
  cofinal in $\alpha'$. 
\end{proof}

\begin{definition}
  \begin{enumerate}
  \item
    If $\alpha<\kappa^+$ and $\vec X=\seq{X_{\alpha'}:\alpha'<\alpha}$
    is a sequence of subsets of $\kappa$ then the \emph{diagonal
      intersection} of the sequence $\vec X$ is the set 
    $\dinter_{\alpha'<\alpha}X_{\alpha'}=\set{\nu<\kappa:\forall\alpha'\in
      A_{\alpha,\nu}\;\nu\in X_{\alpha'}}$. 
  \item
    A sequence $\vec X=\seq{X_{\alpha}:\alpha<\kappa^+}$ is \emph{diagonally
    decreasing} if $X_{\alpha}\setminus\lambda\subset X_{\alpha'}$
    whenever $\alpha'\in A_{\alpha,\lambda}$.
  \item
    The sequence $\vec X$ is \emph{continuously diagonally decreasing}
    if, in addition,   $X_{\alpha}
    =\dinter_{\alpha'<\alpha}X_{\alpha'}$ whenever $\alpha$ is a limit ordinal.
  \end{enumerate}
\end{definition}
\begin{proposition}
  $\seq{B_\alpha:\alpha<\kappa^+}$ is continuously diagonally
  decreasing. \qed
\end{proposition}
%----------------------------------------------------------------
\subsection{The requirements $\II_{\ga,\gl}$ and
  $\OO_{\ga,(\gl',\gl]}$}
\label{sec:IIOO}
%----------------------------------------------------------------
As in the forcing in section~\ref{sec:oneclub} for one closed,
unbounded set, the 
conditions in $\Pkp$  will be  finite sets of requirements, ordered by
subset 
(that is, $p\le q$ if $p\supseteq q$).    The counterparts to
$\II_{\gl}$ and $\OO_{(\gl',\gl]}$ 
are relatively straightforward and are described in
definition~\ref{def:require2}; the counterparts to $\CC_M$ are more
complex and will be introduced in subsection~\ref{sec:CCMa}.
As in section~\ref{sec:oneclub}, the subscripts of the three types
of requirements are distinct and hence we can simply identify the
symbols with their subscripts.    

\begin{definition}\label{def:require2}
  \begin{enumerate}
  \item $\II_{\alpha,\lambda}$ is a requirement whenever $\alpha<\gk^+$ and
    $\lambda\in B^*_\alpha$.
  \item $\OO_{\gamma,(\eta',\eta]}$ is a requirement whenever
    $\eta'<\eta<\gk$, and either $\gamma=0$ or 
    $\gamma$ is a
    successor ordinal smaller than $\gk^+$.    
  \end{enumerate}
\end{definition}

As in the forcing in section~\ref{sec:oneclub}, the requirements
$\II_{\alpha,\lambda}$ will be used to determine the new closed
unbounded sets $D_{\alpha}$: if $G$ is a generic set then we will
define $\lambda\in D_{\alpha}$ if and only if there is $p\in G$ with
$\II_{\alpha,\lambda}\in p$.   
The definition of compatibility for these requirements will be
determined by the analogy to the forcing of section~\ref{sec:oneclub},
together with the desire that the sequence of sets $D_{\alpha}$ be
diagonally decreasing: 
The analogy with section~\ref{sec:oneclub} suggests  that
$\II_{\gamma,\lambda}$ should be incompatible with
$\OO_{\gamma,(\eta',\eta]}$ whenever $\eta'<\lambda\le\eta$, and the
desire that the sets be diagonally 
decreasing suggests that if $\gamma\in A_{\alpha,\lambda}$ then 
$\II_{\alpha,\lambda}$ should be incompatible with
$\OO_{\gamma,(\eta',\eta]}$ as well.

The desire that the sequence $(D_{\alpha}:\alpha<\kappa^+)$ be
\emph{continuously} diagonally decreasing motivates the stipulation
that the 
ordinal $\gamma$ in a requirement $\OOO$ cannot be a
nonzero limit ordinal:  
No condition should force that $\lambda\notin D_{\gamma}$, where
$\gamma$ is a nonzero limit ordinal, without  also forcing that
$\lambda\notin 
D_{\gamma'}$ for some
$\gamma'\in A_{\gamma,\lambda}$. 

\subsection{The requirements $\CC_{M,a}$}
\label{sec:CCMa}

The next three definitions give the formal definition of the
requirements $\CC_{M,a}$.   
In addition to the $\square_{\kappa}$ sequence $\vec C$, we fix a well
ordering $\triangleleft$ of $H_{\kappa^+}$, which will be used to provide
Skolem functions for that set.

\begin{definition}\label{def:model}
  As used in this section, a \model\ is a structure $M$ such that
  \begin{inparaenumi}
  \item $M\prec(H_{\kappa^+},{\in},\vec C,\triangleleft)$,
  \item $M\cap\lim(C_{\sup(M)})$ is cofinal in $M$, and
  \item $\otp(C_{\sup(M)})\notin M$.
  \end{inparaenumi}
\end{definition}

For the remainder of this section we will write $M\prec
H_{\kappa^+}$ rather than $M\prec(H_{\kappa^+},{\in},\vec
C,\triangleleft)$, leaving the predicates ${\in}$, $\vec C$ and
$\triangleleft$ to be understood.
Other predicates, when needed for the construction of particular
models, will be specified: thus if $X$ is a \model{} and
$\tau=\sup(X)$ then we may write $M\prec (X,C_{\tau})$ to indicate
that $M$ is elementary with respect to the extra predicate $C_{\tau}$
as well as the standard predicates $\in$, $\vec C$ and $\triangleleft$.

\begin{proposition}\label{thm:Minter}
  If $M$ and $M'$ are \model{s} then $M\cap M'$ is a model.
\end{proposition}
\begin{proof}
  The presence of the  well ordering $\triangleleft$ provides Skolem
  functions which ensure that an intersection of elementary
  substructures is an
  elementary substructure.   Hence $M\cap M'$
  satisfies clause~\ref{def:model}(i). 

  To verify clause~\ref{def:model}(ii), 
  set
  $\bar\alpha=\sup(M\cap M')$, and note that each of $\lim(C_{\bar\alpha})\cap
  M$ and $\lim(C_{\bar\alpha})\cap
  M'$ is  cofinal in $\bar\alpha$.   If $\bar\alpha=\sup(M)$ or
  $\bar\alpha=\sup(M')$ then this is Definition~\ref{def:model}(iii);
  otherwise it follows from Definition~\ref{def:model}(i), together, if
  $\bar\alpha$ is not in the model, with proposition~\ref{thm:gaps}. 
  Fix any
  $\gamma<\bar\alpha$ and let $\alpha\in M\setminus\gamma+1$ and $\alpha'\in
  M'\setminus\gamma+1$ be limit points of $C_{\bar\alpha}$.     Then 
  $C_{\bar\alpha}\cap\alpha=C_{\alpha}$ and 
  $C_{\bar\alpha}\cap\alpha'=C_{\alpha'}$, so the 
  least limit point of $C_{\bar\alpha}\setminus\gamma$ is also the
  least limit point of both $C_{\alpha}\setminus\gamma$ and of
  $C_{\alpha'}\setminus\gamma$, and hence is in $M\cap M'$.

  To verify clause~\ref{def:model}(iii), note that $\bar\alpha\notin
  M\cap M'$, as 
  otherwise we would have $\bar\alpha+1\in M\cap M'$.   If $\bar
  \alpha=\sup(M)$ then $\otp(C_{\bar \alpha})\notin M$ and if $\bar
  \alpha=\sup(M')$ then $\otp(C_{\bar\alpha})\notin M'$, and in either
  case $\bar\alpha\notin M\cap M'$.
  Otherwise set $\alpha=\min(M\setminus\bar\alpha)$ and
  $\alpha'=\min(M'\setminus\bar\alpha)$ and let
  $\nu=\otp(C_{\bar\alpha})$.   Then 
  $C_{\bar\alpha}=C_{\alpha}\cap\bar\alpha=C_{\alpha'}\cap\bar\alpha$, so 
  $\bar\alpha=c_{\alpha,\nu}=c_{\alpha',\nu}$.
  Thus $\nu\in M\cap M'$ would imply $\bar\alpha\in M\cap M'$.
\end{proof}
For most of this subsection, and all of the following two subsections,
we will only be considering countable \model{s}, but in
subsections~\ref{sec:M'a'} and~\ref{sec:finish} we will discuss \model{s}
$M$ of two other types: models $M$ with $\card M <\kappa$ and $M\cap
\kappa\in \kappa$
 (corresponding to the requirement
$\II_{\sup(M),\sup(M\cap\kappa)}$) and transitive models  $M$ of size
$\kappa$.  
We say that a \model{} $M$ of any of these three types is
\emph{\simple{}} if 
$\otp(C_{\sup(M)})=\sup(M\cap\kappa)$.  
We will show in subsection~\ref{sec:finish} that there are
stationarily many \simple{} models of any of these three types..

\begin{definition}
  \label{def:proxy}
  A \emph{proxy}  is a finite set of pairs $(\alpha,\lambda)$ such
  that $\lambda<\kappa$ and $\alpha$ is a limit ordinal less than
  $\kappa^{+}$.     If $a$ is a proxy then we
  write $a(\lambda)=\set{\alpha:\exists\lambda'\leq\lambda\;
  (\alpha,\lambda')\in   a}$. 
\end{definition}

\begin{definition}\label{def:requireM}
  $\CC_{M,a}$ is a requirement if
  $M$ is a countable \model{} and 
  $a$ is a proxy
  such that 
  \begin{inparaenumi}
  \item \label{item:x1}
    If $(\alpha,\lambda)\in a$ then $\lambda<\sup(M\cap\kappa)$ and
    $\alpha>\sup(M)$,
  \item\label{item:x2}
    $\pi_{\alpha}(\zeta,\lambda)\in M$ whenever  $\lambda\in M$,
    $\alpha\in a(\lambda)$ and $\zeta\in M$, 
  \item\label{item:x3}
    if $\lambda\in M$ and $\alpha\in a(\lambda)$ then 
    either $\sup(A_{\alpha,\lambda})\in
    M$ or  $M\cap A_{\alpha,\lambda}$ is cofinal in $M$, and
  \item \label{item:x4}
    if $\lambda\notin B_{\alpha}$ then 
    $\lambda\notin B_{\sup(M\cap A_{\alpha,\lambda})}$. 
  \end{inparaenumi}
\end{definition}
We will write $\CC_M$ for the requirement $\CC_{M,\nothing}$ with an
empty proxy, and we say that a requirement $\CC_M$
is \emph{\simple} if $M$ is.

Note that if $\CC_{M,a}$ is a requirement in this forcing then
$\CC_{M\cap H_\kappa}$ is a requirement in the forcing $P_B$ of
section~\ref{sec:oneclub}.   
The effect of a requirement $\CC_{M}$ in this forcing will be roughly the same
as if the requirement $\CC_{M\cap H_\kappa}\in P_{B_{\alpha}}$ were
used for each  
set $D_{\alpha}$ with $\alpha\in M$.

We will complete this subsection with some further useful observations
about the behavior of the requirements $\CC_M$ and $\CC_{M,a}$;
but first will we will  briefly explain why the proxies are needed. 
We will want to prove, for any simple countable model $M$, that
the condition $\sing{\CC_M}$ is strongly $M$-generic.
To do so we will need to define a witness function $p\mapsto p\cut
M$.   
Consider the special case $p=\sing{\CC_N,\CC_M}$, where $N$ is another
simple model with $M\cap N\in M$.   The analogy with
section~\ref{sec:oneclub} suggests trying $\sing{\CC_N,\CC_M}\cut
M=\sing{\CC_{M\cap N}}$.   The problem with this is that there may be $\xi\in
M$ such that $\xi>\sup(M\cap N)$ but $\xi\in A_{\eta,\lambda}$ for
some $\eta\in N$ and $\lambda\in M\cap N\cap\kappa$.   
In that case
any requirement $\OO_{\xi,(\lambda',\lambda]}\in M$ would be compatible with
$\CC_{M\cap N}$, however I claim that it must be incompatible with
$\CC_N$.   The reason for this deals with the need for a
function $q\mapsto q\cut N$ witnessing that $\sing{\CC_N}$ is strongly
$N$-generic.    
If $\eta$ is the least ordinal in $N$ such that $\xi\in
A_{\eta,\lambda}$ then 
the requirement $\OO_{\xi,(\lambda',\lambda]}$ is incompatible with
$\II_{\eta,\lambda}$, but compatible with $\II_{\eta',\lambda}$ for
any $\eta'\in N\cap\eta$.  
The same should be true of 
$\sing{\OO_{\xi,(\lambda',\lambda]},\CC_N}\cut N$, and the only condition in
$N$ which would have this effect would seem to be 
$\sing{\OO_{\eta,(\lambda',\lambda]}}$.  However 
$\OO_{\eta,(\lambda',\lambda]}$ is not a requirement
since $\eta$ is a non-zero limit  ordinal, and 
hence $\sing{\OO_{\eta,(\lambda',\lambda]}}$ is not a condition.
Since
there is no good choice for
$\sing{\OO_{\xi,(\lambda',\lambda]},\CC_N}\cut N$, 
our definition of the forcing will have to specify that
$\OO_{\xi,(\lambda',\lambda]}$ is incompatible with $\CC_{N}$.

Thus the correct choice of $\sing{\CC_N,\CC_M}\cut M$ must be a condition
which is incompatible
with every requirement $\OO_{\xi,(\lambda',\lambda]}$ as in the last
paragraph.     This 
will be accomplished by setting $\sing{\CC_N}\cut M=\sing{\CC_{N\cap
    M,b}}$ where $b$ is a proxy chosen so that for any requirement
$\OO_{\xi,(\lambda',\lambda]}\in M$ as in the last paragraph there is
some $\eta'\in b(\lambda)$ such that $\xi\in A_{\eta',\lambda}$.
The construction of $\CC_N\cut M$ will be given
in section~\ref{sec:M'a'}, with the construction of the proxy $b$
given in lemma~\ref{thm:CinM}.

\begin{proposition}\label{thm:CaispaM}
  If $\CC_{M,a}$ is a requirement and $\ga\in M\cup\lim(M)$ then 
  $M\cap\ga=\pi_{\ga}\image (M\cap\kappa)$ and $M\cap
  C_{\alpha}=\set{c_{\alpha,\nu}:\nu\in M\cap\otp(C_\alpha)}$. 
\end{proposition}
\begin{proof}
  If $\ga\in M$ then 
  the proposition is immediate since $\pi_\alpha$ and $C_\alpha$ are
  in $M$.  
  If $\ga\in\lim(M)\setminus M$ and $\ga<\sup(M)$ then 
  set $\ga'\deq\min(M\setminus\ga)$.    
  Then $M\cap\alpha'=\pi_{\alpha'}\image (M\cap\kappa)$ by the
  previous sentence, 
  and proposition~\ref{thm:gaps} implies that this is equal to
  $\pi_{\alpha}\image (M\cap\kappa)$.
  The second clause follows from proposition~\ref{thm:gaps} and the
  observation that $M\cap 
  C_\alpha=M\cap C_{\alpha'}$.

  Finally, if $\ga=\sup(M)$ then
  $M\cap\alpha=\union\set{M\cap\nu:\nu\in\lim(C_{\alpha}}=
  \union\set{\pi_{\nu}\image (M\cap\kappa):\nu\in\lim(C_{\alpha)}}
  =\pi_{\alpha}\image(M\cap\kappa)$. 
\end{proof}

Notice that proposition~\ref{thm:CaispaM} implies in particular that
the set of ordinals of 
any requirement $\CC_M$ is determined by $M\cap
\kappa$ together with  $\sup(M)$.
It follows, by using the well ordering $\triangleleft$ specified at
the beginning of section~\ref{sec:CCMa}, that these determine
$\CC_{M}$ itself.

\begin{corollary}\label{thm:Macap}
  If $\CC_{M,a}$ and $\CC_{N,b}$ are requirements then $M\cap
  N=\pi_{\bar\ga}\image(M\cap N\cap\kappa)\subset
  A_{\bar\ga,\sup(M\cap N\cap\kappa)}$ 
  where $\bar\ga=\sup(M\cap N)$. \qed
\end{corollary}

\begin{corollary}\label{thm:shadow1}
  Suppose that $\CC_{M,a}$ and $\CC_{N,b}$ are requirements with
  $M\cap N\cap\kappa\in M$.   Then $M\cap N\in M$.
\end{corollary}
\begin{proof}
  First note that
  $\bar\ga\deq\sup(M\cap N)<\sup(M)$, since 
  $\otp(C_{\bar\ga})\in\lim(M\cap N\cap\kappa)\subset M$ while 
  $\otp(C_{\sup(M)})\notin M$. 
  Then $M\cap
  N=\pi_{\bar\ga}\image (M\cap
  N\cap\kappa)=\pi_{\min(M\setminus\bar\alpha)}\image(M\cap 
  N\cap\kappa)\in M$. 
\end{proof}

\begin{definition}
  If $\CC_{M,a}$ is a requirement then we write $A_{M,a,\lambda}$ for
  $\bigcup\set{A_{\alpha,\lambda}:\alpha\in
  M\cup a(\lambda)}$. 
\end{definition}

The following observation will be used frequently.

\begin{lemma}\label{thm:ASminBig}
  Suppose that $\CC_{M,a}$ and $\CC_{N,b}$ are requirements
  such that  $N\cap\sup(M\cap N\cap\kappa)\subset M$.   Then 
  $A_{M,a,\gl}\cap N\subset M$ for all 
  $\gl<\sup(M\cap N\cap\kappa)$.
\end{lemma}
\begin{proof}
  Suppose that $\eta'\in N\cap A_{\eta,\gl}$, where $\gl<\sup(M\cap
  N\cap\kappa)$ 
  and  $\eta\in M\cup a(\lambda)$.
  By increasing $\gl$ if necessary, we can assume that $\gl\in M\cap N$.
  Then $\gg\deq\otp(A_{\eta',\gl})\in N\cap\gl^+\subset M$, since
  $\card{A_{\eta',\gl}}\le\gl$ and both $N\cap\kappa$ and
  $M\cap\kappa$ are closed under 
  cardinal successor.   However $A_{\eta',\gl}=A_{\eta,\gl}\cap\eta'$ so
  $\eta'=\pi_{\eta}(\gamma,\lambda)\in M$. 
\end{proof}
\begin{lemma}\label{thm:shadow2} 
  Suppose that $M\cap N\cap\kappa\subseteq M$.
  Furthermore, suppose that  $\gl\in 
  M\cap N\cap\kappa$ and  $\ga\in N\cup b(\lambda)$, and let  
  $\bar\alpha\deq\sup(M\cap N)$.
  \begin{enumerate}
  \item \label{item:shadow2a}
    If 
    $A_{\ga,\gl}\cap\bar\ga$ is bounded in $\bar\ga$ then
    $\sup(A_{\ga,\gl}\cap M)\in M\cap N$.
  \item\label{item:shadow2b}
    If $A_{\ga,\gl}\cap\bar\ga$ is unbounded in $\bar\ga$ and
    $\bar\ga<\sup(N)$ then
    $\ga'\deq\min(N\setminus\bar\ga)\in\lim(A_{\ga,\gl})$ and
    $A_{\ga,\gl}\cap M\subset\ga'$.
  \end{enumerate}
\end{lemma}
\begin{proof}
  For clause~1, suppose that $A_{\ga,\gl}$ is
  bounded in $\bar\alpha$ and set $\gg\deq\sup(A_{\ga,\gl}\cap A_{\bar\ga,\gd})$ 
  where $\gd\deq\sup(M\cap N\cap\kappa)$.
  Since $\gl<\gd$,
  Corollary~\ref{thm:Auseful}(1) implies that $\gg\in
  A_{\ga,\gl+1}\cap A_{\bar\ga,\gd}$.   
  Also, since $\bar\ga$ and $\gd$ are limit ordinals there are
  $\bar\ga'\in (\lim(C_{\bar\ga})\cap M\cap N)\setminus\gg$ and $\gd'\in
  M\cap N\cap\kappa$ so that $\gg\in A_{\bar\ga',\gd'}$.   Then
  $\gg=\sup(A_{\ga,\gl}\cap A_{\bar\ga',\gd'})\in N$, and it follows
  by Lemma~\ref{thm:ASminBig} that $\gg\in
  M$ as well.   Thus $\gamma\in M\cap N$, and  it
  remains to   
  show that $M\cap A_{\ga,\gl}\subset\gg+1$.   Suppose to the
  contrary that there is $\eta>\gg$ in $A_{\ga,\gl}\cap M$, and set
  $\gamma'\deq\min(A_{\alpha,\lambda}\setminus\gamma+1)$.    Then
  $\gamma'\in N$, and
  $\gamma'=\min\left(A_{\eta,\lambda}\cup\sing{\eta}\setminus\gamma+1\right)\in
  M$.   Thus $\gamma'>\gamma$ is in $M\cap N$, contradicting the
  choice of $\gamma$.

  Now suppose that the hypothesis  to clause~2 holds.
  First we show that we can assume that $\alpha\in N$:
  Otherwise   $\alpha\in b(\lambda)$, but
  in that case 
  clause~\ref{def:requireM}(\ref{item:x3}) implies that  
  $\alpha''\deq \sup(N\cap A_{\alpha,\lambda})$ is either a member of
  $N$ or else is equal to $\sup(N)$.   If $\alpha''\in N$ then it will
  be sufficient to show that clause~2 holds with $\alpha''$ in place
  of $\alpha$, and if $\alpha''=\sup(N)$ it will be sufficient to
  show that clause~2 holds for any member of
  $(N\setminus\alpha')\cap A_{\alpha,\lambda}$ in place of $\alpha$.

  Now $\sup(\alpha'\cap A_{\alpha,\lambda})\in N$ because  $\alpha'$,
  $\alpha$ and $\lambda$ are in $N$.   Since
  $\bar\alpha\leq\sup(\alpha'\cap A_{\alpha,\lambda})\leq\alpha'=\min(N\setminus\bar\alpha)$ it follows that 
  $\alpha'\in\lim(A_{\alpha,\lambda})$.

  It remains to show that $A_{\alpha,\lambda}\cap M\subseteq\alpha'$.
  Suppose to the contrary that there is some ordinal
  $\eta\in A_{\alpha,\lambda}\cap
  M\setminus\alpha'$.   Then $\alpha'\in M$, either because
  $\eta=\alpha'$ or because
  $\eta>\alpha'$, in which case  $\alpha'\in
  \lim(A_{\alpha,\lambda})\cap\eta=\lim(A_{\eta,\lambda})$, so 
  $\alpha'\in A_{\eta,\lambda+1}$ by corollary~\ref{thm:Alimits} and thus
  $\alpha'\in M$ by lemma~\ref{thm:ASminBig}.   However $\alpha'\in 
  M\cap N$ would imply $\alpha'+1\in M\cap N$, contradicting the
  fact that $\alpha'\ge\bar\alpha=\sup(M\cap N)$.
\end{proof}

%----------------------------------------------------------------

\subsection{Definition of the forcing $\Pkp$}
\label{sec:Pstar}
The definition of the forcing $\Pkp$, given in
definitions~\ref{def:Pscompat} and~\ref{def:Ps} below, 
is very nearly a word for word copy---with the mechanical addition of
the extra subscripts---of definitions~\ref{def:Pcompat}
and~\ref{def:PB} of the forcing $P_{B}$ in section~\ref{sec:oneclub}.
The most significant changes appear in clauses~1 and~2.  The change in
clause~1, which was alluded to at the end of
subsection~\ref{sec:IIOO}, is needed to account for the added
subscripts $\alpha$ in $\II_{\alpha,\lambda}$ and $\gamma$ in $\OOO$.
The change in clause~2, using $M[a]$ in place of $M$, was alluded to in
subsection~\ref{sec:CCMa} and is needed to take account of proxies.

A 
more subtle change comes in the definition of an $M$-fence:  
if $\II_{\alpha,\lambda}$ is an $M$-fence then $\lambda$ is required
to be a member of $M$, but $\alpha$ is not.   We will see in
subsection~\ref{sec:M'a'} that if $M$ is \simple{} then $\alpha$ can also
be taken to be a member of~$M$.
For this reason we will have strongly generic conditions for
simple models, but only for simple models.

Except for these changes, the definition is essentially a word for
word copy of 
definitions~\ref{def:Pcompat} and~\ref{def:PB} with the additional
subscripts mechanically added to the requirements.

\begin{definition}\label{def:Ma}
  We write $M[a]$ for $\set{\II_{\alpha,\lambda}:\lambda\in 
  M\cap B^*_{\alpha}\And\alpha\in M\cup a(\lambda)}$.
\end{definition}

\begin{definition}\label{def:Pscompat}
  \begin{enumerate}
  \item \label{item:sOI}
    Two requirements $\OO_{\gamma,(\eta',\eta]}$ and
    $\II_{\alpha,\lambda}$ are 
    \emph{incompatible} if $\eta'<\lambda\le\eta$ and
    $\gamma\in A_{\alpha+1,\lambda}$; otherwise they are
    \emph{compatible}. 
  \item\label{item:sOM}
    Two requirements $\OO_{\gamma,(\eta',\eta]}$ and $\CC_{M,a}$ are
    \emph{compatible} if either $\OO_{\gamma,(\eta',\eta]}\in M$ or every
    requirement $\II_{\alpha,\lambda}\in M[a]$ is compatible with
    $\OO_{\gamma,(\eta',\eta]}$. 
  \item\label{item:sMI}
    \begin{enumerate}
    \item \label{item:sMIf}
      An \emph{$M$-fence for a requirement $\II_{\alpha,\lambda}$} is
      a requirement 
      $\II_{\alpha',\lambda'}$ with $\lambda'\in M$ such that any
      requirement 
      $\OO_{\gamma,(\eta',\eta]}$ in $M$ incompatible with
      $\II_{\alpha,\lambda}$ is 
      also incompatible with $\II_{\alpha',\lambda'}$.
    \item\label{item:sMIc}
      Two requirements $\CC_{M,a}$ and $\II_{\alpha,\lambda}$ are
      \emph{compatible} if 
      either $\lambda\ge\sup(M\cap\kappa)$ or there exists a $M$-fence
      for $\II_{\alpha,\lambda}$. 
    \end{enumerate}
  \item\label{item:sMN}
    \begin{enumerate}
    \item\label{item:sMNf}
      An \emph{$M$-fence} for a requirement $\CC_{N,b}$ is a finite 
      set $x$ of requirements $\II_{\alpha,\lambda}$, with 
      $\lambda\in M\cap B_\alpha$, with the following property: 
      Suppose that 
      $\OO_{\gamma,(\eta',\eta]}\in M$ is a requirement such that 
      $\eta\geq\sup(M\cap N\cap\kappa)$, 
      $\eta'\geq\sup(M\cap N\cap\kappa)$ 
      if  $M\cap N\in M$, and $\OO_{\gamma,(\eta',\eta]}$ is 
      incompatible with $\CC_{N,b}$.   Then there is some requirement
      $\II_{\alpha,\lambda}\in x$ which is incompatible with
      $\OO_{\gamma,(\eta',\eta]}$.  
    \item\label{item:sMNc}
      A \model{} $M$ \emph{is fenced from}
        a requirement $\CC_{N,b}$ if 
      \begin{inparaenumi}
      \item \label{item:holds}
        either $M\cap N\cap H_\kappa\in M$ or $M\cap N\cap H_\kappa=M\cap
        H_{\sup(M\cap N\cap\kappa)}$, and 
      \item
        there is a $M$-fence for $\CC_{N,b}$.
      \end{inparaenumi}
    \item 
      Two requirements $\CC_{M,a}$ and $\CC_{N,b}$ are
      \emph{compatible} if $M$ is fenced from $\CC_{N,b}$ and $N$ is fenced from
      $\CC_{M,a}$. 
    \end{enumerate}
  \end{enumerate}
\end{definition}

\begin{definition}
  \label{def:Ps}
  A condition $p$ in the forcing $\Pkp$ is a finite set of requirements
  such that each pair of requirements in $p$ is compatible.   The
  order on $\Pkp$ is reverse inclusion: $p'\le p$ if $p'\supseteq p$.
\end{definition}

Although this forcing is somewhat more complicated than the forcing
$P_B$, our exposition will 
parallel the exposition in section~\ref{sec:oneclub}. 
Like $P_B$, the forcing $P^*$ is not separative and we will
write $p'\lex p$ if $p'\forces p\in\dot G$ and $p\eqx p'$ if $p\lex
p'$ and $p'\lex p$.   In addition we introduce the following notation
for a special case of the failure of separation:

\begin{definition}\label{def:inx}
  We say that $\II_{\alpha,\gl}\inx p$ if 
  there is $\II_{\alpha',\gl}\in p$ such that either
  $\alpha\in A_{\alpha'+1,\lambda}$ or $\alpha$ is a limit ordinal and  
  $A_{\alpha,\gl}$ is a subset (and hence an initial segment) of
  $A_{\alpha',\gl}$.
\end{definition}

\begin{proposition}\label{thm:inx}
  If $p\in\Pkp$ and $\II_{\alpha,\lambda}\inx p$ then
  $p\cup\sing{\II_{\alpha,\lambda}}\in\Pkp$.    Hence
  $p\cup\sing{\II_{\alpha,\lambda}}\eqx p$.
%  $p\forces\II_{\alpha,\lambda}\in\bigcup\dot G$.
\end{proposition}
\begin{proof}
  Let $\II_{\alpha',\lambda}\in p$ witness that
  $\II_{\alpha,\lambda}\inx p$.    Then any requirement $\OOO$ which is
  incompatible with $\II_{\alpha,\lambda}$ is also incompatible with
  $\II_{\alpha',\lambda}$, and   it follows that $\II_{\alpha,\lambda}$ is
  compatible with any requirement $\OOO\in p$.   In addition, if $M$
  is a \model{} then any $M$-fence for $\II_{\alpha',\lambda}$ is also a $M$-fence
  for $\II_{\alpha,\lambda}$, and it follows that any requirement
  $\CC_{M,a}\in p$ is compatible with $\II_{\alpha,\lambda}$.   Hence
  $p\cup\sing{\II_{\alpha,\lambda}}\in\Pkp$.

  To see that $p\forces\II_{\alpha,\lambda}\in\bigcup\dot G$, note
  that the
  first paragraph implies that
  $q\cup\sing{\II_{\alpha,\lambda}}\in\Pkp$ for any condition $q\le p$.
\end{proof}

\begin{proposition}
  If $\CC_{M,a}$ is a requirement then any requirement $\RR\in M$ is
  compatible with $\CC_{M,a}$.
  Thus $p\cup\sing{\CC_{M,a}}\in\Pkp$ for any   $p\in M\cap\Pkp$. 
\end{proposition}
\begin{proof}
  Any requirement $\II_{\alpha,\lambda}\in M$ is its own fence for
  compatibility with $\CC_{M,a}$, and any requirement
  $\OO_{\alpha,(\eta',\eta]}\in M$ is compatible with $\CC_{M,a}$.
  If $\CC_{N,b}\in M$ then, since $M\cap N=N$, the empty set
  $\emptyset$ is both a 
  $M$-fence for $\CC_{N,b}$ and an $N$-fence for $\CC_{M,a}$.
\end{proof}

Unlike the case in the forcing  $P_B$, the fences specified in
definition~\ref{def:Pscompat} are
not unique. 
In the next two lemmas we will give an alternate characterization of
compatibility, and show that if any $M$-fence exists then 
there is a unique minimal $M$-fence: 

\begin{comment}
  If $A\cap
  M$ is bounded in $M$, then 
  we write $\sup^{M,a}(A)$ for $\min((M\cup
  a(\lambda))\setminus\sup(M\cap A)$. 
  \footnote{==== This was originally $\sup^{M,a}(A)$.   I think that
  this is correct in that we don't need the proxies here.}
\end{comment}

\begin{proposition}\label{thm:minimal}
  Suppose that $M$ is a \model{} and $\II_{\alpha,\lambda}$ is a 
  requirement with $\lambda<\sup(M\cap\kappa)$.   Set
  \begin{equation*}
  \lambda'=\min(M\setminus\lambda)\qquad\text{and}\qquad
  \alpha'=\sup(\set{\gamma+1:\gamma+1\in M\cap A_{\alpha+1,\lambda}}).
  \end{equation*}
%%   \begin{equation*}
%%     \alpha' =
%%     \begin{cases}
%%       \alpha&\text{if $\alpha\in M$ and $\alpha$ is a successor,}\\
%%       \sup(\set{\gamma+1:\gamma+1\in A_{\alpha+1,\lambda}})&\text{otherwise.}
%%     \end{cases}
%%   \end{equation*}
  Then $\II_{\alpha,\lambda}$ is compatible with $\CC_{M,a}$ if and
  only if $\lambda'\in B^*_{\alpha'}$.

  Furthermore, in this  case
  $\II_{\alpha',\lambda'}$ is a $M$-fence for 
  $\II_{\alpha,\lambda}$ which is 
  minimal in the sense that  if $\II_{\alpha'',\lambda''}$ is any
  other $M$-fence for $\II_{\alpha,\lambda}$ then 
  \begin{inparaenumi}
  \item $\lambda''=\lambda'$,
  \item $\alpha'\leq\alpha''$, and
  \item $\II_{\alpha',\lambda'}\inx\sing{\II_{\alpha'',\lambda''}}$.     
  \end{inparaenumi}
\end{proposition}
Notice that  $\alpha'=\alpha$ if $\alpha=0$ or $\alpha$
is a successor ordinal in $M$.
We will call the fence $\II_{\alpha',\lambda'}$ of
proposition~\ref{thm:minimal}
the \emph{minimal $M$-fence} for $\II_{\alpha,\lambda}$.

\begin{proof}
  First, suppose that $\lambda'\in B^*_{\alpha'}$, so that 
  $\II_{\alpha',\lambda'}$ is a requirement.
  If $\OOO$ is a requirement in $M$ which is incompatible with
  $\II_{\alpha,\lambda}$ then, because $\gamma\in
  A_{\alpha+1,\lambda}$ and $\gamma$ cannot be a limit
  ordinal, the choice of
  $\alpha'$ ensures that $\gamma\in A_{\alpha'+1,\lambda}$.
  Also $\eta'<\sup(M\cap\lambda)\le\lambda\le\lambda'\le\eta$, so
  $\OOO$ is incompatible with $\II_{\alpha',\lambda'}$.
  It follows that $\II_{\alpha',\lambda'}$ is a
  $M$-fence for $\II_{\alpha,\lambda}$, and hence
  $\II_{\alpha,\lambda}$ is compatible with $\CC_{M,a}$.
  \smallskip
  
  For the other direction, suppose that $\II_{\alpha,\lambda}$ is
  compatible with $\CC_{M,a}$ and let  $\II_{\alpha'',\lambda''}$ be an
  arbitrary $M$-fence for $\II_{\alpha,\lambda}$.
  First we observe that $\lambda''=\lambda'$;  otherwise pick $\eta\in
  M\cap\lambda$ such that $\eta>\lambda''$ if $\lambda''<\lambda$.
  Then the 
  requirement $\OO_{0,(\eta,\lambda']}$ is incompatible with
  $\II_{\alpha,\lambda}$ but is compatible with
  $\II_{\alpha'',\lambda''}$.

  If $\alpha'$ is a successor ordinal then it must be a member of
  $M$.   
  In that case $\OO_{\alpha',(\eta,\lambda']}$ is a requirement in $M$
  which is incompatible with $\II_{\alpha,\lambda}$ and hence must be
  incompatible with $\II_{\alpha'',\lambda'}$, and it follows that
  $\alpha'\in A_{\alpha''+1,\lambda}$.
  
  If $\alpha'$ is a limit ordinal then let  $S$ be the set of ordinals
  $\gamma+1$ such that 
  $\OO_{\gamma+1,(\eta,\lambda']}$ is in 
  $M$ and incompatible with $\II_{\alpha,\lambda}$. 
  Then each ordinal $\gamma+1\in S$ must be a member of
  $A_{\alpha''+1,\lambda}$.
  Since $S$ is cofinal in $\alpha'$, it follows that  $\alpha'$ is a
  limit point of 
  of $A_{\alpha'',\lambda}$, and hence $\alpha'\le\alpha''$
  and $\II_{\alpha',\lambda'}\inx\sing{\II_{\alpha'',\lambda'}}$. 
\end{proof}

% =================================================

% The next two results are probably already in here somewhere.   I'm
% putting them here now for the record.    If they're not already here,
% they probably need revision and moving (probably they go before the
% last two.

% ================

\begin{proposition}\label{thm:min-is-enough}
  Suppose that $M$ is a model, $\alpha,\lambda\in M$, and $\gamma\in
  A_{\alpha,\lambda}\setminus M$.   Then $\gamma\in
  A_{\min(M\setminus\gamma),\lambda}$. 
\end{proposition}
\begin{proof}
  Set $\alpha'=\min(M\setminus\gamma)$.  Then
  $\sup(A_{\alpha,\lambda}\cap\alpha'\in M$, since $\alpha,\lambda$
  and $\alpha'$ are.   But
  $\gamma\leq\sup(A_{\alpha,\lambda}\cap\alpha'\leq\alpha'$, and since
  $\alpha'=\min(M\setminus\gamma)$ it follows that
  $\alpha'=\sup(A_{\alpha,\lambda}\cap\alpha')\in\lim(A_{\alpha,\lambda})$.
  Thus Lemma~\ref{thm:Agaxi-cohere} implies that
  $A_{\alpha',\lambda}=A_{\alpha,\lambda}\cap\alpha'$, so $\gamma\in
  A_{\alpha',\lambda}$. 
\end{proof}
\begin{lemma}\label{thm:abar-is-enough}
  Let $M$ and $N$ be models which satisfy
  $\lim(M)\cap\lim(N)=\lim(M\cap N)$, and 
  set $\bar\alpha=\sup(M\cap N)$ and $\delta=\sup(M\cap
  N\cap\kappa)$. 
  Then for any $\sing{\alpha,\lambda}\subset N$ with
  $\delta\leq\lambda<\kappa$ we have $A_{\alpha,\lambda}\cap
  M\cap\bar\alpha\subseteq A_{\bar\alpha,\lambda}$. 
\end{lemma}
\begin{proof}
  Fix $\gamma\in A_{\alpha,\lambda}\cap M\cap\bar\alpha$.
  By proposition~\ref{thm:min-is-enough} we can assume
  that $\alpha=\min(N\setminus\gamma)<\bar\alpha$. 
  We will show, by induction on $\nu$, that $\gamma\in
  A_{\nu,\lambda}$ for all $\nu\geq\alpha$ in 
  $(M\cap N\cap\bar\alpha)$.  Since $\lim(C_{\bar\alpha})\cap (M\cap
  N)$ is cofinal in $\bar\alpha$ it will follow that $\eta\in
  A_{\bar\alpha,\lambda}$. 
  
  Fix such an ordinal $\nu$.
  If $\alpha\in\lim(C_{\nu})$ then $\gamma\in
  A_{\alpha,\lambda}=A_{\nu,\lambda}\cap\alpha$, so we can assume that
  $\alpha\notin\lim(C_{\nu})$.  

  Set $\nu'=\min(C_{\nu}\setminus\alpha)$, so $\nu'\in N$.
  Since $\nu'\notin\lim(C_{\nu})$ there is $\nu''\in
  C_{\nu}\cup\sing{0}$ such that 
  $\nu'=\min(C_{\nu}\setminus\nu''+1)$.
  Then $\nu''\in N$ and it follows by the
  minimality of $\nu'$ that
  $\nu''<\alpha\leq\nu'$.  Since $\alpha=\min(N\setminus\gamma)$
  it follows that $\nu''<\gamma\leq\alpha$.
  This implies that
  $\nu'=\min(C_{\nu}\setminus\gamma)\in M$, so $\nu'\in M\cap N$
  and  the induction hypothesis implies that $\gamma\in
  A_{\nu',\lambda}$.
  Furthermore, since $\nu'\in M\cap N$ the least ordinal  $\lambda'$
  such that $\nu'\in 
  A_{\nu,\lambda'}$ is also in $M\cap N$, so 
  $\lambda'<\delta$ and hence $\nu'\in A_{\nu,\delta}\subseteq
  A_{\nu,\lambda}$.   Thus 
  $\gamma\in A_{\nu',\lambda}=A_{\nu,\lambda}\cap\nu'$.  
\end{proof}

We now consider the compatibility of requirements $\CC_{M}$ and
$\CC_{N}$.   It is easy to see that if $\CC_{M}$ and $\CC_{N}$ are
compatible then $\CC_{M\cap H_{\kappa}}$ and $\CC_{N\cap H_{\kappa}}$
are compatible in the forcing $P_{B_0}$ of section~\ref{sec:oneclub}.
In particular $M\cap\kappa$ and $N\cap\kappa$ fall into the pattern of
figure~\ref{fig:PB}: a common  initial segment which is followed by a
finite alternating sequence of disjoint intervals.
For pairs $M$ and $N$ which satisfy
Definition~\ref{def:Pscompat}(\ref{item:holds}), so that $M\cap N\cap
\kappa$ is an initial segment of at least one of $M$ and $N$, 
this can be
concisely expressed by the  statement
$\lim(M\cap\kappa)\cap\lim(N\cap\kappa)=\lim(M\cap N\cap\kappa)$.
The following proposition shows that this  equality also holds above
$\kappa$: 
\begin{proposition}\label{thm:limMcapN}
  If $M$ and $N$ are countable \model{s} such that
  $\lim(M\cap\kappa)\cap\lim(N\cap\kappa)=\lim(M\cap N\cap\kappa)$
  then $\lim(M)\cap\lim(N)=\lim(M\cap N)$.
\end{proposition}
\begin{proof}
  Suppose $\alpha\in \lim(M)\cap\lim(N)$.   Then 
  Proposition~\ref{thm:CaispaM} implies that $M\cap
  C_\alpha=\set{c_{\alpha,\nu}:\nu\in M\cap \otp(C_\alpha)}$ and $N\cap
  C_\alpha=\set{c_{\alpha,\nu}:\nu\in N\cap \otp(C_\alpha)}$.
  Since $C_\alpha\cap M$ is a cofinal subset of $M$ and
  $C_\alpha\cap N$ is a cofinal subset of $N$, it follows that
  $\otp(C_\alpha)\in\lim(M)\cap\lim(N)\cap\kappa=\lim(M\cap
  N\cap\kappa)$.   Thus $\set{c_{\alpha,\nu}:\nu\in M\cap
  N\cap\otp(C_\nu)}\subset M\cap N$ is cofinal in $\alpha$.
\end{proof}
\begin{lemma}\label{thm:minMfenceN}
  Suppose $\CC_{N,b}$ is a requirement,  $M$ is a countable
  \model{}, and the models $M$ and $N$ satisfy
  Definition~\ref{def:Pscompat}(\ref{item:holds}).  
  Let $y$ be the set of requirements
  $\II_{\alpha,\lambda}$ such that 
  \begin{enumerate}
  \item $\sup(M\cap N\cap\kappa)\le\lambda=\sup(N\cap\lambda')<\lambda'$ for
    some  
    $\lambda'\in M\cap\kappa$, and
  \item either
    \begin{inparaenumi}
    \item  $\alpha\in b(\lambda)$, 
    \item
      $\alpha=\sup(N\cap M)$, or 
    \item 
      $\alpha=\min(N\setminus\alpha')>\alpha'$ for 
      some $\alpha'\in M\setminus\sup(M\cap N)$.
    \end{inparaenumi}
  \end{enumerate}
   Then there is a $M$-fence for $\CC_{N,b}$ if and only 
   $\lim(M\cap N)=\lim(M)\cap\lim(N)$ and  each of
   the requirements 
   in $y$ is compatible with $\CC_{M}$.  
   
   Furthermore, in this case let $x$ be the set of minimal 
   $M$-fences  for requirements in $y$. Then $x$ is a 
   $M$-fence for $\CC_{N,b}$, and $x$ is minimal in the sense that if
   $x'$ is any other $M$-fence for $\CC_{N,b}$ then
   $\II_{\alpha,\lambda}\inx x'$ for any $\II_{\alpha,\lambda}\in
   x$. 
\end{lemma}

\begin{proof}
  Note that every member $\II_{\alpha,\lambda}\in y$ is a
  requirement since $\cof(\lambda)=\omega$ and hence $\lambda\in
  B^*_{\alpha}$.    Let us say that a requirement $\OOO\in M$
  \emph{clashes} with $\CC_N$ if it is incompatible with $\CC_N$,
  $\eta\geq\sup(M\cap N\cap\kappa)$, and 
  $\eta'\geq\sup(M\cap N\cap\kappa)$  if $M\cap N\in M$.
  Thus a $M$-fence for $\CC_N$ is a finite set $x$ of requirements
  $\II_{\alpha,\lambda}$ such that  $\lambda\in M$ and any requirement
  $\OOO\in M$  
  which clashes with $\CC_N$ is incompatible with some member of $x$.

  We begin by showing that every requirement
  $\OOO\in M$ 
  which clashes with $\CC_N$ is incompatible with some member of $y$.
  To this end let $\OOO$ be a requirement in $M$ which clashes with
  $\CC_N$,
  and let $\II_{\alpha_{0},\lambda_{0}}$ be a
  requirement in $N[b]$ which is incompatible with
  $\OOO$.

  Set $\lambda=\sup(N\cap\eta)$.  Then $\II_{\alpha_{0},\lambda}$ is
  a requirement since $\cof(\lambda)=\omega$, and $\gamma\in
  A_{\lambda_0,\alpha_0}\subseteq A_{\lambda,\alpha_0}$ so
  $\OOO$
  is incompatible with $\II_{\alpha_0,\lambda}$.

  If
  $\alpha_0\in b(\lambda)$ then $\II_{\alpha_0,\lambda}\in y$.
  If $\sup(M\cap N)<\alpha_0\in N$ then set
  $\alpha=\min(N\setminus\gamma)$.   Then $\gamma\in
  A_{\alpha,\lambda_0}$ by proposition~\ref{thm:min-is-enough}, and
  since $A_{\alpha,\lambda_0}\subseteq A_{\alpha,\lambda}$ it follows
  that  $\OOO$ is
  incompatible with $\II_{\alpha,\lambda}\in y$.

  Thus we can assume that $\gamma<\sup(M\cap N)$.  If
  $\lambda_0\geq\delta\deq\sup(M\cap N\cap\kappa)$ then
  Lemma~\ref{thm:abar-is-enough} 
  implies that $\gamma\in A_{\bar\alpha,\lambda_0}\subseteq
  A_{\bar\alpha,\lambda}$, where $\bar\alpha=\sup(M\cap N)$, and hence
  $\OOO$ is incompatible with $\II_{\bar\alpha,\lambda}\in y$.

  The only remaining case has $\gamma<\sup(M\cap N)$ and
  $\lambda_0<\delta$.  Since
  $\lambda_0\in(\eta',\eta]$, this implies that $\eta'<\delta$ and by
  Definition~\ref{def:Pscompat} we must have $M\cap N\notin M$, so
  $M\cap\sup(M\cap N\cap\kappa)\subset N$.  By Lemma~\ref{thm:ASminBig}
  (with $M$ and $N$ switched) it follows that
  $\gamma\in A_{\alpha_0,\lambda_0}\cap M\subseteq A_{N,b,\lambda_0}\cap
  M\subset N$.   Thus $\gamma\in M\cap N\subset
  A_{\bar\alpha,\delta}\subseteq A_{\bar\alpha,\lambda}$, so again
  $\OOO$ is incompatible with 
  $\II_{\bar\alpha,\lambda}\in y$. 
  \smallskip

  This completes the proof that any requirement $\OOO\in M$
  which clashes with $\CC_N$ is incompatible with some requirement
  $\II_{\alpha,\lambda}\in y$.   Now suppose that each requirement
  $\II_{\alpha,\lambda}$ in
  $y$ is compatible with $\CC_M$, and let $x$ be the set of minimal
  $M$-fences for members of $y$.   Than any requirement $\OOO\in M$
  which clashes with $\CC_N$ is incompatible with some member
  $\II_{\alpha,\lambda}$ of $y$ and hence with its minimal $M$-fence
  $\II_{\alpha',\lambda'}\in x$.   If, in addition,
  $\lim(M)\cap\lim(N)=\lim(M\cap N)$ then $y$, and hence $x$, is
  finite:
  If $y$ were infinite then there would be an ordinal in
  $\lim(M)\cap\lim(N)\setminus\lim(M\cap N)$ either as the limit of
  infinitely many cardinals $\lambda$ from clause~1, or else as the
  limit of infinitely many ordinals $\alpha$ from clause~1(iii).
  
  This completes the proof that if each member of $y$ is compatible
  with $\CC_M$ and $\lim(M)\cap\lim(N)=\lim(M\cap N)$ then there is an
  $M$-fence for $\CC_N$.
  \smallskip{}

  Now we verify the final paragraph of the lemma.
  Let $\II_{\alpha,\lambda}$ be
  any member of $y$, let $\II_{\alpha',\lambda'}$ be the minimal
  $M$-fence for $\II_{\alpha,\lambda}$, and suppose that  $x'$ is an  $M$-fence
  for $\CC_N$.   If  $\II_{\alpha',\lambda'}\ninx x'$ then by
  proposition~\ref{thm:minimal} there is a requirement
  $\OO_{\gamma,(\eta,\lambda]}\in M$ which is compatible with $x'$ but
  not with $\II_{\alpha,\lambda}$.  Now pick $\bar\lambda\in N\cap\lambda$ 
  such that $\bar\lambda>\eta$ and $\cof(\bar\lambda)=\omega$.    If
  $\alpha$ was given by clause~2(i) or~2(iii) then set
  $\bar\alpha=\alpha$; otherwise  pick $\bar\alpha\in
  \lim(C_{\alpha})$ such that $\gamma<\bar\lambda$.   Then
  $\II_{\bar\alpha,\bar\lambda}$ is a member of $N$ and is
  incompatible with $\OO_{\gamma,(\eta,\lambda']}$, contradicting the
  assumption that $x'$ is an $M$-fence for $\CC_N$. 
  This completes the proof that the fence $x$ is minimal among all
  $M$-fences for $\CC_N$.
  \smallskip{}

  The last paragraph shows something more:
  it did not assume that
  members of $y$ are compatible with $\CC_M$ or that
  $\lim(M)\cap\lim(N)=\lim(M\cap N)$, and    hence it implies that if there
  exists an $M$-fence $x'$ for $\CC_N$ then $x'$ must include a
  minimal $M$-fence for each member of $y$.  This implies that 
  each member of $y$ has an $M$-fence, and hence is compatible with
  $\CC_M$.   Also, since $x'$ is finite it follows that $y$ is finite, but
  it is easy to see that this implies that
  $\lim(M)\cap\lim(N)=\lim(M\cap N)$.
  This completes the proof of the right to
  left direction of the equivalence, and hence of
  lemma~\ref{thm:minMfenceN}. 
\end{proof}

%----------------------------------------------------------------
\subsection{Completeness}  
\label{sec:complete}
%----------------------------------------------------------------

 At the end of this subsection  we will give a complete
characterization, for any condition $p\in\Pkp$,  of the set of 
requirements
$\II_{\alpha,\lambda}$  such that
$p\forces\II_{\alpha,\lambda}\in\bigcup\dot G$.
For the proof of theorem~\ref{thm:mainthm}, however, we will not use
this characterization but rather two intermediate results.
The first of these will be needed in order to define the witness $p\mapsto
p\cut M$ for the strong genericity of a countable \simple{} \model{} $M$:

\begin{definition}\label{def:compl}
  If $p\in\Pkp$ and $\CC_{M,a}\in p$ then $\cp{M}p$ is the  set
  of all requirements
  $\II_{\alpha,\lambda}$ such that $\alpha=\min(M\setminus\alpha')$
  for some requirement $\II_{\alpha',\lambda}$ which is a minimal
  $M$-fence for some requirement in $p$. 
\end{definition}

Notice that every member of $\cp{M}p$ is a member of $M$, and is a
$M$-fence for the minimal $M$-fence from which it was defined and
hence for the requirement which demanded that minimal $M$-fence.  In
general $\cp{M}p$ need not include a complete set of $M$-fences for
members of $p$, since a minimal $M$-fence $\II_{\alpha',\lambda}$ may
have $\alpha'=\sup(M)$.   We will see later that if $M$ is simple then
this cannot happen.

\begin{lemma}\label{thm:compl}
  Suppose that $p\in\Pkp$ and $\CC_{M,a}\in p$. Then
  $p\cup\cp{M}p\in\Pkp$, 
  $p\cup\cp{M}p\eqx p$, and $\cp{M}{p\cup\cp{M}p}=\cp{M}{p}$.
\end{lemma}

The second asserts that the forcing $\Pkp$ does in fact add new
closed unbounded sets $D_{\alpha}$:
\begin{definition}\label{def:Da}
  If $G\subset\Pkp$ is generic and $\alpha<\kappa^+$ then we write
  $D_{\alpha}=\set{\lambda<\kappa:\II_{\alpha,\lambda}\in\bigcup G}$. 
\end{definition}
\begin{lemma}\label{thm:Da}
  The sequence $\vec D=\seq{D_{\alpha}:\alpha<\kappa^+}$ is a
  continuously diagonally decreasing sequence of closed unbounded
  subsets of $\kappa$.
\end{lemma}
The difficulty here is in showing that the sets $D_\alpha$ are closed;
the rest of lemma~\ref{thm:Da} can easily be proved with the machinery
already developed.

The proof of  lemma~\ref{thm:compl} will be given after the next two
lemmas, which contain the substance of the 
proof. 

\begin{lemma}\label{thm:compf}
  Suppose $p\in\Pkp$ and $\CC_{M,a}\in p$, and let $x$ be the set of
  minimal $M$-fences for requirements in  $p$.   Then $p\cup x\in
  \Pkp$,  $p\cup x\eqx p$, and $p\cup x$ includes an $M$-fence for
  every requirement in $p\cup x$.
\end{lemma}
\begin{proof}
  Let $\II_{\alpha,\lambda}$ be a 
  minimal 
  $M$-fence for one of the requirements  $\II_{\chi,\tau}$ or
  $\CC_{N,b}$ in $p$.   We will show that $\II_{\alpha,\lambda}$
  is compatible 
  with all requirements $\OOO$ and $\CC_{M',a'}$ in $p$.
  Since $\II_{\alpha,\lambda}$ is also a fence for any  $q\le p$ in $\Pkp$, 
  this will imply that $q\cup\sing{\II_{\alpha,\lambda}}\in P^{*}$.
  It 
  follows that  $p\forces\II_{\alpha,\lambda}\in\bigcup\dot G$, that is,
  that $p\eqx p\cup\sing{\II_{\alpha,\lambda}}$.    
  This will be sufficient to prove the lemma, since it follows by 
  an easy induction  that $p\cup x\eqx p$.

  First we show that $\II_{\alpha,\lambda}$ is compatible with every
  requirement $\OOO\in p$.
  In the case that $\OOO\notin M$ the compatibility of $\OOO$ and
  $\CC_{M,a}$ 
  implies that $\OOO$ is compatible with every requirement
  $\II_{\alpha',\lambda'}\in M[a]$, and since $\lambda\in M$ and
  $\alpha\in M\cup\lim(M)$ 
  this implies  that $\OOO$ is compatible with
  $\II_{\alpha,\lambda}$.
  Thus we can assume that $\OOO\in M$.    If $\II_{\alpha,\lambda}$ is
  is the minimal fence for $\II_{\chi,\tau}\in p$ then $\OOO$ is
  compatible with $\II_{\chi,\tau}$, since both are in $p$, and by the
  minimality of $\II_{\alpha,\lambda}$ it follows that $\OOO$ is
  compatible with $\II_{\alpha,\lambda}$.
  On the other hand, if
  $\II_{\alpha,\lambda}$ is a minimal fence for $\CC_{N,b}\in p$ then it
  follows from $\eta\ge\lambda\geq\sup(M\cap N\cap\kappa)$ that
  $\OOO\notin N$, and hence the compatibility of $\OOO$ with
  $\CC_{N,b}$ implies that $\OOO$ is compatible with every requirement
  $\II_{\alpha',\lambda'}\in N[b]$.   Again, the minimality of
  $\II_{\alpha,\lambda}$ then implies that $\OOO$ is compatible with
  $\II_{\alpha,\lambda}$. 

  \smallskip{}

  Now we show that $\II_{\alpha,\lambda}$ is compatible with any
  requirement $\CC_{M',a'}\in p$.
  The proof proceeds by verifying the final statement of the lemma,
  by showing if $\II_{\alpha,\lambda}$ is not its own $M'$-fence then
  the minimal $M'$-fence for $\II_{\alpha,\lambda}$ is the same as
  the minimal
  $M'$-fence for some requirement in $p$.

  If $\lambda\ge\sup(M\cap M'\cap\kappa)$ then any $M'$-fence for
  $\CC_{M,a}$   includes a 
  $M'$-fence for $\II_{\alpha,\lambda}$, so we can assume that
  $\lambda<\sup(M\cap M'\cap\kappa)$.
  If $\lambda\in M'$ then $\II_{\alpha,\lambda}$ is its own
  $M'$-fence, so we can assume that $\lambda\notin M'$.   Hence $M\cap
  M'\cap\kappa\in M$.
  We will show that any $M'$-fence for $\II_{\chi,\tau}$ or for
  $\CC_{N,b}$ is or includes a $M'$-fence for $\II_{\alpha,\lambda}$.

  To this end, suppose $\OOO$ is some requirement in $M'$ which is
  incompatible 
  with $\II_{\alpha,\lambda}$.    Since $\lambda<\sup(M\cap M'\cap\kappa)$ we
  can assume that $\eta<\sup(M\cap M'\cap\kappa)$, so that $M\cap M'\in M$
  implies that $\eta',\eta\in M$.
  Since $\OOO$ is incompatible with $\II_{\alpha,\lambda}$, we have
  $\gamma\in A_{\alpha,\lambda}\cap M'$ and it follows by
  lemma~\ref{thm:ASminBig} that $\gamma\in M$.
  Hence $\OOO$ is in $M$.

  If $\II_{\alpha,\lambda}$ is the minimal $M$-fence for
  $\II_{\chi,\tau}$ then the minimality of 
  $\II_{\alpha,\lambda}$ implies that $\OOO$ is  incompatible with
  $\II_{\chi,\tau}$, and hence is incompatible with the
  $M'$-fence for $\II_{\chi,\tau}$.
  This shows that any $M'$-fence for $\II_{\chi,\tau}$ is a $M'$-fence
  for $\II_{\alpha,\lambda}$, and thus implies that
  $\II_{\alpha,\lambda}$ is compatible with $\CC_{M',a'}$.

  If $\II_{\alpha,\lambda}$ is a minimal $M$-fence for $\CC_{N,b}$
  then it follows
  similarly that $\OOO$ is incompatible with some
  $\II_{\alpha',\lambda'}\in N[b]$.
  We will show that $\eta\ge\sup(M'\cap N\cap\kappa)$.   It then follows
  that $\OOO$ is 
  incompatible with some member of any $M'$-fence for $\CC_{N,b}$, and
  this implies that any $M'$-fence for $\CC_{N,b}$ includes a $M'$-fence
  for $\II_{\alpha,\lambda}$ and 
  hence  completes the proof of lemma~\ref{thm:compl}.

  Suppose to the contrary that $\eta<\sup(M'\cap N\cap\kappa)$.
  If $M'\cap N\cap\kappa$ is an initial segment of $M'$,
  that is, $M'\cap\sup(M'\cap N\cap\kappa)\subseteq N$,
  then $\eta\in M\cap M'\cap\sup(M'\cap N\cap\kappa)\subset M\cap N$; 
  however this is impossible since
  the fact that $\II_{\alpha,\lambda}$ is the minimal $M$-fence for
  $\CC_{N,a}$ 
  implies that  $\sup(M\cap N\cap\kappa)\le\lambda<\eta$.
  Hence we must have $M'\cap N\in M'$, and it follows that
  $\sup(N\cap\lambda)\in M'\cap\lambda\subset M$.
  However by lemma~\ref{thm:minMfenceN}, the fact that 
  $\II_{\alpha,\lambda}$ is in the minimal 
  $M$-fence for $\CC_{N,b}$ implies that
  $\sup(N\cap\lambda)<\lambda=\min(M\setminus\sup(N\cap\lambda))$, so
  that $\sup(N\cap\lambda)\notin M$.
\end{proof}

\begin{lemma}\label{thm:nDg}
  Suppose that $p\in\Pkp$, $\alpha<\kappa^+$ is a nonzero limit
  ordinal and $p\forces\lambda\notin\dot D_{\alpha}$.  Then there is a
  successor ordinal $\gamma\in A_{\alpha,\lambda}$ such that
  $p\forces\lambda\notin \dot D_{\gamma}$.

  Furthermore if $\CC_{M,a}\in p$, $\lambda\in M\cup\lim(M)$, and
  $\alpha\in M\cup a(\lambda)$ 
  then the least such ordinal $\gamma$ is a member of $M$.
\end{lemma}

\begin{proof}
  The hypothesis that $p\forces\lambda\notin\dot D_{\alpha}$ could
  hold either because
  $\lambda\notin B^*_{\alpha}$, so that
  $\II_{\alpha,\lambda}$ is not a requirement, or because
  $\II_{\alpha,\lambda}$ is incompatible with some requirement in $p$.
  
  If $\lambda\notin B^*_{\alpha}$ then there is a successor
  $\gamma\in A_{\alpha,\lambda}$ such that $\lambda\notin B^*_{\gamma}$,
  and hence 
  $\nothing\forces\lambda\notin\dot D_{\gamma}$.
  Furthermore, 
  if $\CC_{M,\alpha}$ is as in the second paragraph then
  $\lambda\in M$, since $\lambda\notin B^{*}_{\alpha}$ implies that
  $\cof(\lambda)>\omega$, and hence $\lambda\notin\lim(M)$.
  Then it is easy to
  see that
  there is some such $\gamma$ in $M$, using elementarity if
  $\alpha\in M$ and definition~\ref{def:requireM}(\ref{item:x4}) if
  $\alpha\in a(\lambda)$.
  
  Thus we can assume that $\II_{\alpha,\lambda}$ is a requirement
  and that
  $\II_{\alpha,\lambda}$ is incompatible with some requirement in
  $p$.   Now if a requirement $\OOO\in p$ is incompatible with
  $\II_{\alpha,\lambda}$ 
  then $p\le\sing{\OOO}\forces\lambda\notin\dot D_\gamma$.
  If the hypothesis of the second paragraph holds then 
  the compatibility of $\OOO$ with $\CC_{M,a}$ implies that $\OOO\in
  M$ and thus $\gamma\in M$.

  The last possibility is that $\II_{\alpha,\lambda}$ is incompatible
  with 
  $\CC_{N,b}\in p$.
  Then $\lambda'\deq\min(N\setminus\lambda)\notin
  B^*_{\beta}$ where $\beta=\sup\set{\xi+1:\xi+1\in
    A_{\alpha+1,\lambda}\cap N}$.
  It follows that there is a successor ordinal $\gamma\in N\cap
  A_{\alpha,\lambda}$ such that $\lambda'\notin B^*_{\gamma}$.
  Thus $p\le\sing{\CC_{N,b}}\forces\lambda\notin \dot D_{\gamma}$. 

  Now suppose that the hypothesis of the  second paragraph holds.
  If $\lambda\ge\sup(M\cap   N\cap\kappa)$
  then any $N$-fence for $\CC_{M,a}$ includes a $N$-fence for
  $\II_{\alpha,\lambda}$, contradicting the assumption that
  $\II_{\alpha,\lambda}$ is incompatible with $\CC_{N,b}$.  Thus we must
  have $\lambda<\sup(M\cap 
  N\cap\kappa)$.    Also $\lambda\notin N$, or else $\II_{\alpha,\lambda}$
  would be its own $N$-fence, and   hence we must have $N\cap M\in M$.
  Then lemma~\ref{thm:ASminBig} implies that  
  $\gamma\in A_{\alpha,\lambda}\cap N\subset M$.   
  \smallskip{}

  This completes the proof of lemma~\ref{thm:nDg}, except that under
  the hypothesis of the second paragraph we
  have only shown that there exists a successor ordinal $\gamma\in
  A_{\alpha,\lambda}\cap M$ such that $p\forces\lambda\notin\dot
  D_{\gamma}$, not that the least such $\gamma$ is a member of $M$.
  Now let $\gamma'$ be the least ordinal in
  $A_{\alpha,\lambda}\cap M$ such that $p\forces\lambda\notin\dot
  D_{\gamma'}$.   Then $\gamma'$ cannot be a limit ordinal, since in
  that case we could apply the lemma with $\gamma'$ in place of
  $\alpha$.  Thus $\gamma'$ must be a successor ordinal, say
  $\gamma'=\gamma''+1$, but then $\gamma''\in M$ and hence
  $p\not\forces\lambda\notin\dot D_{\gamma''}$.   Thus $\gamma'$ is the
  least member $\gamma$ of $A_{\alpha,\lambda}$ such that
  $p\forces\lambda\notin\dot D_{\gamma}$.
\end{proof}

\begin{proof}[Proof of lemma~\ref{thm:compl}]
  We need to show that every requirement $\II_{\alpha,\lambda}\in
  \cp{M}p$ is compatible with $p$.    If $\II_{\alpha,\lambda}$ is not
  compatible with $p$ then $p\forces\lambda\notin \dot D_{\alpha}$,
  and by lemma~\ref{thm:nDg} it follows that $p\forces\lambda\notin
  \dot D_{\gamma}$ for some $\gamma\in M\cap A_{\alpha,\lambda}$.
  However by the definition of $\cp{M}p$ we have
  $\alpha=\min(M\setminus \alpha')$ where $\II_{\alpha',\lambda}$ is a
  minimal $M$-fence for some requirement in $p$.   Now
  $A_{\alpha,\lambda}\cap M=A_{\alpha',\lambda}\cap M$, so it follows
  that $\gamma\in A_{\alpha',\lambda}$ and hence
  $\II_{\alpha',\lambda}$ is incompatible with $p$; however this
  contradicts lemma~\ref{thm:compf}. 
\end{proof}

One more lemma is needed for the proof of  lemma~\ref{thm:Da}.

\begin{lemma}
  \label{thm:makecompP}
  Suppose that $\CC_{M,a}\in p$, and either
  \begin{inparaenumi}
  \item     $\II_{\alpha,\lambda}=\II_{\sup(M),\sup(M\cap\kappa)}$,
    or else
  \item
      $\lambda=\sup(M\cap\lambda')<\lambda'$ for some
      $\II_{\alpha,\lambda'}\inx p$ 
      with $\lambda'\in M$ and $\alpha\in M\cup\lim(M)\cup a(\lambda)$.
  \end{inparaenumi}
  Then $\II_{\alpha,\lambda}$ is compatible with $p$, and indeed
  $p\cup\sing{\II_{\alpha,\lambda}}\eqx p$.
\end{lemma}
\begin{proof}
  As in previous lemmas, it will be sufficient to show that 
  $\II_{\alpha,\lambda}$ is compatible with $p$, since this implies
  that $\II_{\alpha,\lambda}$ is compatible with any $p'\leq p$ and
  hence $p\lex p\cup\sing{\II_{\alpha,\lambda}}$. 
  \smallskip{}

  Since $\cof(\lambda)=\omega$, $\lambda\in B^*_{\alpha}$ and hence
  $\II_{\alpha,\lambda}$ is a requirement.    We first show that
  $\II_{\alpha,\lambda}$ is compatible with any requirement $\OOO\in
  p$.
  In case~(i), where $\lambda=\sup(M)$,  any requirement
  $\OOO$ which is incompatible with $\II_{\alpha,\lambda}$ is
  incompatible with $\CC_{M,a}$, and hence is not in $p$.

  In case~(ii), with $\lambda=\sup(M\cap\lambda')$ where
  $\lambda'\in M$ and 
  $\alpha\in M\cup\lim(M)\cup a(\lambda)$,  any requirement
  $\OOO$  incompatible 
  with $\II_{\alpha,\lambda}$ is incompatible with requirements
  $\II_{\alpha'',\lambda''}\in M[a]$ and hence must be a member of $M$,
  but this implies that $\eta\ge\lambda'$ so that $\OOO$ is
  incompatible with $\II_{\alpha,\lambda'}\inx p$.   

 Thus
  $\II_{\alpha,\lambda}$ is compatible with every requirement $\OOO\in
  p$.   Now we show that $\II_{\alpha,\lambda}$ is compatible with any
  requirement   $\CC_{N,b}$ in $p$.   If  
  $\lambda\ge\sup(M\cap N\cap\kappa)$ then any requirement $\OOO\in N$
  which is incompatible with $\II_{\alpha,\lambda}$ is incompatible
  with requirements $\II_{\alpha'',\lambda''}\in M[a]$ and hence must be
  incompatible with some member of any $N$-fence for $\CC_{M,a}$.
  Hence the $N$-fence for $\CC_{M,a}$ includes a $N$-fence for
  $\II_{\alpha,\lambda}$. 

  Thus we can assume that $\lambda<\sup(M\cap N\cap\kappa)$.   In
  particular, 
  $\lambda\not=\sup(M\cap\kappa)$, so $\II_{\alpha,\lambda}$ comes
  from clause~(ii) for   some  requirement $\II_{\alpha,\lambda'}\inx p$.
  If
  $N\cap\sup(M\cap N\cap\kappa)\subset M$ then
  $\min(N\setminus\lambda)=\min(N\setminus\lambda')$, so any $N$-fence for
  $\II_{\alpha,\lambda'}$ is an $N$-fence for
  $\II_{\alpha,\lambda}$.   Otherwise $M\cap N\in N$, so 
  $\lambda\in N$ and hence $\II_{\alpha,\lambda}\in p'$ is its own
  $N$-fence.  
\end{proof}

\begin{proof}[Proof of lemma~\ref{thm:Da}]
\newcommand{\OOah}{\OO_{\alpha,(\eta',\eta]}}
  To see that $D_\alpha$ is unbounded in $\kappa$, let $p$ be any
  condition in $\Pkp$ and suppose $\zeta<\kappa$.
  Pick $\tau>\zeta$ of cofinality $\omega$ so that $\tau>\eta$ for
  all requirements $\OO_{\alpha,(\eta',\eta]}\in p$ and
  $\tau>\sup(M\cap\kappa)$ for all 
  $\CC_{M,a}\in p$.   Then $\II_{\alpha,\tau}$ is a requirement,
  $p'\deq p\cup\sing{\II_{\alpha,\tau}}\le p$, and
  $p'\forces\tau\in\dot 
  D_{\alpha}\setminus\eta$.

  Proposition~\ref{thm:inx} implies that if $\lambda\in D_{\alpha}$
  and $\alpha'\in A_{\alpha,\lambda}$ then $\lambda\in D_{\alpha'}$,
  so  $\vec D$ is diagonally decreasing.
  Lemma~\ref{thm:nDg} implies that if $\alpha$ is a limit ordinal and 
  $p\forces\lambda\in\dot 
  D_{\alpha'}$ for all $\alpha'\in 
  A_{\alpha,\lambda}$ then $p\forces \lambda\in\dot D_{\alpha}$, so
  $\vec D$ is continuously 
  diagonally decreasing. 
  \smallskip{}

  Thus it only remains to show that  $D_\alpha$ is
  closed for each $\alpha<\kappa^{+}$.
  We will show that for any condition $p$ and ordinals $\alpha$ and
  $\lambda$ such that 
  $p\nforces\lambda\in \dot D_\alpha$, there is a
  requirement $\OOah$, compatible with $p$, such
  that $\eta'<\lambda\leq\eta$.
  Then 
  $p\cup\sing{\OOah}$ is a condition extending $p$ which forces that
  $\dot D_{\alpha}\cap\lambda\subseteq\eta'$, so that 
  $\lambda$ is not a limit point of $D_{\alpha}$.

  By extending $p$ if necessary, and taking $\alpha$ to be minimal, we
  may assume that $p$ forces that
  $p\forces\lambda\in\left(\bigcap_{\alpha'\in A_{\alpha,\lambda}}\dot
  D_{\alpha'}\right)\setminus  \dot  D_{\alpha}$.
  It follows by lemma~\ref{thm:nDg} that $\alpha$ is either $0$ or 
  a successor ordinal, say $\alpha=\alpha_0+1$. 
  The case $\alpha=0$ is identical to lemma~\ref{thm:OOout}, so we
  will assume $\alpha>0$.  By further extending $p$ if necessary, we may
  assume that there is an ordinal
  $\tau>\lambda$ such that $p\forces \tau=\min(\dot
  D_{\alpha}\setminus\lambda)$.  

  Let $Y=\set{\CC_{M,a}\in p:\alpha\in M\And\tau\in\lim(M)}$.
  If $Y=\emptyset$ then set $\eta=\lambda$.   Otherwise note that for
  any two members $\CC_{M,a}$ and $\CC_{M',a'}$ of $Y$, the fact that 
  $\tau\in\lim(M)\cap\lim(M')=\lim(M\cap M')$ implies that $\sup(M\cap
  M'\cap\kappa)\geq\tau$ and hence one of $M\cap\tau$ and $M'\cap\tau$
  is contained in the other.  Thus 
  $\set{M\cap \tau:\CC_{M,a}\in Y}$ is linearly ordered by
  $\subseteq$.  Pick $\CC_{M,a}\in Y$ with
  $M\cap\tau$ minimal, and 
  let $\eta=\min(M\setminus\lambda)$.    The desired requirement will be
  $\OOah$ for some suitably choosen
  $\eta'<\lambda$.
  
  If $Y\not=\emptyset$ then 
  the choice of $\eta$ ensures that $\OOah$ is compatible with any
  requirement $\CC_{M',a'}\in Y$ so long as $\eta'\in M$.
  It remains to show that $\eta'$ can be choosen so that $\OOah$ is
  also 
  compatible with the requirements $\II_{\gamma,\xi}\in p$ and
  $\CC_{M',a'}\in p\setminus Y$.
  \smallskip{}

  Let $\II_{\gamma,\xi}$ be a requirement in $p$.   If
  $\alpha\notin A_{\gamma+1,\xi}$ then
  $\II_{\gamma,\xi}$ is compatible with $\OOah$ for any
  $\eta'<\lambda$, so we can assume that $\alpha\in A_{\gamma+1,\xi}$.
  It follows that $p\forces\xi\in\dot D_{\alpha}$, so the choice of
  $\tau$ ensures that $\xi<\lambda$ or $\xi\geq\tau$.   If
  $\xi\geq\tau$ then $\OOah$ is compatible with $\II_{\alpha,\lambda}$
  for any $\eta'<\lambda$, so we can assume that $\xi<\lambda$.
  Then $\OOah$ is compatible with $\II_{\gamma,\xi}$ for any
  $\eta'\in\lambda\setminus\xi$.   If $Y=\emptyset$ then we are done;
  otherwise we need to show that $\sup(M\cap\lambda)>\xi$.   To see
  this, note that if $\xi'=\min(M\setminus\xi)$ then
  $p\forces\xi'\in\dot D_{\alpha}$ because of the $M$-fence for
  $\II_{\gamma,\xi}$.
  Since $\xi'\leq\eta<\tau$ it follows that $\xi'<\lambda$. 

  It remains to consider requirements $\CC_{M',a'}\in p\setminus Y$.
  We first show that if $\CC_{M',a'}\in p$ and $\alpha\notin M'$ then
  $\OOah$ is compatible with $\CC_{M',a'}$ for any $\eta'<\lambda$. 
  Suppose to the contrary that
  $\OOah$ is incompatible with $\CC_{M',a'}$.
  Then there is some ordinal $\alpha'\in M'\cup a'(\eta)$
  such that $\alpha\in
  A_{\alpha',\lambda}$.   The least such ordinal $\alpha'$ is a
  limit ordinal, and  $p\forces\eta\notin\dot D_{\alpha'}$, so
  lemma~\ref{thm:nDg} 
  implies that 
  there is a successor ordinal  $\gamma\in 
  M'\cap A_{\alpha',\lambda}$ such that $p\forces\eta\notin\dot
  D_{\gamma}$.   By the choice of $\alpha'$, we must have
  $\gamma<\alpha$ and hence $\gamma\in
  A_{\alpha,\lambda}=A_{\alpha_0+1,\lambda}$, and it follows that
  $p\forces\eta\notin\dot D_{\alpha_0}$.   If $\eta=\lambda$ this
  contradicts the choice of $\alpha$.   If $\eta>\lambda$ then
  $\II_{\alpha_0,\eta} \inx p$ as the $M$-fence for
  $\II_{\alpha_0,\lambda}\inx p$, and so again $p\forces\eta\in\dot
  D_{\alpha_0}$.   This contradiction completes the proof that
  $\OOah$ is compatible with $\CC_{M',\alpha'}$.

  The only remaining requirements to consider are 
  $\CC_{M',a'}\in p\setminus Y$ with $\alpha\in M'$.   Now note that
  lemma~\ref{thm:compl}  implies that $p\forces\sup(M'\cap\tau)\in\dot
  D_{\alpha}$,   so 
  $\sup(M'\cap\tau)<\lambda$.   
  If $Y=\emptyset$ then it follows that $\OOah$ is compatible with
  $\CC_{M',a'}$ so long as  $\eta'>\sup(M'\cap\lambda)$. 
  If $Y\not=\emptyset$ then  we must show that
  $\sup(M\cap\lambda)>\sup(M'\cap\lambda)$, so that $\eta'$ can be
  choosen to be a member of $M$.
  Suppose first that 
  $\delta\deq\sup(M\cap   M'\cap\kappa)\geq\tau$.   Since $\eta\in
  M\setminus M'$ this implies that $M\cap M'\in M$, and in particular
  $\sup(M\cap M'\cap\lambda)\in M$, so $\sup(M\cap
  \lambda)>\sup(M'\cap\lambda)$.   Now suppose that $\delta<\tau$, and
  hence $\delta<\lambda$.   If $M\cap M'\in M$ and
  $\delta=\sup(M'\cap\lambda)$ then again $M'\cap\lambda\subset M$
  and hence $\sup(M\cap\lambda)>\sup(M'\cap\lambda)$.
  Otherwise, $p\forces \min(M\setminus\sup(M'\cap\lambda))\in\dot
  D_{\alpha}$ because of the $M$-fence for $\CC_{M',\alpha}$, so as in
  the case of $\II_{\gamma,\xi}$ it follows that
  $\min(M\setminus\sup(M'\cap\lambda))<\lambda$ and hence
  $\sup(M\cap\lambda)>\sup(M'\cap\lambda)$. 
  \smallskip{}
  
  It follows that if $\eta'<\lambda$ is choosen so that $\eta'>\xi$
  for all $\II_{\gamma,\xi}\in p$ with $\xi<\lambda$ and
  $\II_{\alpha,\xi}\inx\sing{\II_{\gamma,\xi}}$, and 
  $\eta'>\sup(M'\cap\lambda)$ for all $\CC_{M',a'}\in p\setminus Y$
  with $\alpha\in M'$, then $\OOah$ is
  compatible with all requirements in $p\setminus Y$;
  furthermore, we have shown that if $Y\not=\emptyset$ then such 
  ordinals $\eta'$ can be found in $M\cap\lambda$.   Such a choice of
  $\eta'$ gives a
  condition $\OOah$ compatible with $p$ which forces that $\dot
  D_{\alpha}$ is bounded in $\lambda$, and    it follows that 
  $D_\alpha$ is closed. 
\end{proof}

This completes the proof of lemmas~\ref{thm:compl} and~\ref{thm:Da}.
In the remainder of this subsection, which is not needed for the proof
of theorem~\ref{thm:mainthm}, we briefly explain how the proofs of
these  lemmas can be used to give a characterization, for an
arbitrary condition $p$, 
of the pairs $(\alpha,\lambda)$ such that 
$p\forces\lambda\in\dot 
D_{\alpha}$. 

This characterization generates the set of such pairs through four
steps. 
We write
$A^{*}_{\alpha,\lambda}$ for 
the intersection of $A_{\alpha,\lambda}$ with
$\kappa^{+}\setminus\lim(\kappa^+)$, that is, $A^*_{\alpha,\lambda}$
contains only $0$ and the successor ordinals from $A_{\alpha,\lambda}$.

\paragraph{Step 1.}
 By lemma~\ref{thm:compl}, we can assume without loss of
  generality that $p$ includes, for each requirement $\CC_{M,a}\in p$,
  the minimal $M$-fence for each requirement in $p$.

\paragraph{Step 2.}
Suppose $\CC_{M,a}\in p$, $\lambda\in M\cup\lim(M)$, and 
$\alpha\in M\cup a(\lambda)$.   
If $\II_{\alpha',\lambda}\inx p$, where $\alpha'=\sup(M\cap
A^*_{\alpha,\lambda})$, then $p\forces\lambda\in\dot D_{\alpha}$. 

This follows from lemma~\ref{thm:nDg}.
With some care it can be shown that any condition $p$ can be extended
to a condition $p'\le p$ with $p'=^*p$ so that $I_{\alpha,\lambda}\inx
p$ 
for each pair $(\alpha,\lambda)$ as in this step.    A key point in
the argument is that the requirements $\CC_{M,a}\in p$ should be considered
in the order of their size: if $\lambda\in M\cap M'\in M'$ then the pairs
$(\alpha,\lambda)$ from $\CC_{M,a}$ should be dealt with before those
from $\CC_{M',a'}$.

\paragraph{Step 3.}
If $\CC_{M,a}\in p$, $\lambda\in M$,
and $\lambda'=\sup(M\cap\lambda)<\lambda$ then
$p\forces\lambda'\in\dot D_{\alpha}$ whenever $\alpha\in M\cup
a(\lambda)$ and $p\forces\lambda\in\dot D_{\alpha}$. 

This follows from lemma~\ref{thm:makecompP}, and it is straightforward
to verify that any condition $p$ can be extended to a condition $p'\le
p$ with $p'=^*p$ such that $\II_{\alpha,\lambda'}\inx p'$ for all
$(\alpha,\lambda')$ as in this step.

\paragraph{Step 4.}
If $\alpha'=\sup(A^{*}_{\alpha,\lambda})$ and
$p\forces\lambda\in\dot D_{\alpha'}$ then $p\forces\lambda\in\dot
D_{\alpha}$. 

This follows immediately from the fact that the sequence $\vec D$ is
continuously decreasing.   This situation is actually an artifact of
our definition of the 
sets $A_{\alpha,\lambda}$: 
if these sets had been defined to be closed under successor then it
would never happen that 
$\sup(A_{\alpha,\lambda})>\sup(A^*_{\alpha,\lambda})$.  
Such a change would make this characterization more natural.
\smallskip{}

The proof of lemma~\ref{thm:Da} shows that for any condition $p$
which has been extended as described in steps~1-3, the only pairs
$(\alpha,\lambda)$ for which 
$p\forces\lambda\in\dot D_{\alpha}$ 
are those such that $\II_{\alpha,\lambda}\inx p$
and those coming from step~4 in which $\II_{\alpha',\lambda}\inx p$.

\subsection{Strongly generic conditions}
\label{sec:M'a'}

Earlier we described three types of \simple{} \model{}s:  in addition
to the countable models we have uncountable models $X\prec
H_{\kappa^+}$ of size less less than $\kappa$ with $X\cap H_\kappa$
transitive, and transitive models $X\prec H_{\kappa^+}$ with
$\kappa\subset X$.
The main result of this section asserts that each of these has a
strongly generic condition: 

\begin{lemma}
  \label{thm:smc}
  Suppose that $X$ is a \simple{} \model{}, and that
  $\sup(X\cap\kappa)\in B_{\sup(X)}$ if $\omega<\card{X}<\kappa$.
  Set  
  \begin{equation*}
    p^X=
    \begin{cases}
      \sing{\CC_X}&\text{if $X$ is countable,}\\
      \sing{\II_{\sup(X),\sup(X\cap\kappa)}}&\text{if $X\cap\kappa\in\kappa$,}\\
      \nothing&\text{if $X$ is transitive.}
    \end{cases}
  \end{equation*}
  Then $p^{X}$ is a tidy strongly $X$-generic condition.
\end{lemma}

The reason for requiring that $X$ be a  \simple{} \model{} is given by the
following observation, which will be used to define the function
$p\mapsto p\cut X$ witnessing
strong genericity.
\begin{proposition}\label{thm:AbndM}
  Suppose that $X$ is a \simple{} \model{}, $\alpha<\kappa^+$, and
  $\lambda<\sup(X\cap\kappa)$.   Then $A_{\alpha,\lambda}\cap X$ is
  bounded in $X$.
\end{proposition}
\begin{proof}
  Suppose to the contrary that
  $A_{\alpha,\lambda}$ is unbounded in $X$.   Set
  $\bar\alpha=\sup(X)\leq\sup(A_{\alpha,\lambda})\leq\alpha$ 
  and  let $\alpha'$ be the 
  least member of 
  $A_{\alpha+1,\lambda}\setminus \bar\alpha$.
  Then  $A_{\alpha,\lambda}\cap
  \alpha'=A_{\alpha',\lambda}\subset c_{\alpha',\lambda}$.   Now
  $\bar\alpha\in\lim(C_{\alpha'})$ since $A_{\alpha,\lambda}$ is
  unbounded in $X$, so
  $C_{\bar\alpha}=C_{\alpha'}\cap\bar\alpha$.   However,  since $X$ is
  \simple{} we have 
  $\otp(C_{\bar\alpha})=\sup(X\cap\kappa)>\lambda$ and hence
  $\sup(X\cap
  A_{\alpha',\lambda})\le c_{\alpha',\lambda}=c_{\bar\alpha,\lambda}<\bar\alpha$, 
  contrary to assumption. 
\end{proof}

In order to make use of this fact we extend  to arbitrary models some
of the notation previously 
associated to countable models $M$.  Recall that 
$A^{*}_{\alpha,\lambda}=
A_{\alpha,\lambda}\cap(\kappa^+\setminus\lim(\kappa^+))$.   

\begin{definition}\label{def:Xfence}
  If $X$ is an uncountable model then $\II_{\alpha',\lambda}$ is an
  $X$-fence for $\II_{\alpha,\lambda}$ if $\lambda\in X$ and every 
  requirement $\OOO\in X$ which is incompatible with
  $\II_{\alpha,\lambda}$ is also incompatible with $\II_{\alpha',\lambda}$.

  We say that $\II_{\alpha',\lambda}$ is the \emph{minimal $X$-fence}
  for $\II_{\alpha,\lambda}$ if $\alpha'=\sup(A^*_{\alpha,\lambda}\cap
  X)$.

  We write $\cp{X}p$ for the set of requirements
  $\II_{\alpha'',\lambda}$ such that
  $\alpha''=\min(X\setminus\alpha')$ where $\II_{\alpha',\lambda}$ is
  the minimal $X$-fence for some requirement $\II_{\alpha,\lambda}\in
  p$. 
\end{definition}
Note that the definitions are identical to those given previously for
countable models (except that if $M$ is countable then $\cp{M}p$ also
includes $M$-fences for 
requirements $\CC_{N,b}\in p$).  

The $X$-fences from definition~\ref{def:Xfence} have the same
properties as $M$-fences.   

\begin{proposition}\label{thm:cpXp}
  If  $X$ is  a \model{} of any type, and $p\le p^X$ if $X$ is
  countable, then  
  $p\cup\cp{X}p\in\Pkp$, and $p\cup\cp{X}p\eqx p$.  
  Furthermore if $X$ is simple then  $\cp{X}p$
  includes an
  $X$-fence for every requirement 
  $\II_{\alpha,\lambda}\in p\cup\cp{X}P$ with
  $\lambda<\sup(X\cap\lambda)$, and if $X$ is countable as well as
  simple then $\cp{X}p$   also 
  includes an $X$-fence for every requirement $\CC_{N,b}\in p$. 
\end{proposition}
\begin{proof}
  The first statement was proved   for countable models as
  lemma~\ref{thm:compl}, so we can assume   
  that $X$ is uncountable and hence $X\cap\kappa$ is transitive, 
  Suppose that $\II_{\alpha,\lambda}\in \cp{X}p$, say that
  $\alpha=\min(X\setminus\alpha'')$ where
  $\alpha''=\sup(A^*_{\alpha',\lambda}\cap X)$ for some
  $\II_{\alpha',\lambda}\in p$.
  Then $\II_{\alpha'',\lambda}\inx p$, so $\II_{\alpha'',\lambda}$ is
  compatible with $p$.   But since $X\cap\lambda$ is transitive and
  $\lambda$ and $\alpha$ are in $X$, the set $A_{\alpha,\lambda}$ is a
  subset of $X$ and hence is contained in $A_{\alpha'',\lambda}$.
  Thus any fence for $\II_{\alpha'',\lambda}$ is also a fence for
  $\II_{\alpha,\lambda}$. 

  For the second statement, if $X$ is
  uncountable then $\lambda<\sup(X\cap\lambda)$ implies that
  $\lambda\in X$, and if $X$ is countable then the compatibility of
  $p$, together with the assumption that $p\leq p^X$, imply that every
  stated requirement has a fence which is a requirement
  $\II_{\alpha,\lambda}$ (or a finite set of such requirements) with
  $\lambda\in X$.   Thus it is enough to 
  show that if $\II_{\alpha',\lambda}$ is the minimal $X$-fence for
  any requirement $\II_{\alpha,\lambda}$ with $\lambda\in X$ then
  $\alpha'<\sup(X)$.    Since $\alpha'=\sup(A^{*}_{\alpha,\lambda}\cap
  X)$, this follows from proposition~\ref{thm:AbndM}. 
\end{proof}

We are now ready to start the proof of lemma~\ref{thm:smc}.   The
function $p\mapsto p\cut X$ witnessing the strong $X$-genericity of
$p^X$ is defined by the equation
\begin{equation*}\label{eq:witness}
  p\cut X= (p\cap X)\cup\cp{X}p\cup\set{\icut{\CC_{M,a}} X:\CC_{M,a}\in
  p\And M\cap X\in X} 
\end{equation*}
where $\icut{\CC_{M,a}} X$ is given by the following lemma:
\begin{lemma}
  \label{thm:CinM}
  Suppose that $X$ is a \simple{} \model{}
  and $\CC_{M,a}$ is a 
  requirement compatible with $p^X$ such that $M\cap X\in X$.
  Then there is a requirement $\icut{\CC_{M,a}}X$ in $X$ 
  such that
  \begin{inparaenumi}
  \item every requirement $\RR$ in $X$ which is compatible with
    $\sing{\icut{\CC_{M,a}}X}\cup\cp{X}{p}$ 
    is also compatible with 
    $\CC_{M,a}$, and 
  \item every requirement $\RR$ which is
    compatible with 
    $\sing{\CC_{M,a}}\cup p^X$ is also compatible with $\icut{\CC_{M,a}}X$.
  \end{inparaenumi}
 \end{lemma}

The proof of lemma~\ref{thm:CinM} will take up most of this
subsection.   We first show that lemma~\ref{thm:smc} follows from
lemma~\ref{thm:CinM}.

\begin{proof}[Proof of lemma~\ref{thm:smc} from lemma~\ref{thm:CinM}]
First we verify that  $p\cut X\in\Pkp$, that is, that any two
requirements in $p\cut X$ are compatible.
Any two requirements in $p\cup\cp{X}p$ are
compatible by proposition~\ref{thm:cpXp}, and if $\CC_{M,a}\in p$ and
$M\cap X\in X$ then $\icut{\CC_{M,a}}X$ is compatible with every requirement
in $p\cup\cp{X}p$ by clause~\ref{thm:CinM}(ii).   Finally, if
$\CC_{N,b}$ is another member of $p$ such that $X\cap N\in X$ then
the compatibility of $\icut{\CC_{M,a}}X$ with $\icut{\CC_{N,b}}{X}$
follows from clause~\ref{thm:CinM}(ii) together with the fact that 
$\icut{\CC_{M,a}}X$ is compatible with $p\le\sing{\CC_{N,b}}\cup
p^X$.

Next we verify that 
the function $p\mapsto p\cut X$  is tidy.  Suppose that $p,p'\le
p^{X}$ are compatible conditions.    Then $p\land p'=p\cup p'$, and
$(p\land p')\cut X=p\cut X\cup p'\cut X=p\cut X\land p'\cut X$ since 
each member of $(p\land p')\cut X$ is determined by the model $X$
together with 
a single requirement from $p\cup p'$. 

\smallskip

It remains to show that $p\mapsto p\cut X$ witnesses that $p^X$ is
strongly generic.   We need to show that any condition $q\le p\cut X$ in $X$ is
compatible with $p$, and for this it is enough to show that 
$q$ is compatible with every 
requirement $\RR\in p$. 

In the case $\RR=\II_{\alpha,\lambda}\in p$ and $\lambda<\sup(X\cap\kappa)$
 there is an $X$-fence for $\II_{\alpha,\lambda}$ in $\cp {X}p\subseteq p\cut
X$, and any requirement in $X$ which is compatible with this
$X$-fence is compatible with $\II_{\alpha,\lambda}$.

Now consider  $\RR=\OOO\in p$.  If $\OOO\in X$ then $\OOO\in p\cap
X\subseteq p\cut X\subseteq q$, so it will be enough to show that if
$\OOO\notin X$ then $\OOO$ is 
compatible with every 
requirement which is a member of  $X$. 
Now if $\OOO$ is incompatible with any requirement in $X$ then $\OOO$ 
is incompatible with a 
requirement of the form $\II_{\alpha,\lambda}\in X$. 
In the case that $X$ is countable it then follows from the definition
of compatibility of 
$\OOO$ with $\CC_X$ that $\OOO\in X$, so we can assume that $X$ is 
uncountable.
It follows that $X\cap\kappa$ is transitive, and since
$\II_{\alpha,\lambda}\in X$ it follows that $A_{\alpha,\lambda}\subset
X$ and therefore $\gamma\in X$.   If $\card X=\kappa$ then
$\eta',\eta\in \kappa\subset X$, so $\OOO\in X$.    Otherwise we have
$\eta'<\lambda<\sup(X\cap\kappa)$, and therefore
$\eta<\sup(X\cap\kappa)$ since $\OOO$ is compatible with
$p^X=\sing{\II_{\sup(X),\sup(X\cap\kappa)}}$, so again $\OOO\in X$.

In the case $\RR=\CC_{M,a}$ with $M\cap X\in X$,  
clause~\ref{thm:CinM}(i) asserts that $\CC_{M,a}$ is compatible with any
requirement in $X$ which is compatible with $p\cut X$.
\smallskip{}  

It only remains to consider the case  $\RR=C_{M,a}$ when  $M\cap
X\notin X$.  In this case  $X$ must be  countable and $M\cap
X\cap H_{\kappa}=X\cap 
H_{\delta}$, where $\delta\deq\sup(M\cap X\cap\kappa)$.    If
$\II_{\alpha,\lambda}\in q$ and $\lambda<\delta$ then
$\lambda\in M$, and in this case $\II_{\alpha,\lambda}$ is its own
$M$-fence.  If  
$\sup(M\cap\kappa)>\lambda\geq\delta$ then the $M$-fence for $\CC_X$,
required for the compatibility of $\CC_{M,a}$ with $\CC_X$,  is a
$M$-fence for $\II_{\alpha,\lambda}$.
Hence any requirement $\II_{\alpha,\lambda}\in q$ 
is compatible with $\CC_{M,a}$.

Now we show that $\CC_{M,a}$ is compatible with any requirement
$\OOO\in q$ by showing that if $\OOO$ is a requirement in $X$ which is
compatible with $p\cut X$ but incompatible with some requirement
$\II_{\alpha,\lambda}\in M[a]$, then $\OOO\in M$.
First, we must have $\eta<\delta$, as otherwise $\OOO$ would be
incompatible with the $X$-fence for $\CC_{M,a}$, which is a member of
$\cp{X}p\subseteq p\cut X$.  Thus
$\sing{\eta',\eta}\subset M$.  Next, we have $\gamma\in
A_{\alpha+1,\lambda}$, where $\lambda\in M$ and $\alpha\in M\cup
a(\lambda)$.    If $\gamma<\alpha$ then, since $M\cap X\notin X$, 
lemma~\ref{thm:ASminBig} implies that $A_{\alpha,\lambda}\cap
X\subset M$, so $\gamma\in M$. 
If $\gamma=\alpha$, on the other hand, then $\gamma=\alpha\in
M[a]=M\cup a(\lambda)$, and 
$a(\lambda)$ is a set of nonzero limit ordinals while $\gamma$ is
either zero or a successor ordinal.   Thus it again follows that
$\gamma\in M$.  

Finally, suppose that $\CC_{N,b}\in q$.
Since $N\in X$ we have  $N\cap 
M\cap\kappa=N\cap(M\cap X)\cap\kappa=N\cap\delta\in X$, and since
$^{\omega}\delta\cap X\subset 
X\cap H_{\delta}=X\cap M\cap H_{\kappa}$ it follows that $N\cap 
M\cap\kappa\in M$.
Hence $\CC_{M,a}$ and $\CC_{N,b}$ satisfy clause~\ref{item:holds} of
definition~\ref{def:Pscompat}.  
We can obtain a $N$-fence for $M$ by taking the
minimal $N$-fences for the members of the minimal $X$-fence
for $M$, which is contained in $p\cut X$.
The $M$-fence for $X$ is also a
$M$-fence for $\CC_{N,b}$.   Hence $\CC_{M,a}$ is compatible with $\CC_{N,b}$.
\end{proof}

As a preliminary to the proof of lemma~\ref{thm:CinM}, we give a
structural characterization of the desired 
requirement $C_{M',a'}=\icut{C_{M,a}}X$.
Recall that $A^*_{M,a,\lambda}$ is the set containing $0$ together
with the successor ordinals from
$A_{M,a,\lambda}=\bigcup\set{A_{\alpha,\lambda}:\alpha\in M\cup 
    a(\lambda)}$. 
\begin{lemma}\label{thm:samefence}
  Suppose that $X$ and $\CC_{M,a}$ are as in lemma~\ref{thm:CinM},
  and that 
  $\CC_{M',a'}$ is a requirement such that $M'=M\cap X$, $a'\in X$, 
  and $A^*_{M',a',\gl}\cap X=A^*_{M,a,\gl}\cap X$ for all
  $\gl\in M'\cap\kappa$.     Then $\CC_{M',a'}$ satisfies the conclusion of
  lemma~\ref{thm:CinM}.
\end{lemma}
Following the proof of lemma~\ref{thm:samefence} we will 
construct such a requirement $C_{M',a'}$. 
\begin{proof}
  The proof breaks into 3 cases, numbered from 1 to 3, depending
  whether  the
  requirement $\RR$ has the form
  $\II_{\alpha,\lambda}$, $\OOO$ or $\CC_{N,b}$.
  Furthermore, each of these three cases has two subcases, which are
  labeled~(a) 
  and~(b) to correspond to the two clauses in the conclusion of
  lemma~\ref{thm:CinM}.

  Note that the hypothesis of lemma~\ref{thm:CinM} implies that
  $M'\cap \kappa$ is an initial 
  segment of $M\cap \kappa$.

  \CASE{1a}
  First suppose that $\RR=\II_{\ga,\gl}\in X$ and
  $\II_{\alpha,\lambda}$ is compatible 
  with $\CC_{M',a'}$.   We must show that $\II_{\alpha,\lambda}$ is
  compatible with $\CC_{M,a}$. 
  If $\sup(M\cap\kappa)>\gl\ge\sup(M\cap X\cap\kappa)$ then any
  $M$-fence for $X$ includes an 
  $M$-fence for $\II_{\alpha,\lambda}$, so we can assume that
  $\gl<\sup(M\cap X\cap\kappa)=\sup(M'\cap\kappa)$.   Set
  $\lambda'\deq\min(M\setminus\lambda)=\min(M'\setminus\lambda)$. 
  Now $A_{\alpha,\lambda}\cap
  M'=A_{\alpha,\lambda}\cap (M\cap X)=A_{\alpha,\lambda}\cap M$ since
  by lemma~\ref{thm:ASminBig} $M\cap X\in X$ implies that
  $A_{\alpha,\lambda}\cap M\subset X$.  Thus any $M'$ fence for
  $\II_{\alpha,\lambda}$ is also a $M$-fence for
  $\II_{\alpha,\lambda}$. 
  
  \CASE{1b}
  Now suppose that  $\RR=\II_{\ga,\gl}$ is compatible with 
  $\sing{\CC_{M,a}}\cup p^{X}$.   We will show that $\RR$ is compatible
  with $\CC_{M',a'}$.  
  This is immediate if $\lambda\ge\sup(M'\cap\kappa)$, so we
  assume that $\lambda<\sup(M'\cap\kappa)$.  Then
  $\lambda'\deq\min(M'\setminus\gl)=\min(M\setminus\lambda)$ and
  $A_{\alpha,\lambda}\cap M'\subset M$ since $M'\subset M$.
  Hence any $M$-fence for $\II_{\alpha,\lambda}$ is also an $M'$-fence
  for $\II_{\alpha,\lambda}$.

  \smallskip 

  \CASE{2a}
  Here we assume that 
   $\RR=\OOO\in X$ and 
    $\RR$ is
  compatible with $\sing{\CC_{M',a'}}\cup \cp{X}p$, and we 
  will show that  $R$ is compatible with $\CC_{M,a}$.   
  This is immediate unless there is some $\II_{\alpha,\lambda}\in
  M[a]$ which is incompatible with $\OOO$, and in this case we must
  have $\lambda<\sup(M\cap X\cap\kappa)$, or else $\OOO$ would be
  incompatible with the $X$-fence for $\CC_{M,a}$, which is contained
  in $\cp{X}p$.  Thus $\lambda\in
  X\cap M=M'$.    Furthermore $\gamma\in A^*_{M,a,\lambda}\cap
  X=A^*_{M',a',\lambda}\cap X$, so $\OOO$ is incompatible with a
  requirement in $M'[a']$ and hence is a member of $M'\subset M$.
  Thus $\OOO$ is compatible with $\CC_{M,a}$.
  \CASE{2b}
  Now suppose that $\RR=\OOO$ is compatible with
  $\sing{\CC_{M,a}}\cup p^X$.   
  Then $\OOO$ is compatible with $\CC_{M',a'}$ unless there is some
  $\II_{\alpha,\lambda}\in M'[a']$ which is incompatible with $\OOO$.
  Since $M'$ and $a'$ are in $X$ this implies that
  $\II_{\alpha,\lambda}\in X$, and since $\OOO$ is compatible with
  $p^{X}$ it follows that $\OOO\in X$.   Since $A^*_{M,a,\lambda}\cap
  X=A^*_{M',a',\lambda}\cap X$ it follows that $\OOO$ is incompatible
  with a requirement in $M[a]$, and hence $\OOO\in M$ since $\OOO$ is
  compatible with $\CC_{M,a}$.   Hence $\OOO\in M\cap X=M'$.

  \smallskip
 
  \CASE{3a}
  Suppose that $\RR=\CC_{N,b}\in X$ and $\RR$ is
  compatible with $\CC_{M',a'}\cup\cp{X}p$. 
  We need to show 
  that $\RR$ is   
  compatible with $\CC_{M,a}$.
  Since $\CC_{N,b}$ is compatible with $\CC_{M',a'}$, $M\cap
  N\cap\kappa=M'\cap N\cap\kappa$ is either a member of or an initial segment
  of $N$.
  Also, since $M'\cap\kappa$ is
  an initial segment of $M$ the set $M\cap N\cap\kappa$ is also either
  a member of or an initial segment of $M$ according as  it is
  a member or initial segment of $M'$.

  Let $x'$ be a $M'$-fence for $\CC_{N,b}$ and let $x$ be a $M$-fence
  for $\CC_{X}$.   Then $x'\cup x$ is a $M$-fence for
  $\CC_{N,b}$:  If $\OOO\in M$ is incompatible with some
  $\II_{\alpha,\lambda}\in N[b]$ then $\OOO$ is incompatible with some
  member of $x$ if $\lambda>\sup(M\cap X\cap\kappa)$, and otherwise $\OOO$ is
  incompatible with some member of $x'$.

  In the other direction, let $x'$ be an $N$-fence for
  $\CC_{M',a'}$, let $y$ be an $X$-fence for $\CC_{M,a}$ which is
  contained in 
  $\cp{X}p$, and let $x$ be the set of minimal $N$-fences for
  members of $y$.  
  Then $x\cup x'$ is a
  $N$-fence for $\CC_{M,a}$: 
  Let  $\OOO\in N$ be a 
  requirement which is incompatible with some requirement
  $\II_{\alpha,\lambda}\in M[a]$.
  If $\lambda\ge\sup(M\cap N\cap\kappa)$ then $\OOO$ must be
  incompatible with 
  some member of $x$.
  If $\lambda<\sup(M\cap N\cap\kappa)$ then
  $\lambda\in M'$ since $N\in X$ implies that $\sup(M\cap
  N\cap X)\leq\sup(M\cap X\cap\kappa)$, and 
  $M\cap X\cap \kappa=M\cap\sup(M\cap X\cap\kappa)$.   Thus 
  the fact 
  that  $A^*_{M,a,\lambda}\cap
  X=A^*_{M',a',\lambda}\cap X$ implies that $\OOO$ is incompatible with
  some member of $M'[a']$ and hence with some member of $x'$.

  \CASE{3b}
  Finally, suppose  $\RR=\CC_{N,b}$ is compatible with
  $\sing{\CC_{M,a}}\cup p^X$. 
  We need to show that $\RR$ is also compatible with $C_{M',a'}$.
  First, $N\cap M'\cap\kappa=N\cap (M\cap X)\cap\kappa$.   This is an
  initial segment of $N\cap M$; thus it is
  either a member or initial segment of $M'$ depending on whether
  $N\cap M\cap\kappa$ is a member or initial segment of $M$, and it
  is a a member or initial segment of $N$ depending on whether
  $N\cap M\cap\kappa$ is a member or initial segment of $N$. 
  
  If $x$ is any  $M$-fence for $C_{N,b}$ then
  $\set{\II_{\alpha,\lambda}\in x:\lambda\in M'}$ is  an $M'$-fence for
  $\CC_{N,b}$.  
  An
  $N$-fence for $\CC_{M',a'}$ can be obtained by taking the union of an
  $N$-fence for $\CC_X$ and an $N$-fence for $\CC_{M,a}$: If
  $\II_{\gamma,\xi}\in M'[a']$ then $\lambda\in M'=M\cap X$ and either
  $\gamma\in M'\subseteq M$ or $\gamma\in a'(\lambda)\subset X$.
\end{proof}

\begin{proof}[Proof of lemma~\ref{thm:CinM}]
  It remains to construct a pair $M',a'$ 
  satisfying the hypothesis of lemma~\ref{thm:samefence}.
  We already have $M'=M\cap X$.   In order to construct $a'$ we will
  define a sequence of proxies $a(i)$ and $b(i)$ by recursion on $i$,
  each of which  satisfies the following recursion hypotheses:
  \begin{enumerate}
  \item
    \begin{inparaenum}[(a)]
    \item \label{item:CinM1a}
      $\CC_{M',a(i)}$ is a requirement,
    \item\label{item:CinM1b}
      $\CC_{M,b(i)}$ satisfies
      Definition~\ref{def:requireM}(\ref{item:x2},\ref{item:x4}), and
    \item\label{item:CinM1c}
      for any $\nu\in M'\cap\kappa$ and any $\alpha\in b(i)(\nu)$,
      either $A_{\alpha,\nu}\cap M'$ is unbounded in $M'$ or
      $\sup(A_{\alpha,\nu}\cap X)\in M'$.
    \end{inparaenum}
  \item\label{item:CinM2} $a(i)\in X$.
  \item\label{item:CinM3}
    $A^*_{M',a(i)\cup b(i),\nu}\cap X=A^{*}_{M,a,\nu}\cap X$ for all
    $\nu<\sup(M'\cap\kappa)$.
  \item\label{item:CinM4}
    Set $d(i)=\set{\alpha:\exists\lambda\,(\alpha,\lambda)\in
      b(i)}$.  Then $d(i+1)\lessdot d(i)$ where 
    $\lessdot$ is  the ordering  of 
    $[\kappa^+]^{<\omega}$ defined by $d'\lessdot d$ if
    $\max(d'\dinter d)\in d$.   
  \end{enumerate}
  The ordering $\lessdot$ is a well order, so clause~3 implies that
  there is some $k<\omega$ such that $b(k)=\nothing$.   We will 
  set $a'=a(k)$.   Then $\CC_{M',a'}$ is a requirement by clause~(1) of
  the recursion hypothesis, it is in $X$ by clause~2, and it
  satisfies $A^{*}_{M',a',\lambda}\cap X=A_{M,a,\lambda}\cap X$ for
  $\lambda\in M'\cap \kappa$ by clause~3.  Hence $\CC_{M',a'}$
  satisfies the conclusion of lemma~\ref{thm:CinM}. 

  Note that clause~\ref{item:CinM1b} is a modification of
  clause~\ref{item:x3} of the Definition~\ref{def:requireM} of a
  requirement of the type $\CC_{M,a}$. 

  \CASE{$i=0$}
  The recursion starts  with $a(0)=\nothing$ and 
  \begin{equation*}
    b(0)= b
    \cup\set{(\alpha,0): \alpha\in M\setminus\sup(M')
      \And   X\cap\alpha\not\subset\sup(M\cap \alpha)
    }. 
  \end{equation*}
  The set $b(0)$ is finite, since by proposition~\ref{thm:limMcapN}
  there can be only finitely many $\alpha>\sup(M\cap X)$ in $M$ such
  that $X\cap\alpha\not\subset\sup(M\cap\alpha)$.

  Clause~2 of the recursion hypotheses is immediate and clause~4 does
  not apply, so we only need to verify clauses~1 and~3.

  Clause~\ref{item:CinM1a} is immediate since $a(0)=\emptyset$. 
  Since the requirement $\CC_{M,a}$
  satisfies \ref{def:requireM}(\ref{item:x2},\ref{item:x4})  and
  $b(0)(\lambda)\subseteq M\cup
  a(\lambda)$  for $\lambda\in M\cap\kappa$, clause~\ref{item:CinM1b}
  holds for $b(0)$.
  Finally, clause~\ref{item:CinM1b} follows from 
  Lemma~\ref{thm:shadow2}(\ref{item:shadow2a}).

  Now we verify clause~4 of the recursion hypothesis: 
  \begin{CLAIM}
    $A_{M',a(0)\cup b(0),\nu}=A_{M,a,\nu}$ for all $\nu\in M'\cap\kappa$.
  \end{CLAIM}
  \begin{proof}
    We have $a(0)=\nothing$, and it is clear that
    $A_{M',b(0),\nu}\subseteq 
    A_{M,a,\nu}$.  
    Since $a\subseteq b(0)$  it only remains to show that
    $A_{M,\emptyset,\nu}\cap X\subseteq 
    A_{M',b(0),\nu}$.  Suppose $\gamma\in X\cap A_{\alpha,\nu}$ where
    $\alpha\in M$ and $\nu\in M'\cap\kappa$.  
    If $\gamma<\sup(M\cap X)$ and $A_{\alpha,\nu}\cap X$ is bounded in
    $\sup(M\cap X)$ then clause~1 of
    lemma~\ref{thm:shadow2} implies that $\alpha'\deq\sup(X\cap
    A_{\alpha,\nu})\in M\cap X=M'$, and then $\gamma\in
    A_{\alpha'+1,\nu}\subseteq A_{M',\emptyset,\nu}$.
    If $\gamma<\sup(M\cap X)$ and $A_{\alpha,\nu}\cap X$ is unbounded in
    $X$ then $\gamma\in A_{\alpha',\nu}\subseteq A_{\alpha,\nu}$ for
    any $\alpha'\in A_{\alpha,\nu}\cap X\setminus\gamma$.   Thus we
    can assume that $\gamma\geq\sup(M\cap X)$.   Then 
    $(\alpha',0)\in b(0)$ where $\alpha'\deq\min(M\setminus\gamma)$,
    so $A_{\alpha',\nu}\subseteq A_{M',b(0),\nu}$, and it follows by 
    proposition~\ref{thm:min-is-enough} that
    $\gamma\in A_{\alpha,\nu}\cap\alpha'=A_{\alpha',\nu}\subseteq
    A_{M',b(0),\nu}$. 
  \end{proof}
  
  \CASE{$i+1$}
  Now assume that $a(i)$ and $b(i)$ have been defined, and
  $b(i)\not=\nothing$.   To define $a(i+1)$ and $b(i+1)$, 
  let  $(\alpha,\lambda)$ be the lexicographically least member of
  $b(i)$, and set $b'(i)=\set{(\alpha',\lambda')\in
    b(i):\alpha'>\alpha}$.  
  Note that, while there may be more than one
  ordinal $\lambda$ such that $(\alpha,\lambda)\in b(i)$, all but the
  least of these are redundant and may be discarded.   

  We begin with several special cases:  
  If $\lambda\ge\sup(M'\cap\kappa)$ then $(\alpha,\lambda)$ 
  can 
  be discarded since $\nu<\sup(M'\cap\kappa)$; in this case we set
  $a(i+1)=a(i)$ and $b(i+1)=b'(i)$. 
  If $\alpha\in X$ then we set
  $a(i+1)=a(i)\cup\sing{(\alpha,\min(X\setminus\lambda)}$ and
  $b(i+1)=b'(i)$.   If $\alpha\le \alpha'+\omega$ for some limit
  ordinal $\alpha'$ then we set $a(i+1)=a(i)$ and
  $b(i+1)=\sing{(\alpha',\lambda)}\cup b'(i)$.

  The recursion hypotheses are clear in each of these three cases.
  For the 
  remainder we can assume that $\alpha\notin X$ and that all members
  of $C_{\alpha}$ are limit ordinals.
  If $C_{\alpha}$ is bounded in $X$ then 
  set $\eta\deq\max(\lim(C_{\alpha})\cap\bar X)$ where $\bar X$ is the
  closure  $\bar X=X\cup\lim(X)$ of $X$.
  Otherwise, if $C_{\alpha}$ is cofinal in $X$, set
  $\eta=c_{\alpha,\sup(M'\cap \kappa)}$.  Then
  $A_{\alpha,\nu}=A_{\eta,\nu}$ for all $\nu<\sup(M'\cap\kappa)$, and
  $\eta\in X$ since $M'\in X$ and
  $\eta=c_{\alpha',\sup(M'\cap\kappa)}$ for any
  $\alpha'\in\lim(C_\alpha\cap X)\setminus\eta$. 
 
  \begin{CLAIM}
    If $\eta\leq\sup(M')$ then $A_{\alpha,\nu}\cap\eta\subseteq
    A_{M',\emptyset,\nu}$ for all $\nu\leq\sup(M'\cap\kappa)$.
  \end{CLAIM}
  \begin{proof}
    For any $\nu$ such that 
    $A_{\alpha,\nu}$ is cofinal in $M'$ we have
    $A_{\eta,\nu}=A_{\alpha,\nu}\cap\eta\subseteq 
    A_{\alpha,\nu}\cap\sup(M')=A_{\sup(M'),\nu}\subseteq
    A_{M',\emptyset,\nu}$ since $\lim(C_{\sup(M')})\cap M'$ is cofinal
    in $M'$.  
    If $\nu\in M'\cap\kappa$ and $A_{\alpha,\nu}$ is bounded in
    $\sup(M')$     then
    $\xi\deq\sup(A_{\alpha,\nu}\cap X)\in M'$  by clause~\ref{item:CinM1c} of
    the recursion hypothesis.   Then 
    $A_{\alpha,\nu}\cap\eta=A_{\xi,\nu}\cap\eta\subseteq
    A_{M',\emptyset,\nu}$. 
  \end{proof}

  Thus if $\eta\leq\sup(M')$ we can set $a(i+1)=a(i)$. 
  Otherwise set
  $\eta'=\min(X\setminus\eta)$ and $\lambda'=\min(X\setminus\lambda)$,
  and set $a(i+1)=a(i)\cup\sing{(\eta',\lambda')}$.  
  
  If
  $A_{\alpha,\nu}\cap X\subseteq A_{\eta,\nu}$ for all
  $\nu<\sup(M'\cap\kappa)$ then set $b(i+1)=b'(i)$.
  Otherwise let $(\gamma_j:j<m)$ enumerate the
  set $\set{\min(C_{\alpha}\setminus\xi):\eta< \xi\in X\cap
    A^{*}_{\alpha,\sup(M'\cap\kappa)}}$.   
  Note that $m$ is finite  since
  otherwise $\sup_{j<\omega}\gamma_j$ would 
  be in $\lim(C_{\alpha})\cap\bar X$.
  For each $j<m$ let $\chi_j$ be the least ordinal $\chi\in
  X\setminus\lambda$ such that 
  $\gamma_j\in A_{\alpha,\chi}$, and     
  set $b(i+1)=b'(i)\cup \set{(\gamma_j,\chi_j):j<m}$.

  \smallskip{}
  This completes the definition of $a(i+1)$ and $b(i+1)$.   
  Again, Clauses~2 and~4 of the recursion hypotheses are clear.
  Clause~\ref{item:CinM1a} is also immediate unless $a(i+1)\not=a(i)$,
  in which case
  we need to show that each clause of Definition~\ref{def:requireM}
  holds of $(\eta',\lambda')$ for $\CC_{M',a(i+1)}$.
  Clause~\ref{def:requireM}(\ref{item:x1}) is clear.
  For 
  clause~\ref{def:requireM}(\ref{item:x2}), note that if $\nu\in
  M'\cap(\kappa\setminus\lambda')$, $\gamma\in M'\cap\kappa$ and
  $\pi_{\eta'}(\gamma,\nu)$ is defined then 
  $\pi_{\eta'}(\gamma,\nu)\in X$ since $\sing{\nu,\gamma,\eta'}\subset X$.
  In addition
  $\pi_{\eta'}(\gamma,\nu)=\pi_{\eta}(\gamma,\nu)=
  \pi_{\alpha}(\gamma,\nu)\in M$, and hence
  $\pi_{\eta'}(\gamma,\nu)\in M\cap X=M'$.
  For clause~\ref{def:requireM}(\ref{item:x3}),
  $A_{\eta,\nu}$ is an initial segment of $A_{\alpha,\nu}$, and
  clause~\ref{item:CinM2} asserts that either $A_{\alpha,\nu}\cap M'$
  is unbounded in $M'$ or else $\sup(A_{\alpha,\nu}\cap X)\in M'$.
  Since $\eta>\sup(M')$, the first alternative implies that
  $A_{\eta',\nu}$ is cofinal in $M'$.   Since $\eta'$ and $\nu$ are in
  $X$, the second alternative implies that
  $\sup(A_{\eta',\nu}\cap X)=\sup(A_{\eta,\nu}\cap
  X)=\sup(A_{\alpha,\nu}\cap X)\in M'$.  
  For
  clause~\ref{def:requireM}(\ref{item:x4}), if $\nu\notin B_{\eta'}$
  and $\gamma$ is least such that 
  $\gamma\in A_{\eta',\nu}$ and $\nu\notin
  B_{\gamma}$ then $\gamma\in X$, and hence 
  $\gamma\in A_{\eta',\nu}=A_{\alpha,\nu}\cap\eta'$.   Thus $\gamma\in
  M$ by clause~\ref{item:CinM1b} of the recursion hypothesis, so
  $\gamma\in M\cap X=M'$. 

  To verify clauses~\ref{item:CinM1b} and~\ref{item:CinM1c}, 
  use the recursion hypothesis and the fact that
  $A_{\gamma_j,\nu}=A_{\alpha,\nu}\cap\gamma_j$ for
  $\nu\geq\chi_{i+1}$.

  It only remains  to verify clause~3 in the final case of the
  definition: 
  \begin{CLAIM}
    In the final case of  the definition of $a(i+1)$ and $b(i+1)$ we
    have 
    $A^*_{M',a(i)\cup b(i),\nu}\cap X=
    A^*_{M,a,\nu}\cap X$ for all $\nu\in M'\cap\kappa$.
  \end{CLAIM}
  \begin{proof}
     The change from $a(i)\cup b(i)$ to $a(i+1)\cup b(i+1)$
    consists of replacing the single pair $(\alpha,\lambda)\in b(i)$ with
    the finite set 
    $\set{(\gamma_j,\chi_j):j<m}\subseteq b(i+1)$,  together with 
    $(\eta',\min(X\setminus\lambda))\in a(i+1)$ if $\eta>\sup(M')$.
    If $\nu<\lambda$ then none of these contributes any members to
    either of the  
    sets  $A^*_{M',a(i)\cup 
      b(i),\nu}$ or $A^*_{M,a,\nu}$, so 
    it will be sufficient to verify the
    conclusion of the claim for $\nu\in
    M'\cap(\kappa\setminus\lambda)$.
    Since $M'\subseteq X$, this
    implies that $\nu\geq\min(X\setminus\lambda)$.
    Thus it will be sufficient to
    verify that 
    \begin{equation}\label{eq:f}
    A^*_{\alpha,\nu}\cap X=
      \left(A^*_{\eta',\nu}\cup\bigcup\set{A^*_{\gamma_j,\nu}:j<m
          \And\chi_j\le\nu} \right)   \cap  X
    \end{equation}
    for any $\nu>\lambda$ in  $M'$.

    Since $\eta\in \lim(C_{\alpha})\cap\lim(C_{\eta'})$,
    we have 
    $A_{\alpha,\nu}\cap\eta=A_{\eta,\nu}=A_{\eta',\nu}\cap\eta$ for every
    $\nu<\kappa$.
    Furthermore $A_{\gamma_j,\nu}\subseteq A_{\alpha,\nu}$ for all
    $\nu\ge\chi_j$.   Thus it 
    it will be sufficient to show that
    $A^*_{\alpha,\nu}\cap(X\setminus\eta)\subset
    \bigcup\set{A^{*}_{\gamma_j,\nu}:j<m\And\chi_j\le\gamma_j}$.

    Suppose that $\gamma\in A^{*}_{\alpha,\nu}\cap(X\setminus\eta)$.
    Then there is some $j<m$ such that
    $\gamma'_j<\gamma\le\gamma_{j}$, where
    $\gamma'_j=\sup(C_\alpha\cap \gamma_j)$.
    However $\gamma\notin C_{\alpha}$ since it is a successor ordinal,
    so $\gamma<\gamma_j$.   Furthermore $\gamma\in A_{\alpha,\nu}$
    implies that $\gamma_j\in A_{\alpha,\nu}$ and hence
    $\nu\ge\chi_j$, so $\gamma_j\in b(i+1)(\nu)$ and $\gamma\in
    A_{\alpha,\nu}\cap\gamma_j=A_{\gamma_j,\nu}\subseteq A_{b(i),\nu}$.
  \end{proof}
  This completes the proof of lemma~\ref{thm:CinM} and hence of the
  strong genericity lemma~\ref{thm:smc}.
\end{proof}

%================================================================

\subsection{Completion of the proof of theorem~\ref{thm:mainthm}}
\label{sec:finish}
We first verify 
that there are stationarily many  models satisfying the 
hypothesis of  lemma~\ref{thm:smc}:

\begin{lemma}\label{thm:simple}
  \begin{inparaenumi}
  \item The set of transitive \simple{} \model{s} $X\prec H_{\kappa^+}$ is
    stationary. 
  \item The set of countable \simple{} \model{s} $M\prec H_{\kappa^+}$ is
    stationary. 
  \item If $\kappa$ is $\kappa^+$-Mahlo then the set of \simple{} \model{s}
    $Y\prec H_{\kappa^+}$ with $Y\cap\kappa\in B_{\sup(Y)}$
    is stationary.
  \end{inparaenumi}
\end{lemma}
\begin{proof}
  For clause~(i), any transitive set $X\prec H_{\kappa^+}$ with 
  $\cof(\sup(X))=\kappa$ is a simple \model{}.

  For the remaining clauses, let $X$ be   any model as in the last
  paragraph and set $\tau=\sup(X)$.

  For clause~(ii), let $M$ be
  any  countable elementary substructure of 
  the structure $(X,C_{\sup(X)})$.
  Because $C_{\sup(X)}$ was included as a predicate, $M\cap
  C_{\sup(X)}$ is unbounded in $\delta\deq\sup(M)$ and hence
  $C_{\delta}=C_{\sup(X)}\cap\delta$.
  Finally, $c_{\delta,\xi}=c_{\sup(X),\xi}\in M$ for all $\xi\in
  M\cap\kappa$, 
  so $\otp(C_{\delta})=\sup(M\cap\kappa)$  and $\lim(C_{\delta})$ is
  cofinal in $M$.
  Thus $M$ is a \simple{} \model{}. 

  For clause~(iii), let $E$ be the closed and unbounded set of cardinals 
  $\lambda<\kappa$ such that there is a set
  $X_\lambda\prec (X,C_{\sup(X)})$ with $X_{\lambda}\cap
  H_{\kappa}=H_{\lambda}$. 
  As in the last paragraph the \model{s} $X_{\lambda}$ are \simple{}. 
  Set $\tau=\sup(X)$.  Since $\kappa$ is $\tau+1$-Mahlo there is a
  stationary set of $\lambda\in E\cap B_{\tau}$.  
  Pick $\lambda\in E\cap
  B_\tau$,  and set $\tau'=\sup(X_\lambda)$.
  Then $A_{\tau,\lambda}=A_{\tau',\lambda}$, so
  $f_{\tau}(\lambda)=f_{\tau'}(\lambda)$, and since $\lambda\in
  B_\tau$ it follows that $\lambda\in B_{\tau'}$ as
  well.    Thus the set $Y=X_{\lambda}$ satisfies clause~(iii).
\end{proof}
\begin{corollary}\label{thm:presat}
  The forcing $\Pkp$ has the $\kappa^+$-chain
  condition and is $\omega_1$-presaturated.

  If $\kappa$ is $\kappa^+$-Mahlo then
  $\Pkp$ is $\kappa$-presaturated.
\end{corollary}
\begin{proof}
  The proof is immediate from
  lemma~\ref{thm:smccc}, corollary~\ref{thm:smccc1}, lemma~\ref{thm:smc}
  and lemma~\ref{thm:simple}.
\end{proof}

\begin{corollary}\label{thm:cards}
  If $\kappa$ is $\kappa^+$-Mahlo and $G$ is a generic subset of
  $\Pkp$ then $\gw_1^{V[G]}=\gw_1^{V}$, 
  $\gw_{2}^{V[G]}=\gk$, and all cardinals larger than $\kappa$ are preserved.
\end{corollary}
\begin{proof}
  By corollary~\ref{thm:presat},  $\Pkp$ is $\omega_1$-presaturated,
  $\kappa$-presaturated and has the $\kappa^+$-chain condition.
  Hence these three cardinals, and all
  cardinals greater than $\kappa^+$, are
  preserved, and
  it only remains to show that all
  cardinals between $\gw_1$ and $\gk$ are collapsed.  This follows
  by the  proof of the corresponding lemma~\ref{thm:collapse} from
  section~\ref{sec:oneclub}, using $B_0$ in place of $B$, $\DD_0$ in
  place of $D$, and $\II_{0,\ga}$ instead of $\II_\ga$. 
\end{proof}

\begin{corollary}\label{thm:Igw2NS}
  If $\kappa$ is $\kappa^+$-Mahlo then every subset of $\Cof(\gw_1)$
  in $V[G]$ in $I[\gw_2]$ is nonstationary.
\end{corollary}
\begin{proof}
  Assume to the contrary that  $A\deq\seq{a_\xi:\xi<\gk}$ is a sequence of
  countable subsets of $\gk$ in $V[G]$ such that 
  the set $B(A)\cap\Cof(\omega_1)$ is 
  stationary, where $B(A)$ is the set defined in
  definition~\ref{def:Ikplus}.  
  Let $\dot A$ be a name for $A$.
  Fix a transitive \simple{} \model{} $X\prec
  (H_{\kappa^+},\dot A)$, so that 
  $\forces\dot A\in V[\dot G\cap X]$, and as in the proof of
  lemma~\ref{thm:simple} let $E$ be the set of
  $\lambda<\kappa$ such that 
  there is a model $X_{\lambda}\prec (X,\dot A,C_{\sup(X)})$ with
  $X_{\lambda}\cap H_{\kappa}=H_{\lambda}$.
  Then $E$ contains a  closed and unbounded subset of $\kappa$.
  Since $B(A)\cap\Cof(\omega_1)$ is stationary there is an ordinal
  $\lambda\in E\cap D_{\tau+1}\cap B(A)\cap\Cof(\omega_1)$,  where
  $\tau=\sup(X)$. 
  Then $X_{\lambda}$ is simple, $\II_{\sup(X_{\lambda}),\lambda}\in
  \bigcup G$, and 
  $\sing{\II_{\sup(X_{\lambda}),\lambda}}\forces\forall\nu<\lambda\;\dot
  a_{\nu}\in V[\dot G\cap X_{\lambda}]$.   Thus $a_{\nu}\in V[G\cap
  X_{\lambda}]$ for all $\nu<\lambda$.

  Now let $c\subset\lambda$ witness that $\lambda\in B(A)$.  Thus
  $\otp(c)=\omega_1$, $\bigcup c=\lambda$, and
  $c\cap\beta\in\set{a_\nu:\nu<\lambda}\subset V[G\cap X_{\lambda}]$ for all
  $\beta<\lambda$.    It follows by lemma~\ref{thm:approx} that $c\in
  V[G\cap X_{\lambda}]$.
  \smallskip{}

  We complete the proof by showing that this is impossible. 
  Let $\dot c$ be a $\Pkp\cap X_{\lambda}$-name for $c$.
  For a closed unbounded set $E'$ of cardinals $\lambda'<\lambda$
  there is a model $X'\prec 
  (X_{\lambda},\dot c)$ with $X'\cap H_{\kappa}=H_{\lambda'}$.   Since
  $\lambda\in B_{\tau+1}\cap D_{\tau+1}$ there is a cardinal
  $\lambda'\in  E'\cap 
  D_{\tau}$.   As in the previous argument, $X'$ is a \simple{} 
  \model{} and  $\sing{\II_{\tau,\lambda'}}\in G$ is a strongly
  $X'$-generic condition.   Since $\otp(c)=\omega_1\subset X'$, it
  follows that 
  $\sing{\II_{\tau,\lambda'}}\forces\dot c\subset X'$,
  contradicting the fact that $c$ is cofinal in $\lambda$.
\end{proof}
This completes the proof of theorem~\ref{thm:mainthm}.

\section{Discussion and questions}
Several related questions and ideas are discussed in the paper
\cite{mitchell.acus}, and we will only summarize some of them here.

The first problem is whether these techniques can be applied at larger
cardinals.   One easy answer to this problem is given for any
regular cardinal $\gk$ by substituting ``of size less than $\gk$'' for
``finite'' and using \model{}s of size $\gk$ instead of countable
\model{}s.   The resulting forcing adds closed unbounded subsets of
$\gk^{++}$ and demonstrates the consistency of the statement 
that every subset of $\Cof(\kappa^+)$ in $I[\gk^{++}]$ is
nonstationary.   

No such generalization is known for cardinals $\gk^+$ where
$\gk$ is a limit cardinal.   
This problem is of particular interest in the case when $\kappa$ is a
singular cardinal.  
Shelah has shown that if $\kappa$ is singular then $I[\kappa^+]$
includes a stationary 
subset of $\Cof(\lambda)$ for every regular $\lambda<\kappa$, but it
is open whether $\Cof(\lambda)\in I[\kappa^+]$ for any regular
$\lambda$ in the interval
$\omega_1<\lambda<\kappa$. 

Another natural question is whether the techniques of this paper can
be applied at multiple cardinals, giving a model in which, for
example, neither 
$I[\omega_2]\cap\Cof(\omega_1)$ nor $I[\omega_3]\cap\Cof(\omega_2)$
contain a nonstationary set.   This problem seems to be quite
difficult, and a useful test problem comes from considering the much
simpler argument, alluded to at the
send of section~\ref{sec:oneclub} and given in 
\cite{mitchell.acus}, which uses the techniques 
of this paper to give a model with no $\omega_2$-Aronszajn trees.
Can this construction by used to duplicate the results of
\cite{abraham83.at} by obtaining, from
a supercompact cardinal $\kappa$ and a weakly compact cardinal
$\lambda>\kappa$, a model with no $\omega_2$- or $\omega_3$-Aronszajn
trees?  
Two approaches to this problem have been attempted.
The first, an iteration of the basic method
analogous to Abraham's construction in \cite{abraham83.at}, initially seemed
quite promising; however the author has withdrawn previous claims to
have such a proof.
The second approach would 
operate simultaneously on both cardinals by using forcing with
finite conditions as in the present technique, but 
containing as requirements models of size less than $\kappa$ (that is,
less than $\omega_2$ in the generic extension) as well as 
countable models.    
This would give a structure analogous to a gap-2 morass.
The combinatorics of this approach are are quite complicated.

It seems that plausible that a solution for the problem concerning
$I[\omega_2]$ and $I[\omega_3]$ will require solutions to \emph{both}
approaches to the Aronszajn tree problem, with the second of the two
approaches being used to provide a structure at $\lambda$ like the 
$\square_{\kappa}$ sequence needed in this paper.

A third question is whether it is possible for $I[\gw_2]$ to be
$\gw_3$-generated, that is, that $I[\gw_2]$ cannot be normally
generated by any of its subsets of size less than $\omega_2$.
Note that the continuum hypothesis implies that $I[\gw_2]$ is trivial,
that is, $\gw_2\in I[\gw_2]$, and 
this paper presents a model in which $I[\gw_2]$ is generated by
$\Cof(\gw)$.
Either the model of 
section~\ref{sec:oneclub} or the original model \cite{mitchell.atit}
with no Aronszajn trees on $\gw_2$ give an example in which the
restriction of $I[\gw_2]$ to 
$\Cof(\gw_1)$ is generated by the single set
$\set{\nu<\gw_2:\cof^{V}(\nu)=\gw_1}$.    If $I[\gw_2]$ is generated
by fewer than $\gw_3$ many sets then it is generated by the diagonal
intersection of these sets, so if $2^{\gw_2}=\gw_3$ then
the only remaining possibility is that $I[\gw_2]$ requires $\gw_3$
generators. 

It is likely that it is possible to obtain such a model  by using
the techniques of this paper to add closed, unbounded subsets
$D_{\ga,\gl}\subset \gl\cap B^*_\ga$ for $\ga<\gk^+$ and $\gl\in
B_\ga$, with the sets $\set{D_{\ga,\gl}:\gl\in B_\ga}$ forming a
$\square_{\gk}$-like tree.   A witness that $A\deq B_{\ga+1}\setminus
B_\ga\in I[\gw_2]$ would then be given by $\set{D_{\ga,\gl}\cap
  C_\gl:\gl\in A}$, where $C_\gl$ is a closed, unbounded subset of
$\gl$ such that $C_\gl\cap B_\ga=\nothing$.

% ===========================================

%\newcommand\MR[1]{}
 \bibliographystyle{amsalpha}
\bibliography{logic}

\providecommand{\bysame}{\leavevmode\hbox to3em{\hrulefill}\thinspace}
\providecommand{\MR}{\relax\ifhmode\unskip\space\fi MR }
\providecommand{\MRhref}[2]{%
  \href{http://www.ams.org/mathscinet-getitem?mr=#1}{#2}
}
\providecommand{\href}[2]{#2}
\begin{thebibliography}{FMS88}

\bibitem[Abr83]{abraham83.at}
Uri Abraham, \emph{Aronszajn trees on {$\aleph \sb{2}$} and {$\aleph \sb{3}$}},
  Ann. Pure Appl. Logic \textbf{24} (1983), no.~3, 213--230. \MR{MR717829
  (85d:03100)}

\bibitem[AS83]{abraham-shelah.fcus}
Uri Abraham and Saharon Shelah, \emph{Forcing closed unbounded sets}, Journal
  of Symbolic Logic \textbf{48} (1983), no.~3, 643--657.

\bibitem[Bau84]{baumgartner.apfa}
James~E. Baumgartner, \emph{Applications of the proper forcing axiom}, Handbook
  of set-theoretic topology, North-Holland, Amsterdam, 1984, pp.~913--959.
  \MR{86g:03084}

\bibitem[FMS88]{formagshe.martins-max-i}
Matt Foreman, Menachem Magidor, and Saharon Shelah, \emph{Martin's maximum,
  saturated ideals, and non-regular ultrafilters, {I}}, Annals of Mathematics
  (2nd series) \textbf{127} (1988), 1--47.

\bibitem[Fri06]{friedman.ffc}
Sy-David Friedman, \emph{Forcing with finite conditions}, Set theory, Trends
  Math., Birkh\"auser, Basel, 2006, pp.~285--295. \MR{MR2267153}

\bibitem[Ham03]{hamkins:ggft}
Joel~David Hamkins, \emph{Extensions with the approximation and cover
  properties have no new large cardinals}, Fund. Math. \textbf{180} (2003),
  no.~3, 257--277, arXiv:math.LO/0307229. \MR{MR2063629}

\bibitem[Kos00]{koszmider:scuf}
Piotr Koszmider, \emph{On strong chains of uncountable functions}, Israel
  Journal of Mathematics \textbf{118} (2000), 289--315.

\bibitem[Mit04]{mitchell.wvsi}
William~J. Mitchell, \emph{A weak variation of {S}helah's {$I[\omega\sb 2]$}},
  J. Symbolic Logic \textbf{69} (2004), no.~1, 94--100. \MR{2 039 349}

\bibitem[Mit05]{mitchell.acus}
\bysame, \emph{Adding closed unbounded subsets of {$\omega\sb 2$} with finite
  forcing}, Notre Dame J. Formal Logic \textbf{46} (2005), no.~3, 357--371.
  \MR{MR2162106}

\bibitem[Mit06]{mitchell04:nhal}
\bysame, \emph{On the {H}amkins approximation property}, Ann. Pure Appl. Logic
  \textbf{144} (2006), no.~1-3, 126--129. \MR{MR2279659}

\bibitem[Mit73]{mitchell.atit}
William~J. Mitchell, \emph{Aronszajn trees and the independence of the transfer
  property}, Ann. Math. Logic \textbf{5} (1972/73), 21--46. \MR{47 \#1612}

\bibitem[She91]{shelah:rsss}
Saharon Shelah, \emph{Reflecting stationary sets and successors of singular
  cardinals}, Arch. Math. Logic \textbf{31} (1991), no.~1, 25--53, [Sh351].
  \MR{93h:03072}

\bibitem[Tod85]{todorcevic.dsct}
Stevo Todor{\v{c}}evi{\'c}, \emph{Directed sets and cofinal types}, Trans.
  Amer. Math. Soc. \textbf{290} (1985), no.~2, 711--723. \MR{MR792822
  (87a:03084)}

\bibitem[Zap96]{zapletal.ccf}
Jindrich Zapletal, \emph{Characterization of the club forcing}, Papers on
  general topology and applications (Gorham, ME, 1995), Ann. New York Acad.
  Sci., vol. 806, New York Acad. Sci., New York, 1996, pp.~476--484.
  \MR{97m:03084}

\end{thebibliography}
\end{document}